\documentclass[twoside,a4paper,11 pt,openright]{report}

\usepackage{geometry}
\usepackage[latin1]{inputenc}
\usepackage[T1]{fontenc}
\usepackage{amssymb}
\usepackage{euscript}
\usepackage{amsxtra}
\usepackage{ae}
\usepackage[all]{xy}
\usepackage{enumerate}
\usepackage[lite]{amsrefs}
\usepackage{amsthm}

\usepackage{beton}
\usepackage{eulervm}

\theoremstyle{plain}
\newtheorem{thm}{Theorem}
\newtheorem{lemma}[thm]{Lemma}
\newtheorem{lemmaanddef}[thm]{Lemma and Definition}
\newtheorem{prop}[thm]{Proposition}
\newtheorem{cor}[thm]{Corollary}

\theoremstyle{definition}
\newtheorem{defn}[thm]{Definition}
\theoremstyle{remark}
\newtheorem{rmk}[thm]{Remark}

\newtheorem{exples}[thm]{Examples}
\newcommand{\lra}{\leftrightarrows}
\newcommand{\C}{{\mathbb{C}}}
\newcommand{\R}{{\mathbb{R}}}
\newcommand{\Z}{{\mathbb{Z}}}

\newcommand{\N}{{\mathbb{N}}}
\newcommand{\rationals}{{\mathbb{Q}}}

\newcommand{\Ker}{\operatorname{\mathrm{Ker}}}

\newcommand{\Gr}{\mathrm{Gr}}
\newcommand{\classs}{\mathsf B}
\newcommand{\classb}{\mathsf E}
\newcommand{\partialdirac}{\mathrm{\partial}\!\!\!/}
\newcommand{\dirac}{\mathrm{D}\!\!\!\!/\thinspace}
\newcommand{\Dirac}{\mathcal{D}}
\newcommand{\inv}{^{-1}}

\newcommand{\pt}{\mathrm{pt}}
\newcommand{\id}{\mathrm{id}}
\newcommand{\sign}{(-1)^{q/2}}
\newcommand{\rk}{\mathrm{rk\ }}

\newcommand{\ch}{\operatorname{ch}}

\newcommand{\res}{\mathrm{res}}

\newcommand*{\q}{\mathsf{q}}
\newcommand*{\Q}{\mathsf{Q}}

\newcommand*{\PD}{\mathrm{PD}}
\newcommand*{\Tot}{\mathrm{Tot\ }}
\newcommand*{\Totpr}{{\mathrm{Tot}}^{\Pi}\mathrm{ }}
\newcommand*{\Totsum}{{\mathrm{Tot}}^{\oplus}\mathrm{ }}
\newcommand*{\KK}{\operatorname{KK}}

\newcommand*{\HH}{\operatorname{HH}_*}

\newcommand*{\HP}{\operatorname{HP}}
\newcommand*{\HPL}{\operatorname{HL}}
\newcommand*{\HA}{\operatorname{HA}}
\newcommand*{\CC}{\operatorname{\mathcal{CC}}}
\newcommand*{\CCP}{\operatorname{\mathcal{CP}}}
\newcommand*{\HC}{\operatorname{HC}}

\newcommand*{\RKKG}{R\!\operatorname{KK}_*^\Gamma}
\newcommand*{\K}{\operatorname{K}}
\newcommand*{\RK}{R\!\operatorname{K}}
\newcommand*{\Homol}{\operatorname{H}}
\newcommand*{\Homolc}{\operatorname{H}_\mathrm{c}}
\newcommand*{\HomoldR}{\operatorname{H}_\mathrm{dR}}
\newcommand*{\Homollf}{\operatorname{H}^\mathrm{lf}}

\newcommand*{\SL}{\operatorname{SL}}

\newcommand*{\SO}{\operatorname{SO}}
\newcommand*{\U}{\operatorname{U}}
\newcommand*{\Spin}{\operatorname{\mathsf{Spin}}}
\newcommand*{\Spinor}{\mathsf{Spin}}
\newcommand*{\Cliff}{\operatorname{\mathsf{Cliff}}}
\newcommand*{\SU}{\operatorname{SU}}
\newcommand*{\Fl}{\operatorname{Fl}}
\newcommand*{\End}{\mathop\mathrm{End}}
\newcommand*{\Morita}{\operatorname{\mathsf{Morita}}}

\newcommand*{\Aroof}{\operatorname{\widehat{A}}}

\newcommand*{\ev}{\mathrm{ev}}
\newcommand*{\odd}{\mathrm{odd}}

\newcommand*{\cpspace}{\C\mathbb{P}}
\newcommand*{\proj}{\mathbb{P}}
\newcommand*{\vare}{\xi}
\newcommand*{\diag}{\mathrm{diag}}

\newcommand*{\alTor}{\mathbb{T}}

\DeclareMathOperator{\Hom}{Hom}

\DeclareMathOperator{\Ind}{Ind}

\DeclareMathOperator{\AutTop}{Aut_{Top}}

\DeclareMathOperator{\Image}{Im}
\DeclareMathOperator{\Sec}{Sec}
\DeclareMathOperator{\rest}{res}

\newcommand*{\Hilm}{{\mathcal{E}}}

\newcommand*{\Hil}{{\mathcal{H}}}
\newcommand*{\EG}{\classb\Gamma}
\newcommand*{\Enerve}{E.\Gamma}
\newcommand*{\BG}{\classs\Gamma}

\newcommand*{\Cred}{C^*_{\mathrm{\it{r}}}}
\newcommand*{\Ko}{\mathcal{K}}

\newcommand*{\rcross}{\mathbin{\rtimes_{\mathrm{r}}}}
\newcommand*{\rlcross}{\mathbin{\ltimes_{\mathrm{r}}}}
\newcommand*{\protens}{\mathbin{\widehat{\otimes}}}
\newcommand*{\alcross}{{\rtimes}}

\newcommand*{\allcross}{{\ltimes}}

\newcommand*{\altens}{\mathbin{\otimes_{\mathrm{alg}}}}
\newcommand*{\cross}{\mathbin{\rtimes}}
\newcommand*{\defeq}{\mathrel{:=}}
\newcommand*{\btdR}{\mathcal{D}_\bullet}
\newcommand*{\CINF}{\mathop{\mathcal{C}^\infty}}
\newcommand*{\CCINF}{\mathop{\mathcal{C}_\mathrm{c}^\infty}}
\newcommand*{\Todd}{\mathop{\mathrm{Todd}}}
\newcommand*{\Sch}{\mathop{\mathfrak{S}}}
\newcommand*{\Defo}{\mathop{\mathfrak{D}}}
\newcommand*{\Cone}{\mathop{\mathcal{C}}}
\renewcommand*{\S}{\mathop{\mathcal{S}}}
\newcommand*{\smKo}{{\mathcal{K}^\infty}}
\newcommand*{\smA}{\mathcal{A}}
\newcommand*{\hqv}{\EG\times_\Gamma \partial(\EG)}

\newcommand*{\Pfaff}{\text{Pfaff}}
\newcommand*{\Tr}{\text{Trace}}
\newcommand*{\ad}{\text{ad}}
\newcommand*{\Ad}{\text{Ad}}

\usepackage{fancyhdr}
\usepackage{abstract}

\fancypagestyle{Selbergemacht}{
\lhead[\thepage]{Synopsis}
\rhead[]{\thepage}
}
\fancypagestyle{Selbergemacht2}{
\lhead[\thepage]{Introduction}
\rhead[]{\thepage}
}
\fancypagestyle{Selbergemacht4}{
\lhead[\thepage]{Bibliography}
\rhead[]{\thepage}
}
\fancypagestyle{Selbergemacht5}{
\lhead[\thepage]{A.4 A deformation version of the smooth Dirac-dual Dirac method}
\rhead[]{\thepage}
}
\begin{document}
\title{\fontfamily{cmr}\fontseries{bx}
\Huge
K-Theory and Cyclic Cohomology\\
of Crossed Products with Higher Rank Lattices\\
\vspace{4.1em}
\LARGE
Mathias Fuchs\\
\vspace{2.7em}
2007
}
\author{}
\date{}
\maketitle
\pagestyle{empty}
\thispagestyle{empty}
\mbox{}
\newpage
\section*{}
\begin{center}
Abstract
\end{center}
In \cite{CoCC}*{p. 255}, Connes shows that the class of unit in the crossed product of a discrete cocompact torsion-free subgroup of $SL(2,\R)$ with the algebra of continuous functions on $\mathbb{\R}P^1$ has torsion. However, he shows that this fails when $\R$ is replaced by $\C$.\\
In order to have a closer look at this phenomenon, we set out to generalize to the following case. We consider the $C^*$-algebra crossed product $C(\partial(G/K))\rcross\Gamma$ of a discrete torsion-free subroup $\Gamma$ of a connected semisimple Lie group $G$ with the algebra of continuous functions on the boundary of the symmetric space $G/K$ of the Lie group, where $K$ is a maximal compact subroup. We identify the map on the geometrical side of the Baum-Connes conjecture associated to the map on the crossed product induced by $\partial(G/K)\to\pt$ with a part of the topological Gysin sequence associated to the cosphere bundle of the tangent bundle of $\Gamma\backslash G/K$. It follows that the class of unit in the $\K$-theory of the crossed product is of finite order if and only if the rank of $G$ is equal to that of $K$, thus Connes's result. In case it is torsion we show that its order is the absolute value of the Euler characteristic of the discrete group if there is a $\Gamma$-equivariant Spin$^c$-structure on $G/K$. We make use of the injectivity of the assembly map.\\
In the second part, we carry out a new proof of the Kaplansky conjecture on the absence of nontrivial idempotents in the group rings $\C\Gamma$, based on the dual question to the one posed above: Does the class of the trace in cyclic cohomology lift to the crossed product $\C\Gamma$ with an algebra of functions on a boundary of the symmetric space? In fact, the idea is that this would allow to calculate the value of the pairing of the trace with the idempotents of $\C\Gamma$ on the crossed product rather than on $\C\Gamma$. This would be useful since the K-theory of the crossed product with the Furstenberg boundary is known thanks to the positive answer to the Baum-Connes conjecture for solvable groups, due to Kasparov (\cite{Kas}).\\
We therefore review topologically Nistor's computation of the homogeneous part of the periodic cyclic cohomology of crossed products $\Gamma\allcross A$ of torsion-free discrete groups $\Gamma$ with a complex $\Gamma$-algebra $A$. We use periodic cyclic cohomology associated to bornological algebras. Let $G$ be a complex semisimple Lie group and $B$ be a minimal parabolic subgroup of $G$. Applied to torsion-free discrete subgroups $\Gamma$ of $G$ and the algebra $A$ of smooth functions on the Furstenberg boundary $G/B$, the construction yields under a weak topological condition, together with the classical splitting principle, the surjectivity of the map $\HP^*\bigl(\Gamma\allcross\mathcal{C}^\infty(G/B)\bigr)_{<1>}\longrightarrow\HP^*(\C\Gamma)_{<1>}$ on the homogeneous parts of periodic cyclic cohomology of the map of smooth convolution groupoids associated to $G/B \to \pt$.\\
Let now $G$ be $\SL(n, \C)$ for arbitrary $n\ge 2$. Then there is, in analogy to the $C^*$-algebraic case, a dual Dirac morphism from the algebraic $\K$-theory group $\K_0(\mathcal{R}\otimes(\Gamma\allcross\mathcal{C}^\infty(G/B)))$ to the topological $\K$-theory group with compact supports $\K_c^{\dim B/T}(\Gamma\backslash G/T)$, $\mathcal{R}$ being the ring of smooth compact operators, and $T$ being a maximal torus of $B$. More important, an analogous statement holds for analytic cyclic homology. One sees that this $\K$-theory map has a dual map with respect to the index pairing which is an isomorphism.\\
In the proof, we show a strong version of equivariant Bott periodicity. In fact, we construct for any finite real representation space $V$ of $T$ and any algebra $A$ endowed with an action of $T$ an invertible element in $\HA(A\alcross T,(A\protens\CCINF(V))\alcross T)$, where $\HA$ denotes bivariant analytic cyclic homology (\cite{Meyerthesis}). By a deformation method, we deduce the existence of an isomorphism
$$
\HA_*(\Gamma\allcross\CINF(G/B))\cong\HA_{*+\dim B/T\!\! \mod\! 2}(\Gamma\allcross\CCINF(G/T)).
$$
As a consequence, the trace functional can be evaluated on the whole of $\K_0(\C\Gamma)$ by passing to the crossed products and can, with the help of standard index theorems, be seen to be integer if $\Gamma$ is cocompact. We thus obtain another proof of the Kaplansky conjecture for the considered groups than those of Lafforgue (\cite{Lafforgue}) and Puschnigg (\cite{Puschnigg}).
\pagenumbering{roman}
\setcounter{page}{2}
\tableofcontents
\pagestyle{fancy}
\chapter*{Synopsis}
{{\renewcommand{\leftmark}{Synopsis}
\addcontentsline{toc}{chapter}{Synopsis}
The goal of this synopsis is to provide a guide to the text in the sense that it wants to make clear its logical structure on the large scale. For instance, we recommend to consult this synopsis during the lecture of the technical theorems contained in the fourth chapter, in order to get information on what they will be needed for.\\
The difference to the abstract is that this synopsis is quite more detailed, and that it does not only want to give informations on the content but that it rather wants to make it easier to understand. The difference to the introduction is that the synopsis does not contain any motivation or aside remarks, and only very few recalls of well-known notions and facts.
\section*{K-theory}
\addcontentsline{toc}{section}{K-theory}
Let $\Gamma$ be a uniform lattice in a semisimple Lie group $G$, let $\EG=G/K, \BG=\Gamma\backslash G/K$, where $K$ is a maximal compact subgroup of $G$. All these spaces being endowed with their canonical Riemannian structures, we denote by $S\BG$ the cosphere bundle associated to to the tangent bundle of $\BG$, and $\partial\EG$ the visibility boundary of $\EG$. The homotopy quotient $\EG\times_\Gamma\partial\EG$ identifies with $S\BG$. Suppose that $\BG$ is Spin$^c$. Let $q=\dim\BG$. Then the Baum-Connes assembly map for the crossed product $\Gamma\rlcross C(\partial\EG)$ can be written under the form $\mu:\K^{*+q}(\EG\times_\Gamma \partial\EG)\to\K_*(\Gamma\rlcross C(\partial\EG))$. There is a diagram
$$
\xymatrix{
\K^q(\BG)\ar[r]^\mu\ar[d]^{\pi^*}&\K_0(\Cred\Gamma)\ar[d]\\
\K^q(S\BG)\ar[r]^-{\mu^{G/B}}&\K_0(\Gamma\rlcross C(\partial\EG))
}
$$
which by some general properties of the outer Kasparov product commutes.\\
The map $\pi^*$ that appears at the left hand side is rationally injective iff the Euler class $e(\BG)$ vanishes.\\
Motivated by the literature, one sets out to study the question if the class of unit in $\K_0(\Cred\Gamma)$ is torsion, in dependence on $\Gamma$ or, better yet, only on $G$. $[1]\in\K_0(\Cred\Gamma)$ resp. in $\K_0(\Gamma\rlcross C(\partial\EG))$ is the image of the fundamental class $[\BG]$ associated to the Spin$^c$ structure in $\K^q(\BG)$ under the assembly map. So $[1]$ is torsion iff $e\ne 0$.\\
This latter statement holds iff the Euler-Poincar\'e form $\xi$ on
$G/K$ does not vanish, which is a consequence of the Chern-Weil
theory of characteristic classes. The curvature tensor on $G/K$ which
determines $\xi$ is, up to sign, equal to that of the compact dual
$G'/K'$ of the symmetric space $G/K$, so $\xi\ne 0$ iff the Euler
characteristic $\chi(G'/K') \ne 0$. A classical theorem of Hopf and
Samelson of compact homogeneous spaces states that $\chi(G'/K')\ne 0$
iff those ranks of $G'$ and $K'$ differ. Furthermore, one has rank
$G'$=rank $G$ and rank $K'$= rank $K$, so $[1]$ is torsion if rank
$G$= rank $K$.
\section*{Cyclic Cohomology}
\addcontentsline{toc}{section}{Cyclic Cohomology}
In the same situation as above, one restricts to complex semisimple Lie groups $G$ (which implies that $G$=2 rank $K$), and one replaces the visibility boundary $\partial\EG$ by the Furstenberg boundary $G/B$ of $G/K$, and one asks the dual question to the one discussed above: Does the class of the trace $[\tau]\in\HP^0(\C\Gamma)$ lift to $\HP^0(\Gamma\allcross \CINF(G/B))$, using the algebraic crossed product? The goal in what follows is to develop, starting from a positive answer, a proof of the Kaplansky conjecture for $\C\Gamma$. One shows the existence of a commutative diagram analogous to the one for K-theory
$$
\xymatrix{
\Homol^\ev(\BG)\ar[r]&\HP^0(\C\Gamma)\\
\Homol^{\ev+\dim G/B}(\EG\times_\Gamma G/B)\ar[u]^{\pi_!}\ar[r]&\HP^0(\Gamma\allcross\CINF(G/B))\ar[u]
}
$$
but making use of quite different techniques, i. e. here one inspects the cyclic bicomplex. $\pi_!$ denotes integration over the fibre. Having put in elementary algebraic topology, one sees that $\pi_!$ is surjective iff $\pi^*:\Homol^*(\BG)\to\Homol^*(\EG\times_\Gamma G/B)$ is. By a reformulation of the classical splitting principle, applied to the vector bundle $\EG\times_\Gamma \C^n\to\BG$ if $G=\SL(n,\C)$, one sees that this latter statement always holds. As a consequence, the trace class always lifts. In the more general case where $G\ne\SL(n,\C)$, the splitting principle applies whenever a crucial property on the topology of $G$ can be verified.\\
We want to calculate not only the homogeneous part of $\HP^*(\Gamma\allcross\CINF(G/B))$ but also analytic cyclic homology of this algebra because only the latter one behaves well under passage to $C^*$-algebraic completions. This is necessary because only the completion of the crossed product has an accessible K-theory, allowing to approach the Kaplansky conjecture. In fact, $\Gamma\rlcross C(G/B)$ is $\KK-$-equivalent to $C(\Gamma\backslash G/T)$, $T$ being a maximal torus of $B$, because $B$ is solvable, so its Kasparov $\gamma$-element equals 1. More precisely, we recall that $\K_*(\Gamma\rlcross C(G/B))$ is calculated as follows by the Dirac-dual Dirac method of Kasparov.\\
Let $T$ be the canonical maximal torus of $B$
$$
T=\{(z_i)\in\U(1)^{\times n}\ |\ \Pi z_i=1\}\subset\SL(N,\C).
$$
Let $\beta_B$ be the dual Dirac element for $B$, $\beta_B\in\KK^B_{\dim B-\dim T}(\C,C_0(B/T))$, invertible following Kasparov. By this one concludes that
\begin{align*}
j_{\Gamma\times B}(\Ind_{B\uparrow G}(\beta_B))\in\thinspace&\KK_{\dim B-\dim T}(\Gamma\rlcross C(G)\rtimes B,\Gamma\ltimes C_0(G)\otimes C_0(B/T)\rtimes B)\\
=\thinspace&\KK_{\dim \EG}(\Gamma\rlcross C(G/B),C(\EG\times_\Gamma G/B))
 \end{align*}
is invertible as well. One wants to mimic that in the case of smooth algebras. The technical theorem which allows to get off the ground is the following equivariant Bott periodicity:
\begin{itemize}
\item
Let $T\to\SO(V)$ be a finite orthogonal representation of the torus $T=\mathbb{S}^1\times\dots\times\mathbb{S}^1$, $A$ a bornological algebra with a smooth action of $T$. Then there exist an invertible element in the bivariant analytic cyclic homology $\HA_{\dim V}\bigl((A\protens \CCINF(V))\alcross T,A\alcross T\bigr)$.
\item\label{Spin}
One wants to replace the group $T$ by any Lie group such that the representation on $A\protens\CCINF(V)$ comes from a representation on the factors, and such that the one on $\CCINF(V)$ factorises via an orthogonal representation of a torus $T$ on $V$.
\end{itemize}
In other words, we can merely construct an analogue of the descent homomorphism of the invertible in $\KK^T_{\dim V}(C_0(V),\C)$, rather than an analogue to an invertible in $\KK^B_{\dim V}(C_0(V),\C)$.\\
Let us look at the action of $B$ on $B/T$. There is a smooth path of group homomorphisms $(\alpha_t:B\to B)_{t\in[0,1]}$ such that $\alpha_0=(B\to T\to B),\alpha_{0,5}=(B\to TA\to B),\alpha_1=\id_B$. Consequently, there are smooth algebras $\Gamma\allcross\CCINF(G\times B/T)\alcross_{\rho\times(\zeta\circ\alpha_t)}B$ where $\rho:B\to\End(\CCINF(G))$ and $\zeta:B \to \End(\CCINF(B/T))$ are the standard representations. Moreover, there is a corresponding algebra of sections of the algebra bundle over $[0,0.5]$ resp. $[0.5,1]$. Denote it by $S_{[0,0.5]}$ resp. $S_{[0.5,1]}$. (Attention: The homotopy equivalence of the spaces $T$ and $B$ and the homotopy of actions $B\circlearrowleft B/T$ and $B\to T\circlearrowleft B/T$ \emph{does not} yield a homotopy of the crossed product algebras.) One can with the help of a deformation \emph{d\'eformation}, construct a bivariant element in $\HA_{\dim B-\dim T}(\Gamma\allcross\CINF(G/B),\Gamma\allcross\CCINF(G/T))$, and show that the latter has the Dirac element as inverse. In fact, there is a commutative diagram (the one appearing in the proof of \ref{thmwithdiagram}) from which one concludes that
\begin{align*}
\alpha:\HA_*(\CCINF(\Gamma\backslash G/T))=\thinspace&\HA_*(\Gamma\allcross\CCINF(G\times B/T)\alcross B)\\
\longrightarrow\thinspace&\HA_*(\Gamma\allcross\CCINF(G)\alcross B)\\
=\thinspace&(\Gamma\allcross\CINF(G/B))
\end{align*}
is surjective. The injectivity follows from a comparison with the case of $C^*$-algebras, where there is an isomorphism with \emph{local cyclic homology} due to Michael Puschnigg $\HA_*(A)\cong\HPL_*(A)$ as well as an isomorphism $\HPL_*(A)\cong\K_*(A)\otimes\C$ whenever $A$ is $\KK$-equivalent to a commutative $C^*$-algebra, as $\Gamma\rlcross C(G/B)$ is, thanks to the existence of a bivariant Chern-Connes character $\KK\to\HPL$.\\
Let us give some indications on how the announced Bott periodicity above is shown. A \emph{quasi-homomorphism} between two algebras $A$ and $B$ is an algebra $C$ which contains $B$ as an ideal, together with two homomorphisms $f:A\to C, f':A\to C$ such that $\Image(f-f')\subset B$. One writes $(f,f'): A\to C\rhd B$. Making use of half-exactness of $\HA$, one finds that a quasi-homomorphism $(f,f')$ induces a bivariant element in $\HA_0(A,B)$. One then shows:
\begin{itemize}
\item
The dual-Dirac element
$$
\beta_B\in\KK^B(\C,C_0(B/T))
$$
for $B$ can be realized by a $T$-equivariant quasi-homomorphism between the smooth subalgebras
 $$
 \C\rightrightarrows \CINF(B/T)\otimes M_2\rhd\CCINF(B/T)\otimes M_2.
 $$
\item
Its inverse in $\KK^T$, the Dirac element $\alpha_B$, can be realized as a $B$-equivariant quasi-homomorphism between the smooth subalgebras.
$$
 \CCINF(B/T)\rightrightarrows \mathcal{L}(L^2 B/T)\otimes M_2\rhd \ell^p(L^2 B/T)\otimes M_2.
$$
\end{itemize}
It is then shown that the diagram
$$
\xymatrix{
\q^2\C\ar[rd]\ar[r]^{\q\beta}&\q\CCINF(\R^2)\otimes M_2\ar[r]^\alpha&\ell^p\otimes M_4\\
&\C\ar[ru]
}
$$
commutes up to smooth homotopy.\\
Moreover, it is necessary to understand the effect of $\alpha$, and thereby that of $\beta$, on $\HP(\Gamma\allcross\allowbreak\CINF(G/B))$. This is made possible due to the existence of the spectral sequence, which allows to reduce the study of $\HP(\Gamma\allcross A)$ \`to that of $\HP(A)$. So the Chern-Connes character of the Thom isomorphism comes in there, and this is well known from the classical proof of the Atiyah-Singer index theorem.\\
Having calculated $\HA_*(\Gamma\allcross G/B)$, together with some lemmas which assure that the various constructions are compatible with $\HP$ and $\KK$, one looks at the diagram \ref{Schlussdiagramm}, which involves the index theorem of Connes-Moscovici and makes sure that the idempotent $e\in\C\Gamma$ satisfies, with Atiyah's $L^2$-index theorem, $\tau(e)=\tau(\mu(\Dirac))=\Ind\Dirac$ for some differential operator $\Dirac$ on $\BG$. So $e\in\{0,1\}$.
\thispagestyle{Selbergemacht}
}
}
\clearpage{\pagestyle{empty}\cleardoublepage}
\pagenumbering{arabic}
\chapter{Introduction}
{\pagestyle{Selbergemacht2}
{
\section*{General Part}
\addcontentsline{toc}{section}{General Part}}
The $C^*$-algebras associated to Lie groups or to their discrete subgroups provide a rich source of examples for $\K$-theory of $C^*$-algebras. One approach to their study is to replace $C^*\Gamma$ by the crossed product algebra $C(M)\cross\Gamma$, where $M$ is a compact $\Gamma$-space. This ``interpolates'' between the cases of a free and proper action, where $C(M)\cross\Gamma$ is strongly Morita equivalent to $C(M/\Gamma)$ and the action on a point, where $C(\pt)\cross\Gamma=C^*\Gamma$.\\
In the first part of this work, we will consider the following situations: Fix some connected  semisimple Lie group $G$ and a discrete not necessarily cocompact but torsion-free subgroup $\Gamma$. Let $K$ be a maximal compact subgroup of $G$, and consider the homogeneous space $G/K$. This is a space with non-positive sectional curvature homeomorphic to some $\R^n$ on which $\Gamma$ acts freely and properly. So $G/K$ is a model for $\EG$, and $\Gamma\backslash G/K$ is a model for $\BG$. In what follows, we shall always mean by $\EG$ and $\BG$ these particular realizations with their smooth structure. We have the classical facts on existence and uniqueness of geodesics, and the boundary $\partial(G/K)$ of $G/K$ is defined as follows: Any two geodesic half-rays are said to be equivalent iff their distance is bounded for any $t$. The quotient space can be identified with a fiber of the sphere bundle associated to the tangent bundle of $G/K$ by assigning to a geodesic $\gamma$ its direction $\frac{d}{dt}\gamma(t)|_{t\downarrow 0}$ in the origin, thus $\partial (G/K)$ identifies with the $n-1$-sphere $\mathbb{S}^{n-1}$.\\
All this is a generalization of the classical example, where $G=PSL(2,\R)$, $K=PSO(2)$, so $G/K$ is the upper half (Poincar\'e) plane and the quotient $\Gamma\backslash G/K$ is a Riemann surface of genus $g \ge 2$. So if $\Gamma$ is cocompact, it equals the fundamental group
$$
\pi_1(\BG)=<a_1,b_1,\dots,a_g,b_g|a_1b_1a_1^{-1}b_1^{-1}\dots a_gb_ga_g^{-1}b_g^{-1}>
$$
of the Riemann surface. The boundary of the upper half plane is the one-point-compactification of the real line, i. e. the circle, and the action of $G$ on the boundary coincides with the natural one on $\R P^1=\mathbb{S}^1$. All this has a complex analog $G=PSL(2,\C)$, $K=PSU(2)$, where again the natural action of $PSL(2,\C)$ on $\C P^1=\mathbb{S}^2$ is that on the boundary of the symmetric space.\\
These actions of $\Gamma$ on the boundary are far from being proper, and, in these examples, the orbits are in fact dense in $\partial (G/K)$. So the topological quotient is degenerate, and, according to Connes' philosophy (\cite{Co2}), the crossed product is studied instead.\\
The algebras $C(M)\cross\Gamma$ can be introduced from other points of view as well:\\
Connes has shown how to associate to a foliation $\mathcal{F}$ of a manifold $V$a non-commutative $C^*$-algebra $C^*(V,\mathcal{F})$. If we take $\mathcal{F}$ to be the foliation of the unit sphere bundle $S\EG$ of $\EG=G/K$ defined by the equivalence relation introduced above on the geodesic half-rays, then it passes to $S\BG$ because $G$ acts by isometries, and then the equivalence classes make up a foliation $\mathcal{F}$ of the total space of $S\BG$. The $C^*$-algebra $C^*(V,\mathcal{F})$ will then be Morita equivalent to $C(M)\rcross\Gamma$.\\
As a third point of view, we can associate to $\Gamma$, in case it is word-hyperbolic, its Gromov boundary $\partial\Gamma$ and then we will have $\partial(\EG)=\partial\Gamma=\mathbb{S}^{n-1}$, and the actions coincide.\\
We will make use of the Baum-Connes-assembly map (\cite{Co3})
$$\mu: \RKKG(C_0(\EG),C(M))\longrightarrow \K_*(C(M)\rcross\Gamma)$$
whose injectivity for the cases we consider has been established by Kasparov (\cite{Kas}).\\
In \cite{CoCC}, Connes proves the following fact: If $G=PSL(2,\R)$ and $\Gamma$ is a torsion-free uniform lattice in $G$, then the class of unit $[1]\in \K_0(C(M)\cross \Gamma)$ is a torsion element. On the other hand, in the complex analog, the class of unit is of infinite order. The main purpose of this paper is to ``explain'' this phenomenon by generalizing it to arbitrary uniform lattices in connected semi-simple groups. Note that the phenomenon of the class of unit being of finite order is neither shared by commutative $C^*$-algebras nor by $C^*$-algebras of groups. If the group is discrete, for example, there is the canonical trace $\tau: \Cred\Gamma\to\C $ passing to $K_0$ which satisfies $\tau([1])=1$. The idea is to use a reduction of the problem to a geometric situation by means of the assembly map. In \cite{AnDe}, C. Anantharaman-Delaroche computes the order of the class of unit in the example $G=PSL(2,\R)$ and obtains the Euler characteristic of the group. We will show that this holds in full generality, provided $\EG$ has a $\Gamma$-equivariant Spin$^c$-structure.\\
One can also consider actions on trees, in which case Robertson (\cite{Ro}) obtains the analog result on the order of the class of unit again.\\
Our main result is the equivalence of the following statements:
\begin{enumerate}
\item[-] The class of unit is torsion.
\item[-] The Euler characteristic of the group is nonzero.
\end{enumerate}
This result links up with classical facts to the further equivalence to the statements:
\begin{enumerate}
\item[-] rk $\mathfrak{g}=$ rk $\mathfrak{k}$, i. e. there is a Cartan subalgebra of $\mathfrak{g}$ contained in $\mathfrak{k}$.
\item[-] There is a discrete series representation of $G$.
\end{enumerate}
Roughly speaking, the strategy of the proof is as follows: First we show that the class of unit is the image under the assembly map of the K-theoretic fundamental class $[\BG]$ of $\K^*(\BG)$, using the commutativity of Kasparov's outer product. The left hand side of the assembly map is the K-theory of the unit sphere bundle of the Riemannian manifold $\BG$, and the class in question is $\pi^*(\BG)$, where $\pi:S\BG\to\BG$ is the projection. The question if this class is of finite order or not is preserved by the Chern character. This allows to reduce the problem to a look at the Gysin sequence of the tangent bundle of $\BG$. The Gysin sequence explicitly involves the Euler class, and in view of the Gauss-Bonnet-theorem the equality o the first two statements is proved. Next, Hirzebruch's proportionality principle states that $\chi(\Gamma):=\chi(\BG)$ is nonzero if and only if rk $G$ = rk $K$, thus our main result.\\
In order to get information on the exact order of the class of unit, we will establish the Gysin sequence at the level of $\K$-theory. It involves the K-theoretic Euler class, and we show that that one corresponds under the Chern character to the ordinary Euler class. So the order of the class of unit is $\chi(\Gamma)$.\\
We note that for the cases in which the rank of $G$ is one, one can argue as follows: The boundary is a homogeneous space $P\backslash G$, where $P$ is a minimal parabolic subgroup, and there is the isomorphism given by Morita equivalence $\K_*(C(P\backslash G)\rcross\Gamma)\cong P\ltimes_\mathrm{r} C(G/\Gamma)$, and the latter group is isomorphic to $\K^*(T\backslash G/\Gamma)$ as we shall examine thoroughly, and also for higher ranks, in the second part.
Also, we mention  that the considerations made in this paper illustrate a part of the paper of Emerson \cite{Emersonthesis} on non-commutative Poincar\'e duality, and have a certain intersection with \cite{EmersonMeyer}.\\
We proceed as follows: In the first section of the second chapter, we give the setup and definitions we talk about. In the second we give a definition of Kasparov's $\alpha$ element for the group $\Gamma$ in the form in which we shall need it. In the third section we prove our main theorem. In the fourth section, we calculate the order of the class of unit.\\
The second part of this work is concerned with the $\K$-theory and cyclic cohomology of group rings of lattices in semisimple Lie groups. In particular, an approach to the idempotent theorem will be proposed which makes crucial use of the Furstenberg boundary of the globally symmetric Riemannian space $G/K$.\\
The content of the theorem is that the group algebra $\C\Gamma$ contains no idempotents other than $0$ and $1$. It has been a long standing conjecture, and for the groups considered in this paper it has been proved by Lafforgue (\cite{Lafforgue}), using Banach $\KK$-theory, and by Puschnigg (\cite{Puschnigg}), who computes the Chern-Connes character of the $\gamma$-element in local cyclic cohomology, which both involve technically quite advanced theory. The use of the Furstenberg boundary is supposed to provide a simpler and more geometric approach.\\
The theorem is weaker than the analogous statement for the $C^*$-algebraic completion $C^*_\mathrm{r}\Gamma$ where it is known under the name Kadison-Kaplansky conjecture and still unsolved for the class groups considered in this paper. The Kadison-Kaplansky conjecture is an easy corollary of the Baum-Connes-conjecture calculating $\K_0(C^*_\mathrm{r}\Gamma)$. However, the latter is unproved even for $\SL(3,\Z)$ or its principal congruence subgroups, which are torsion-free. Our approach to the theorem, being index-theoretic, is basically the same as the one via the Baum-Connes conjecture, and is the most popular since many years, based on the \emph{canonical trace} $\tau:\C\Gamma\to \C, \sum a_\gamma \gamma\mapsto a_1$ on $\C\Gamma$. Its usefulness for proving the idempotent theorem is given as follows: The theorem holds if the functional is integer-valued on all idempotents of the group ring. To prove this, one first needs to make use of the involution $\sum a_\gamma \gamma\mapsto\sum a_\gamma^*\gamma^{-1}$ on the group ring. Trivially the trace is then positive in the sense that $\tau(x^*x)\ge 0$ for all $x\in\C\Gamma$. Moreover, it is faithful in the sense that equality holds only if $x=0$. Let $e$ be a  projection. Then by the positivity $\tau(e)\ge 0$ and $1-\tau(e)=\tau(1-e)\ge 0$, so $\tau(e)\in\{0, 1\}$ if the trace is integer-valued on idempotents, and then the faithfulness implies $e\in\{0, 1\}$. If $e$ is a mere idempotent, not a projection, then one argues that $e$ is homotopic in the reduced $C^*$-algebra to a projection, where the same argument applies. The homotopy does not affect the value of the trace. Now the problem can be translated into $\K-$theory, because by the definition of algebraic $K_0$, $\tau$ passes to a functional $\K_0(\C\Gamma)\to\C$.\\
So the theorem holds if the trace functional on $\K_0(\C\Gamma)$ is integer.\\
Up to here, there has not yet been any mention of differential geometry or analysis whatsoever, up to the fact that we intend to consider lattices $\Gamma$ in Lie groups. However, one of the most fruitful sources of integrality in mathematics is index theory, and in fact index theory does the clue here, too. It comes in as follows. There is one canonical way of constructing $\K-$theory elements (and in many cases the only one known), namely via the assembly map
$$
\mu:\RK_0(\BG)\to\K_0(\mathcal{R}\Gamma)
$$
where the left hand side is $\lim_{X\subset\BG} \K_0(X)$, the $\K$-homology with compact support of $\BG$. The limit is taken over all compact subsets of $\BG$. Normally, the target of $\mu$ is the $\K$-theory of the reduced $C^*$-algebra, but it factors, in fact, through $\K_0(\mathcal{R}\Gamma)$ (\cites{ConnesMosc, Thom}). Here, $\mathcal{R}$ is supposed to mean the algebra of smooth compact operators. There is a map $\C\Gamma\to\mathcal{R}\Gamma$, and $\tau:\K_0(\C\Gamma)\to\R$ factorizes via $\K_0(\mathcal{R}\Gamma)$, given that $\mathcal{R}$ has a canonical trace, too. Now there is Atiyah's $L^2$-index-theorem which states that the trace of the assembly map of a class represented by an operator $\Dirac$ on $\BG$ is its index: $\tau(\mu([\Hil,\rho,\Dirac]))=\Ind \Dirac$.\\
Here, bounded Fredholm operators $\Dirac:\Hil^+\to\Hil^-$ which satisfy $\rho(f)\Dirac-\Dirac \rho(f)\in\Ko(\Hil^+,\Hil^-)$ for any $f\in C(X)$ are cycles for the even $\K$-homology of a compact space $X$, $\rho$ being an even representation $C(X)\to\mathcal{L}(\Hil)$, and the index of such a $\Dirac$ is the ordinary Fredholm index $\dim\Ker \Dirac-\dim\Ker (\Dirac^*)$.\\
 Therefore the idempotent theorem would follow from the Baum-Connes conjecture, or at least from the surjectivity of $\mu$, which however is an open problem since many years.\\
The idea underlying this paper is to consider a unital coefficient $\Gamma$-algebra $A$ with the following two crucial properties:
\begin{itemize}
\item
The trace is in the image of the restriction map on cyclic cohomology along the natural embedding $\C\Gamma\to \Gamma\allcross A$.
\item
The analytic cyclic homology of $\Gamma\allcross A$ as well as the $\K$-theory of a suitable $C^*$-completion of $\Gamma\allcross A$ are computable, and are equal, after tensorization with $\C$, to what the Baum-Connes conjecture predicts for K-theory.
\end{itemize}
Here and in the entire sequel we shall not distinguish algebraic or $C^*$-algebraic crossed products by notation,  in order to avoid overloaded formulae. We hope that the correct interpretation will always be clear from the context, whenever no explicit explanation is provided. Given these properties, we could then use the simple formula $<x,\tau>=<x,i^*\tau'>=<i_*x,\tau'>$ to calculate the image of the pairing on all $\K$-theory classes of $\mathcal{R}\Gamma$ by evaluating it on all $K$-theory classes of the crossed product that can come from $\K_0(\mathcal{R}\Gamma)$. The right argument will be slightly more complicated, however.\\
It turns out that $A=\CINF(G/B)$ does the job, where $B$ is a minimal parabolic subgroup of $G$, and $G/B$ is the so-called Furstenberg boundary. In case $G=\SL(n.\C)$, $B$ is the subgroup of upper triangular matrices with determinant $1$. The proof of the first property relies on an inspection of Nistor's classical calculation of the homogeneous part of cyclic cohomology of crossed products $\Gamma\allcross A$. For it can be stated as the existence of an isomorphism between group hypercohomology with coefficients in the cyclic cohomology of $A$ and the cyclic cohomology of $\Gamma\allcross A$. If $A=\CCINF(M)$ then under Poincare duality for $M$ the spectral sequence for hypercohomology becomes the Leray-Serre spectral sequence for the fibration $\EG\times_\Gamma M\to\BG$. The desired property then becomes equivalent to the requirement that the map from the cohomology of the base to that of the total space do not annihilate the fundamental class, in other words that the spectral sequence be degenerate. The latter condition is exactly assured by the splitting principle if the maximal compact subgroup $T$ of $B$ is also a maximal abelian subgroup of a maximal compact subgroup $K$ of $G$. This condition, in turn, is satisfied, for example, in the standard example $G=\SL(n,\C)$. In fact, the subgroup of diagonal matrices with entries in $\U(1)$ such that the determinant is one is a maximal compact abelian subgroup for $B$ but also one for $K=\SU(n,\C)$.\\
The second condition on $A$ relies on the fact that $B$ is a \emph{solvable} Lie group. So Kasparov's Dirac-dual Dirac machinery allows to construct an invertible Thom element in $\KK_{\dim B/T}^\Gamma\bigl(C(G/B),C_0(G/T)\bigr)$, thus exhibiting $\K_*(\Gamma\ltimes_\mathrm{r}C(G/B))$ isomorphic to \linebreak $K^{*+\dim B/T}(\Gamma\backslash G/T)$.\\
It is not immediate to draw from this conclusions on the $\K$-theory of the algebraic crossed product since it is difficult to get hands on the map induced by the inclusion. However, it is possible to prove that the inclusion induces an isomorphism \emph{in analytic cyclic homology} by making use of Cuntz' picture of $\K$-theory. Namely there is a smooth bornological version of the functors $\q$ and $\Q$, allowing to find homomorphisms
$$
\beta:\Gamma\allcross \q\CINF(G/B)\to \Gamma\allcross M_n(\CCINF(Y))
$$
and
$$
\alpha:\q(\Gamma\allcross M_n(\CCINF(Y)))\to \Gamma\allcross \q(M_n(\CCINF(Y)))\to\Gamma\allcross \CINF(G/B)\protens\ell^p
$$
which induce isomorphisms in $\HA_*$. $Y$ is a certain homogeneous vector bundle over $G/B$, diffeomorphic to $G/T$, but endowed with another $\Gamma$-action. More precisely, one has to define several variants of Dirac and dual-Dirac morphisms. One has to define equivariant Bott periodicity in the form of an invertible element in $\HA_{*+\dim V}((A\protens \CCINF(V))\alcross B,A\alcross B)$. With some technical arguments, which contain the largest part of technical work in this paper, involving several times smooth homotopy invariance and excision, one can deduce that $\HA_*(\Gamma\allcross\CINF(G/B))\cong\HA_{*+\dim B/T}(\Gamma\allcross\CCINF(G/T))$. This is immensely helpful because the pairing only depends on the class of the idempotent in the latter cyclic homology.\\
All this done, there is a way to make use of the functoriality of the constructions we have done. They allow to put together Connes-Moscovici's index theorem and Atiyah's $L^2$-index theorem to deduce that the Chern character of an idempotent $e$ of $\C\Gamma$ (which of course only determines the pairing) is the same as the Chern character of an idempotent which is in the range of the assembly map. Unfortunately, the application of the Connes-Moscovici index theorem requires $\Gamma$ to be cocompact. This is the only place where we require $\Gamma$ to be cocompact.\\
This part of the paper is structured as follows. In the first section of the third chapter we review Nistor's construction. In the second, we investigate the splitting principle from various sides. In the fourth chapter, the Dirac-dual Dirac method is carried out. In its first section, we show how the Connes foliation index theorem would yield the result, if it applied.In the second., we show how in $\KK$-the necessary calculations are done. The third contains equivariant Bott periodicity. The fourth contains the core part, namely the calculation $\HA_*(\Gamma\allcross \CINF(G/B))\cong\Homol^*(\Gamma\allcross G/T)$. The fifth compares with the spectral sequence. The sixth contains the proof of the Kaplansky conjecture.\\
The results of this paper were communicated at VASBI conference on K-theory and Noncommutative Geometry, Valladolid (2006).\\
I am indebted to my thesis adviser Michael Puschnigg for most of the ideas underlying this paper and for the good cooperation, as well as to Heath Emerson for a fruitful discussion on the subject, during which I was also recommended to use the Furstenberg boundary, and Ralf Meyer for several hints.
\section*{$\KK$-theory and cyclic cohomology theories} 
\addcontentsline{toc}{section}{$\KK$-theory and cyclic cohomology theories}
In this section, we will sketch some of the ideas underlying the various cohomology theories used in this thesis, as well as some of their important properties.\\
Shortly after the definition of topological K-theory of locally compact Hausdorff spaces it was realized that the definition extended quite naturally to the category of possibly noncommutative Banach spaces. The extension to a bivariant $\KK$-functor, at least on separable $C^*$-algebras, was due to three major motivations. In fact, since K-theory is a generalized cohomology theory there exists necessarily a dual K-theory, baptized K-homology, defined in a hardly elegant way as the K-theory of the complement of the space embedded in some high $\R^n$, or homotopically with the help of spectra. It was Atiyah who realized that if the space is a smooth manifold then elliptic differential operators give functionals on K-theory, thus suggesting that these might be cycles for K-homology. However, he was not able to define the right homology relation. There was another main input was from Brown-Douglas-Fillmore who related K-homology to $C^*$-algebra extensions.
Cycles for $KK(A,B)$ are triples $(\Hilm,\pi,F)$, where
\begin{itemize}
\item
$\Hilm$ is a countably generated right $B$-Hilbert module
\item
$\pi:A\to\mathcal{L}\Hilm$ is a representation by Hilbert module maps.
\item
$F:\Hilm\to\Hilm$ is an operator of degree $1$, commuting with $B$.
\item
$\pi(a)(F-F^*)$ is compact for all $a\in A$.
\item
$\pi(a)(F^2-1)$ is compact for all $a\in A$.
\item
$[\pi(a),F]$ is compact for all $a\in A$.
\end{itemize}
The homology relation is, roughly spoken, homotopy. $\KK(A,B)$ are the morphisms in a category of $C^*$-algebras, and the functor from the ordinary category of $C^*$-algebras to that one can be characterized as the universal stable and split-exact functor on $C^*$-algebras. Such a functor is necessarily homotopy invariant by a result of Higson.\\
One of the central problems of $C^*$-algebra K-theory is the computation of the K-theory $K_*(\Cred G)$ of reduced $C^*$-algebras of locally compact groups $G$. It were Baum and Connes who realized that indices of $\Gamma$-equivariant elliptic operators on $G$-compact subsets of the universal classifying space for proper actions defined not only integers but elements of $K_*(\Cred G)$. They readily conjectured all these elements could arise this way, for all locally compact groups, and that there could be nor further equivalence relations among them but those coming from K-homology. Kasparov succeeded in proving the latter if $\Gamma$ is a discrete subgroup of a connected semisimple Lie group $G$. To that purpose, he defined equivariant KK-theory for pairs of separable $G$-$C^*$-algebras, having not merely triples $(\Hilm,\pi,F)$ as cycles but quadruples $(\Hilm,\rho,\pi,F)$ where all of the above axioms hold, and $\rho$ is a representation on the Hilbert module such that
\begin{itemize}
\item
$\rho$ is unitary in the sense
$$
(\rho(\gamma)s)(f)=\rho(\gamma)(s(\gamma\inv)f),\quad\gamma\in G,s\in\Hilm,f\in B.
$$
(note that $\rho$ seldom commutes with $B$, therefore)
\item
$\pi$ becomes covariant representation in the sense that
$$
\rho(\gamma)\pi(a)\rho(\gamma)\inv=\pi(\gamma.a)
$$
\item
$\pi(a)(\gamma F\gamma\inv-F)$ is compact for all $a\in A$.
\end{itemize}
Again, $KK^G$ viewed as a functor from the $G-C^*$-algebras and homomorphisms to $G-C^*$-algebras and $\KK^G$-groups can be characterized as the universal stable and split-exact functor on the category of $G-C^*$-algebras.\\
$\KK^G$ allowed to define the assembly map
$$
\lim_X\KK^G_*(C_0(X),\C)\to\K_*(\Cred G)
$$
the limit being taken over all $G$-compact subsets of the universal space, which for torsion-free discrete groups is ordinary $\classb G$.\\
The group $\K_*(\Cred\Gamma)$ is the receptacle group for the analytic index. Let us assume that a compact smooth $\BG$ exists. Then the class of an operator in the K-homology $\K_*(\BG)$ ) can analogously be seen as its topological index. The Baum-Connes conjecture is then a great extension of the program proposed by the Atiyah-Singer index theorem in the form that the analytical index agrees with the topological one. In K-homology terms and the language of the Baum-Connes conjecture one can state this as the fact that the following three maps coincide:
\begin{itemize}
\item
$K_*(M)\to\K_*(\pt)=\Z$, the map induced by shrinking $M$ to a point.
\item
The Poincar\'e duality followed by the shriek map of the collapsing map: $\K_*(M)\cong\K^*(T^*M)\to\K^*(\pt)=\Z$.
\item
The classifying map for the fundamental group followed by the assembly map followed by the trace functional
$$
\K_*(M)\to\K_*(\classs\pi_1M)\to\K_*(\Cred\pi_1M)\to\R.
$$
\end{itemize}
However, the Atiyah-Singer index theorem only gains its full power after translating it into ordinary cohomology, i. e. after computing the Chern character of the shriek map $\K^*(T^*M)\to\Z$. It is therefore apparent that it is very useful to dispose of of a version of de Rham cohomology for noncommutative algebras. This should be a generalized cohomology theory defined by a natural complex together with a pairing with K-theory. The above discussion implies that traces should be cocycles for that theory, and that the pairing should generalize the trace functional. This theory is, as one can guess, cyclic cohomology, and it was first defined by Connes and independently by Tsygan. See section \ref{IntroCyclic} for short introduction to cyclic homology. It turns out that cyclic homology only behaves well for much smaller algebras than $C^*$-algebras. For example on has $\HC^n(\CINF(M))\cong\Homol_{\ev/\odd}M$ for large $n$, but $\HC^0(C(M))$ is the space of Radon measures on $M$. (One has to use $\HC$ for Fr\'echet algebras, i.e. one has to replace $\otimes$ by $\protens$ for homology, and functionals by continuous functionals in cohomology). This is not surprising in view of the fact that we are dealing with a generalization of de Rham-theory, not Cech theory.\\
$\HC_*$ enjoys desirable properties like excision and smooth homotopy invariance only in a very limited sense. It is therefore better to pass to periodic cyclic homology, defined by the complex $\Pi_\ev \Omega A\leftrightarrows \Pi_\odd \Omega A$. This theory does satisfy the following properties (\cite{Pulocal}).
\begin{itemize}
\item
Invariance with respect to diffeotopies (smooth homotopies)
\item
Invariance under nilpotent extensions
\item
Algebraic Morita invariance
\item
Excision
\end{itemize}
Analytic cyclic homology was defined by Puschnigg and developed in great detail on the category of bornological algebras by Ralf Meyer (\cite{Meyerthesis}). It is possible to define analytic cyclic homology as the homology of the complex $\widehat{\Omega}_\ev A \leftrightarrows \widehat{\Omega}_\odd A$, where $\widehat{\Omega}$ is the bornological completion of the complex of noncommutative differential forms. It is also invariant under topologically nilpotent extensions, satisfies excision, and is Morita invariant. The importance of analytic cyclic homology in this work lies in the fact that we can show that it is not degenerate neither on small algebras like the algebraic crossed product nor on $C^*$-algebras.\\
Moreover, there is a Chern character both from K-theory and K-homology. However, there is no bivariant one. This is remedied by local cyclic homology, due to Puschnigg. $\HPL$ enjoys the following two crucial properties.
\begin{itemize}
\item
There is a bivariant Chern-Connes character $\KK\to\HPL$ in full generality, i. e. a multiplicative natural transformation of bifunctors. This is called a generalized Grothendieck-Riemann-Roch theorem.
\item
On commutative $C^*$-algebras one has $\HPL_*(C_0(M))=\Homol_c^{2\Z+*}(M)$. So whenever a $\KK$-group is computed by a KK-equivalence with a commutative $C^*$-algebra, then so is $\HPL$ by the preceding property.
\item
It is a fact proved in \cite{Pulocal} that for a large class of algebras, that $C^*$algebras belong to, local cyclic homology coincides with analytic cyclic homology (Caution: This is not true in cohomology).
\end{itemize}
\thispagestyle{Selbergemacht2}

}
\clearpage{\pagestyle{empty}\cleardoublepage}
\chapter{Equivariant K-theory of the visibility boundary}
Let $G$ be a real connected semisimple Lie group, $K$ a maximal compact subgroup. Then the homogeneous space $G/K$ is a globally symmetric Riemannian space of non-positive sectionl curvature, see appendix \ref{symmetric space}.\\
Throughout this thesis, $i$ will denote the dimension of $G/K$. Since $G$ acts freely and properly on itself, so does a discrete and torsion-free subgroup $\Gamma$ on $G$ and even on $G/K$. Since $G/K$ is contractible, we have found a model for the classifying space $\EG$ for principal fibre bundles, and we also have $\BG=\Gamma\backslash \EG=\Gamma\backslash G/K$. In the sequel, $\EG$ resp. $\BG$ will stand for the particular models $G/K$ resp. $\Gamma\backslash G/K$ with all extra structures.\\
We are now going to define the boundary of $G/K$.
\begin{defn}
Let $X$ be a manifold of non-positive sectional curvature. Two unit-speed geodesic rays in $X$ are called asymptotic if there is a constant $a\in\R$ such that
$$
d(a_1(t),d(a_2(t))\le a \text{ for all } t\in\R_+.
$$
The boundary $\partial(X)$ is the set of equivalence classes of geodesic rays. This set is topologized by identifying $X\cup \partial X$ with a fiber of the unit disk bundle of $TX$ via the exponential map and $\partial X$ therefore with a fiber of the unit sphere bundle.
\end{defn}
If $\Gamma$ is cocompact, then $\Gamma$ is hyperbolic if and only if the real rank of $G$ is one. This follows from the facts that the inclusion $\Gamma\to G$ is a quasi-isometry, and that the sectional curvature of $G$ attains zero if and only if there is a two-dimensional flat in $G$. If $\Gamma$ is cocompact and hyperbolic, then the Gromov boundary $\partial\Gamma$ coincides with $\partial (G/K)$ (\cite{Ghys}, th\'eor\`eme 37). The action of $\Gamma$ on geodesics passes to the quotient because $\Gamma$ acts by isometries, so we have an action $\Gamma\to \AutTop(\partial(G/K))$. This is the action we will mainly be interested in. We remark that it is always ergodic (see\cite{Zimmer}). Furthermore, in many considered cases all orbits are dense so that the naive quotient $\Gamma\backslash \partial(G/K)$ has a very bad topology. Instead, we will look at the homotopy quotient $\hqv$. This space is defined to be the quotient of $\EG\times\partial(\EG)$ by the diagonal action $\gamma(a,b):=(a\gamma\inv,\gamma b)=(\gamma a, \gamma b)$.
\begin{lemma}
The space $\hqv$ is homeomorphic to the unit sphere bundle $S(\Gamma\backslash G/K)$ of $\BG$.
\end{lemma}
\begin{proof}
The unit sphere bundle of $G/K$ is trivial since $G/K$ is contractible. So it is homeomorphic to $G/K \times \partial (G/K)$. There is a canonical $\Gamma$-equivariant homeomorphism which assigns to a tangent vector its base point and direction. On the other hand, using the definition of the boundary as given above, one observes that under this homeomorphism the action of $\Gamma$ on the unit sphere bundle corresponds to the diagonal action of $\Gamma$ on the product. So the homeomorphism descends to the quotient.
\end{proof}
Our object of study is the $C^*$-algebraic reduced crossed product $\Gamma\rlcross C(\partial(G/K))$ defined by the action described above.\\
We shall make use of the Baum-Connes assembly map (\cite{Co3}) in the form
$$
\mu: \RKKG(C_0(\EG), A)\longrightarrow \K_*(\Gamma\rlcross A),
$$
as described, for example, in \cite{Valette}. Assuming an equivariant Spin$^c$-structure on $\EG$, we shall write and use $\mu$ in the form
$$
\KK_{i+j\!\!\!\mod 2}(\C,\Gamma\rlcross (C_0(\EG)\otimes A))\longrightarrow \K_j(\Gamma\rlcross A)
$$
In our context, this map is split injective (\cite{Kas}).\\
The strategy of the proof is to realize the class of unit ``geometrically'', i. e. to find a canonical preimage $x$ of $[1]$ under $\mu$. Of course, the finiteness of the order of $x$ is sufficient for that of its image $\mu(x)$ under $\mu$ without any assumption on $\mu$. If $\mu$ is rationally injective, the converse also holds. So for the sufficiency of our criterion for the finiteness of the order of $[1]$, we shall only need the mere existence of $\mu$. This was essentially the line of reasoning Connes used in \cite{CoCC}. For the necessity of the criterion, we shall need to make use of $\mu$'s rational injectivity.\\
\pagestyle{fancy}
\section{The Dirac and dual-Dirac elements for $\Gamma$}\label{DdDG}
In this section, we will describe the construction of the Dirac element $\alpha$ for $\Gamma$, $\alpha\in\KK^\Gamma(\Sec\Cliff T^*\EG,\C)$, and how to reduce it to $\KK^\Gamma(C_0(\EG),\C)$. We need $\alpha$ in the latter form because it allows to write the assembly map in most geometric form.\\
Almost all of the contents of this section are based on \cite{Kas}.\\
We shall make use of the classical Dirac element for the group $\Gamma$ as defined by Kasparov, \cite{Kas}*{4.2}. It is an element of
$$
\KK^\Gamma(\Sec(\Cliff(T^*\EG\otimes \C)),\C),
$$
where $\Sec$ denotes the continuous sections of a bundle vanishing at infinity, and $\Cliff(T^*\EG\otimes \C)$ denotes the Clifford algebra bundle of the complexified cotangent bundle.\\
The Hilbert space defining this $\KK$-element is $L^2\Lambda^*T^*\EG$, the space of square-integrable differential forms with respect to the measure defined by the metric. This Hilbert space is equipped with a canonical $\Gamma$-action induced from the $\Gamma$-action on $\EG$ by translations. The algebra $\Sec(\Cliff(T^*\EG\otimes \C))$ can be represented on this Hilbert space by 
$$
\pi(f)\xi=\mathrm{ext}(f)(\xi)+\mathrm{int}(f)(\xi)
$$
for $f$ a one-form on $\EG$. By calculation, one obtains
$$
\pi(f)^2=\mathrm{ext}(f)\mathrm{int}(f)+\mathrm{int}(f)\mathrm{ext}(f)=||f||^2
$$
which assures by the universal property of the Clifford algebra that $\pi$ passes to a representation of $\Sec(\Cliff(T^*\EG\otimes \C))$ on the Hilbert space. Finally the operator $D$ is given by functional calculus with $x\mapsto\frac{x}{\sqrt{1+x^2}}$ applied to $\mathcal{D}=d+d^*$, making it a bounded operator. This defines the data needed for the $\KK$-cycle, and it remains to verify the required properties. $\mathcal{D}$ is essentially self-adjoint, so $D$ is self-adjoint. $\mathcal{D}$ is $\Gamma$-invariant, so $D$ is $\Gamma$-invariant. Let $a\in\Sec(\Cliff(T^*\EG\otimes \C))$. Then $a(1-D^2)=a(\frac{1}{1+\mathcal{D}^2})$ which, roughly speaking, is compact because of the following argument. The function $\frac{1}{1+x^2}$ is norm-limit of compact support-functions. In a compact interval of $\R$ can be only finitely many eigenvalues of an elliptic operator, since its spectrum is discrete, moreover they are all of finite multiplicities. So $\frac{1}{1+\mathcal{D}^2}$ is limit of a sequence of finite rank operators, thus it is compact (more precisely, one has to invoke Rellich's lemma). The hard part of the proof then is to show that $[a,D]$ is compact, too. This is done by exhibiting $[a,D]$ as an integral converging in norm of compact operators. Thus the $\alpha$-element is defined.\\
We will now describe how to reduce the element $\alpha'\in\KK^\Gamma_q(\Sec(\Cliff(T^*\EG\otimes \C)),\C)$, where, as always, $q=\dim\Gamma$, to an element $\alpha^\Gamma\in\KK^\Gamma_q(C_0(\EG,\C))$, and similarly for $\beta'\in \KK^\Gamma_q(\C,\Sec\Cliff T^*\EG)$. This is what we shall need in view of the fact that we have to carry out topological constructions with $\alpha$.\\
In fact, denote by $\Spinor$ the spinor bundle on $\EG$ such that $\Cliff T^*\EG\cong\End\Spinor$. Then the space of sections of $\Spinor$ vanishing at infinity provides a $\Gamma$-imprimitivity bimodule between $\Sec\Cliff T^*\EG$ and $C_0(\EG)$, thus an invertible element
$$
\Morita\in\KK^\Gamma_0(\Sec\allowbreak\Cliff T^*\EG, C_0(\EG))
$$
(see \cite{Kas}*{2. 18}). We thus define
$$
\alpha^\Gamma=\Morita\inv\otimes\alpha'\text{, and }\beta^\Gamma=\beta'\otimes\Morita.
$$
The element $\alpha^\Gamma$ has then following properties. First, $\alpha^\Gamma\otimes\beta^\Gamma=1\in\KK^\Gamma_0(C_0(\EG),C_0(\EG))$, which follows directly from the definitions, the associativity of the Kasparov product and $\alpha'\otimes\beta'=1$. So the right product with $\alpha^\Gamma$ is injective.\\
Second, the element
$$
\rest_e^\Gamma(\beta^\Gamma\otimes\alpha^\Gamma)=\rest_e^\Gamma(\beta')\otimes\rest_e^\Gamma(\alpha')=\rest_e^\Gamma(\gamma)\in\KK_0(\C,\C)\allowbreak =\Z
$$
equals $1\in\Z$, where $\rest_e^\Gamma$ denotes the forgetful map $\KK^\Gamma\to\KK$ defined by restriction along $e\to\Gamma$. The fact that $\rest_e^\Gamma(\gamma)=1$ can be seen as follows. It is an idempotent in $\KK_0(\C,\C)$, so it is either $0$ or $1$. The first possibility is excluded by the calculation
$$
1=\rest(1)=\rest(\alpha' \otimes\beta' \otimes\alpha' \otimes\beta')=\rest(\alpha')\otimes\rest(\gamma)\otimes\rest(\beta').
$$
This means that $\rest_e^\Gamma(\beta^\Gamma)$ and $\rest_e^\Gamma(\alpha^\Gamma)$ are inverse Bott elements, up to sign. By the orientation conventions, this sign is $+1$.\\
Whenever $G$ is simply connected, for instance $G=\SL(n,\C)$, then the $\Spin^c$-condition is automatic (see the proof of theorem \ref{SPINN1} and remark \ref{SPINN2}). Note also that Connes shows in \cite{CoCC}*{p. 37} that in the case $G=\SL(n,\R),K=SO(n)$ for even $n$ the representation $\K\to\SO(\mathfrak{g/k})$ even lifts to Spin($\mathfrak{g/k}$), so in particular it preserves a Spin$^c$-structure.

\section{The main theorem on the class of unit}
\begin{thm}\label{mainKtheorem}
Let $G$ be a connected semisimple Lie group, $K$ a maximal compact subgroup, $G/K$ the associated globally symmetric space of non-positive sectional curvature, $\Gamma$ a torsion-free cocompact lattice of $G$, $\Gamma\to \AutTop(G/K)$ the left $\Gamma$-action coming from the natural action of $\Gamma$ on $G$, $\Gamma\to \AutTop(\partial(G/K))$ its extension to the boundary of the symmetric space. Assume furthermore that $\EG$ is endowed with a $\Gamma$-invariant Spin$^c$-structure, and that $\Gamma$ is cocompact.\\
Then the following are equivalent:
\begin{enumerate}[(1)]
\item
The unit $[1]\in \K_0(C(\partial(\EG))\rcross\Gamma)$ is a torsion element.
\item
The K-theoretic fundamental class $[\BG]\in\K^{q\mod 2}(\BG)$ is mapped by
$$
\pi^*:\quad\K^q(\BG)\to\K^q(S(\BG))
$$
to a torsion class.
\item
The orientation class or fundamental class of $H^*(\BG;\rationals)$ is in the kernel of
$$
\pi^*:\quad\Homol^*(\BG;\rationals)\to \Homol^*(S(\BG);\rationals).
$$
\item The Euler characteristic $\chi(\Gamma)=\chi(\BG)$ of $\Gamma$ does not vanish.
\item The Euler-Poincar\'e measure of $G$ is non-zero.
\item The Lie algebras $\mathfrak{g}$ and $\mathfrak{k}$ of $G$ and $K$ have same rank.
\end{enumerate}
\end{thm}
See appendix \ref{EP} for a definition of the Euler-Poincar\'e-measure.
\begin{rmk}\mbox{ }\\
\vspace{-1em}
\begin{enumerate}
\item
The fifth and sixth criteria are independent on $\Gamma$.
\item
The condition on $G$ appearing in statement $(6)$ is much easier to verify than the first four statements. It holds for any complex semisimple Lie group $G$. In particular, it holds for $G=\SL(n,\C), n\ge 2$. The reason is the following. Since the Lie algebra $\mathfrak{g}$ of of a complex semisimple Lie group $G$ has a compact real form $\mathfrak{g}_0$, which is unique up to isomorphism, and since Cartan subalgebras in the real and complex sense correspond to each other, we have $\rk \mathfrak{g}=\rk \mathfrak{g}_0\otimes\C=2\thinspace\thinspace\rk \mathfrak{g}_0$. So both ranks are always different, since they never vanish. On the other hand, the condition is not met by the real semisimple Lie group $G=\SL(2,\R)$. Here, $\rk \mathfrak{g}= \rk\mathfrak{k}=1$.\\
This achieves one of our aims, namely an ``explanation'' of Connes' calculation (\cite{CoCC}) that the class of unit has torsion in case $G=\SL(2,\R)$, and does not have torsion if $G=\SL(2,\C)$.
\item
We will prove below, that when $\Gamma$ is not cocompact, there is only one possibility, namely that the statements $(1),(2)$ and $(3)$ do not hold for \emph{any} real semisimple connected $G$. Statement $(4)$ does not make sense in that case, and consequently there is nothing to say on statements $(5)$ and $(6)$ then.
\item
The condition on the Lie group $G$ appearing in the theorem also shows up in quite different contexts. For instance, it is equivalent to $G$ having a discrete series representation (\cite{AtiyahSchmid}).
\end{enumerate}
\end{rmk}

\emph{Proof} (of the theorem).\\ \nopagebreak
As we have assumed that $\Gamma$ is cocompact, the Riemannian manifold $\BG=\Gamma\backslash G/K$ is compact. However, we want to the part of the proof that shows $(1)\Leftrightarrow(2)\Leftrightarrow(3)$ to apply also in the noncompact case that we will look at below. So we remind the reader that we always use K-theory and cohomology with compact supports.
\begin{enumerate}
\item[$(1)\Leftrightarrow (2):$]
We set up a commutative diagram
$$
\xymatrix@C=9em{
\K_q(C_0(\BG))\ar[d]^{\pi^*}\ar[r]^{-\otimes  j(\alpha)}&K_0(\Cred\Gamma)\ar[d]\\
\K_q(C_0(\hqv))\ar[r]^{-\otimes j(\alpha\otimes 1_{C(\partial(\EG))})}&\K_0(\Gamma\rlcross C(\partial(\EG)),
}
$$
where for simplicity the isomorphisms coming from the strong Morita equivalence, i. e. the invertible element in $\KK_0(\Gamma\ltimes C_0(\EG\times X),C_0(\EG\times_\Gamma X))$ are suppressed in the notation. Recall that we view $\alpha$ being in $\KK^\Gamma_q(C_0(\EG),\C)$, and that $\pi:\partial(\EG)\to\pt$ defines an element in $\KK^\Gamma_0(\C,\allowbreak C(\partial(\EG))$. Then we compute with the definition of the cup product (this is the case $-\otimes_\C-$ of the Kasparov cup-cap product) and its commutativity the Kasparov element that describes the composition $(-\otimes j(\alpha\otimes \partial(\EG)))\circ\pi^*$:
\begin{multline*}
j(1_{C_0(\EG)}\otimes_\C \pi^*)\otimes_{\Gamma\rlcross C_0(\EG\otimes\partial(\EG))} j(\alpha\otimes 1_{C(\partial(\EG))})=\\
j((1_{C_0(\EG)}\otimes_\C\pi^*)\otimes_{C_0(\EG\otimes\partial(\EG))}(\alpha\otimes 1_{C(\partial(\EG))}))=\\
j(\pi^*\otimes_\C \alpha)=j(\alpha\otimes_\C \pi^*)=j(\alpha)\otimes_{\Cred\Gamma} j(\pi^*)
\end{multline*}
and by the associativity we have shown that the diagram commutes.\\
Now the lower arrow is split injective by Kasparov's equality $\alpha\otimes\beta=1_{C_0(\EG)}$. In fact,
\begin{multline*}
j(\alpha\otimes 1_{C(\partial(\EG))})\otimes j(\beta\otimes 1_{C(\partial(\EG))})=j((\alpha\otimes\beta)\otimes 1_{C(\partial(\EG))})\\
=j(1_{C_0(\EG)}\otimes 1_{})=1
\end{multline*}
and right multiplication with $j(\beta\otimes 1_{C(\partial(\EG))})$ provides a left inverse to the assembly map. The action of the $\gamma$-element $j(\beta\otimes\alpha)\in\KK_0(\Cred\Gamma,\Cred\Gamma)$ on $\K_*(\Cred\Gamma)$ is the projection onto the image.\\
So it only remains to identify a preimage of the class of the unit in $\K_0(\Cred\Gamma)$.\\
The assembly map in the form $\mu:\RK_0(\BG)\to\K_0(\Cred\Gamma)$, which is the easiest case of the definition of Baum, Connes and Higson, factorizes as K-theoretic Poincar\'e duality $\RK_0(\BG)\cong\K^q(\BG)$ followed by $-\otimes j(\alpha)$ (\cite{Puschnigg}, \cite{Kas}). Now Valette (\cite{Valette}) shows that $\mu([i_*])=[1]$, where $[i_*]$ is the class in $\RK_0(\BG)$ that corresponds to the inclusion of a point in $\BG$. So the statement follows from the fact that the orientation class in $\K^q(\BG)$ is the Poincar\'e dual of the point class in $\RK_0(\BG)$.\qed
\item[$(2)\Leftrightarrow(3):$]
The fact that the Chern character $\K^*(X)\to \oplus \Homolc^{*+2i}(X)$ is a rational isomorphism implies that it suffices to show that it maps the $\K$-theoretic fundamental class to the ordinary cohomological one, up to sign. Topology teaches that both fundamental classes can be defined as follows: Denote by $[\mathbb{S}^q]$ the fundamental class of the q-sphere, i. e. its orientation class or the Bott class (be it even or odd). Choose an arbitrary point in $\BG$ and map a neighborhood which is diffeomorphic to the i-disc in such a way to the sphere that the complement of the chosen neighborhood is mapped to the north pole, and the base point is mapped to the south pole, and the resulting map from the one-point compactification of the interior of the neighborhood to the sphere has degree one. Then $[\BG]=f^*([\mathbb{S}^q])$. This follows from mapping the long exact sequence associated to the pair $(\mathbb{S}^q,\pt)$ in any generalized cohomology theory to the one associated to $(\BG,\BG-\mathbb{D}^q\cup\partial\mathbb{D}^q)$. The generator of $\Homol^q(\BG,\BG-\mathbb{D}^q\cup\partial\mathbb{D}^q)=\Homol^q(\BG,\BG-\pt)$ is mapped to the orientation class in $\Homol^q(\BG)$. So the naturality of the Chern character and the well-known fact that the Chern character of the Bott class is, up to sign, the orientation class - which, in turn is proved by reducing it by the long exact sequence relating reduced and unreduced cohomology of the sphere to the case of $\R^q$ - together achieve the proof.\qed
\item[$(3)\Leftrightarrow (4):$]
Let us recall the definition and fundamental properties of the Euler class as well as the Gysin sequence
$$
\dots\to \Homol^m(\BG)\xrightarrow{-\cup e}\Homol^{m+q}(\BG)\xrightarrow{\pi^*}\Homol^{m+q}(S\BG)\to\Homol^{m+1}(\BG)\to\dots,
$$
where $e\in H^i(\BG)$ is the Euler class of $\BG$.\\
So let $\pi:E\to X$ be an arbitrary euclidean vector bundle (the multiple use of the letter $\pi$ should cause no confuse) over a locally compact Hausdorff space $X$, and let $SE$ and $DE$ denote its unit sphere and unit disc bundle, respectively. Let $i=\dim X, n=\rk E$.\\
Let $\Homol^*$ be any generalized cohomology theory with arbitrary supports, and assume that $E$ is orientable with respect to $\Homol$. This means that there is a Thom class $u\in\Homol^n(E,E-x)$, the cohomology of the total space of $E$ with compact supports in the vertical direction, (or $u_{\Homol}$ if the cohomology theory is to be specified) such that $\Homol^*(E,E-X)$ is a free module over $\Homol^*(X)$ (with arbitrary supports) of rank one with generator $u$. If the base space is compact, then $\Homol^*(E,E-X)$ is the compact support cohomology of the total space of $E$. Here, the module structure is defined by the scalar multiplication $\pi^*(-)\cup u$ where $\cup$ is the pairing $\Homol^p(E)\otimes \Homol^m(E,E-X)\to \Homol^{p+m}(E,E-X)$, and $\pi^*:\Homol^p(X)\to\Homol^p(E)$. If one prefers to deal only with compact support cohomology and proper maps, one can write down the module structure as follows. Let $m:E\to M\times E$ be the proper map $e\mapsto (\pi(m),e)$. Then the expression $m^*(-\times u)$ defines the $\Homolc^*(X)$-module structure on $\Homolc^*(E)$, where $-\times -:\Homolc^p(X)\otimes\Homolc^m(E)\to \Homolc^{p+m}(X\times E)$ denotes the exterior or cross-product. If $X$ happens to be compact, then one easily checks that both module structures agree. The second expression enjoys a Thom isomorphism only if it agrees with the first one, i. e. in case $X$ is compact. So we will deal with the noncompact case separately.

So let $X$ be compact.\\
Consider the relative group $\Homol^k(DE,SE)$ in the theory $\Homol$. For example, if $\Homol=\K$, then this is just the 6-term exact sequence associated to the ideal given by the continuous functions on $DE$ that vanish on $SE$ in the algebra of all functions on $DE$.\\
By the cofibration property, this is equal to $\Homol^k(DE/SE,SE)=\Homol^k(E^+,\infty)$, the compact support cohomology of the total space of $E$, which in turn is isomorphic by the Thom isomorphism theorem to the cohomology $\Homol^{k-n}$ of the base space $X$. (Equivalently, the pair (DE,SE) is homotopy equivalent to (E,E-X)). Furthermore, $DE$ can be retracted to $X$, and under all these isomorphisms the maps in the long cohomology exact sequence associated to the pair $DE,SE$ can be rewritten: The map $SE\to DE$ becomes the projection $\pi:SE\to X$ (which is proper). The map $(DE,SE)\to DE$ becomes multiplication in the indicated sense with the Thom class $u\in\Homol^n(E^+,\infty)$, $-\times u:\Homolc^{k-n}(X)\to\Homolc^k(E)$ followed by restriction to $\Homol^k(X)$. So the map $\Homol^{k-n}(X)\to\Homol^k(E)\to\Homolc^k(X)$ is expressed by\\
$$
i^*(m^*(-\times u))=m^*(\mathrm{id}\times i)^*(-\times u)=m^*(-\times i^*(u))=-\cup i^*(u),
$$
which motivates the Euler class.
\begin{defn}
$i^*(u)\in\Homol^k(X)$ is called the \emph{Euler class} of the vector bundle $E$ and commonly denoted $e$ (or $e_{\Homol}$ if the cohomology theory is to be specified). If $X$ is a Riemannian manifold and $E=TX$ then $e$ is called the Euler class of $X$.\\
If $\Gamma$ is a discrete group which is the fundamental group of a smooth compact aspherical manifold $\BG$, then the Euler characteristic in ordinary cohomology of $\BG$ is called the Euler characteristic of $\Gamma$.
\end{defn}
The name (and the well-definedness) comes from the classical \emph{Gauss-Bonnet theorem}: Let $E=TX$. Then
$$
<e,[X]>=\int_X e=\chi(X),
$$
the Euler characteristic of $X$, defined by $\chi(X)=\sum_i (-1)^i a_i$, where $a_i$ is the number of cells in dimension $i$ of a triangulation of $X$, or by $\chi(X)=\sum_i (-1)^i\dim H_i(X)$. In other words, the theorem states that $e=\chi(X)[X]$, because $e$ is concentrated in top degree. We will lateron show an analogous property in K-theory.\\
The preceding discussion, together with the definition of the Euler class, yield the existence of an exact sequence, called the \emph{Gysin sequence}:
$$
\dots\to \Homol^k(X)\xrightarrow{-\cup e}\Homol^{k+n}(X)\xrightarrow{\pi^*}\Homol^{k+n}(SE)\to\Homol^{k+1}(X)\to\dots.
$$
which establishes the equivalence of the statements (3) and (4).\qed
\item[$(4)\Leftrightarrow(5)\Leftrightarrow\\
(6):$] These equivalences are classical. We will give an account of the arguments in appendix \ref{EP}.\qed
\end{enumerate}
Let us now look at the noncompact case. We shall not even require $\Gamma$ to be of finite covolume.
\begin{thm}
In the situation of the preceding theorem, let $\Gamma$ not be cocompact. Then the following statements hold:
\begin{enumerate}[(1)]
\item
The unit $[1]\in \K_0(C(\partial(\EG))\rcross\Gamma)$ is not a torsion element.
\item
The K-theoretic fundamental class $[\BG]\in\K^{q\mod 2}(\BG)$ is mapped by
$$
\pi^*:\quad\K^q(\BG)\to\K^q(S(\BG))
$$
to a class of infinite order.
\item
The (orientation) fundamental class of $H^*(\BG;\rationals)$ is not in the kernel of
$$
\pi^*:\quad\Homol^*(\BG;\rationals)\to \Homol^*(S(\BG);\rationals).
$$
\end{enumerate}
\end{thm}
\emph{Proof.}\\
The proof of the preceding theorem remains valid to show that $(3)\implies (2)\implies (1)$. Let us show that $(3)$ holds for any Lie group $G$. The model $\Gamma\backslash G/K$ for $B\Gamma$ is noncompact as well, and the fibration $\EG \times_\Gamma \partial(\EG)\to\BG$ has a noncompact base space. There need not always be a compact model for $\BG$. For instance, the fundamental group of a Riemannian surface of infinite genus is a subgroup of $\SL(2,\R)$ (However, even if there were such a compact $\BG$ as well, it would not be useful, because $S\BG$ obviously depends on the diffeomorphism type of $\BG$).
\begin{lemma}
Let $X$ be an \emph{open} Riemannian orientable $i$-dimensional manifold, $E=TX$ its tangent bundle. Then the fundamental class of $\Homolc^q(X)$ is not annihilated by $\pi^*:\Homolc^q(X)\to\Homolc^q(SX)$.
\end{lemma}
\begin{proof}
Consider the K-theory 6-term exact sequence associated to the ideal $C_0(TX)$ of $C_0(DX)$
$$
\xymatrix{
{\K^0(TX)}\ar[r]^{i^*}&{\K^0(DX)}\ar[r]^{\pi^*}&{\K^0(SX)}\ar[d]\\
  {\K^1(SX)} \ar[u] & {\K^1(DX)} \ar[l]^{\pi^*} & {K^1(TX)} \ar[l]^{i^*}.
}
$$
Tensoring with $\Q$, replacing $DX$ with the properly homotopy equivalent $X$, and taking into account that the Chern character is a rational isomorphism, we get again an exact sequence
$$
\Homolc^\ev(TX)\xrightarrow{i^*}\Homolc^\ev(X)\xrightarrow{\pi^*}\Homolc^\ev(SX)
$$
and similarly for odd cohomology. $i^*$ and $\pi^*$ preserving the degree, it follows that for any $k$ there is an exact sequence
$$
\Homolc^k(TX)\xrightarrow{i^*}\Homolc^k(X)\xrightarrow{\pi^*}\Homolc^k(SX).
$$
If we rewrite the first map $i^*:\Homolc^k(TX)\to\Homolc^k(X)$ under Poincar\'e duality and use homotopy invariance, we get a map $\Homol_{2i-k}(X)\to \Homol_{i-k}(X)$ that lowers degree by $i=\dim X$. An open $i$-manifold always has the homotopy type of its $i-1$-complex. So $\Homol_i(X)$ is zero. So since $\Homol_q(X)$ can only be non-zero for $0\le q\le i-1$, the map $i^*$ must be zero.\\
We conclude that neither of the statements $(1),(2),(3)$ holds, and the class of unit in $\K_0(\Gamma\rlcross C(\partial(\EG)))$ is not torsion, whenever $\Gamma$ is not cocompact.
\end{proof}
\section{The exact order of the class of unit}
We now want to assume that $\Gamma$ is cocompact and that its Euler characteristic is non-zero (in particular, $i$ is even), so that the unit class is torsion, and want to compute the order of the latter in $\K_0(\Gamma\rlcross C_0(\partial(\EG))).$
\begin{thm} 
In the situation of the main theorem, let $\Gamma$ be cocompact and let the equivalent statements of the theorem hold. Then the class of unit
$$
[1]\in\K_*(C(\partial(\EG)))\rcross\Gamma)
$$
is of order $|\chi(\Gamma)|$, the absolute value of the Euler characteristic.
\end{thm}
In particular, the class of $[1]$ is non-zero.
\emph{Proof.}
The discussion of the Gysin sequence in the proof of the equivalences $(3)\Leftrightarrow(4)$ of the main theorem on the class of unit applies also in the case where the generalized cohomology theory $\Homol$ is $\K$-theory. We thus get a cyclic 6-term exact sequence:
$$
\xymatrix{
  {\K^0(X)} \ar[r]^>>>>>>{-\cup e}  & {\K^n(X)} \ar[r]^{\pi^*} & {\K^n(SE)}\ar[d]\\
  {\K^1(SE)} \ar[u] & {\K^1(X)} \ar[l]^>>>>>>{\pi^*} & {K^1(X)} \ar[l]^>>>>>>>{-\cup e}.
}
$$
(One may have the impression that this concludes the proof. However, the following arguments are still necessary.)\\
Recall that we assumed $\EG$ to be $\Gamma$-invariantly Spin$^c$, or equivalently that $\BG$ is Spin$^c$.\\
First we observe that the Chern character of the K-theoretic Euler class is, up to sign, equal to the cohomological one. In fact,
\begin{multline*}
\ch(e_{\K})=\ch(i^*u_{\K})=i^*(\ch(u_{\K}))=i^*(\tau (\BG)u_{\Homol})=\\
i^*(\sign u_{\Homol}+y)=\sign (i^*(u_{\Homol})+i^*y)=\sign i^*(u_{\Homol})=\sign e_{\Homol},
\end{multline*}
  where $e_K, e_H, u_K, u_H$ are the K-theoretic resp. cohomological Euler resp. Thom classes and $\tau$ is a characteristic class beginning with $\sign$. This means that $\tau=\sign[pt]+y$, $[pt]$ being the pullback of the chosen generator of the cohomology of a point along the map $\BG\to\pt$, and $y$ being a cohomology class whose component zero part vanishes. This follows from the equality (see \cite{Summerschool})
 $$
 \tau=\frac{\sign}{\Todd\overline{T\BG}},
 $$
 where $\Todd$ denotes the Todd genus of a Spin$^c$-vector bundle.\\
 There is a diagram, the first line being exact as a part of the 6-term sequence established above:
 $$
 \xymatrix@C=5em{
 \K^0(\BG)\ar[d]\ar[r]^{\cup e}&\K^0(\BG)\ar[d]\ar[r]&\K_0(\Gamma\rlcross C(\partial(\EG)))\\
 \Homol^{\ev}\ar[r]^>>>>>>>>>>>{\cup \sign e}&\Homol^{\ev}(\BG).
 }
 $$
 and the upper right arrow maps $[\BG]$ to $[1]$.\\
 By the rational injectivity of the Chern character and $\ch([\BG]_{\K})=[\BG]_{\Homol}$, we conclude that $e_{\K}=\chi(\BG)[\BG]+t$, where $t$ is a torsion element. The following lemma shows that $t=0$, because $[\BG]$ is an element of infinite order.\\
 \begin{lemma}
 The element $e_{\K}=i^*(u)$ is in the subgroup of $\K^*(\BG)$ generated by $[\BG]$.
 \end{lemma}
 \begin{proof}
 The space $\BG$ being compact and orientable, it can be realized simplicially as a complex with a single top cell. This means that $\BG/(\BG)^{q-1}$ is homeomorphic to $\mathbb{S}^q$, where $\BG^{q-1}$ is the $(q-1)$-skeleton of a cell decomposition of $\BG$. Denote by $j:\BG^{q-1}\hookrightarrow \BG$ the inclusion and by $f:\BG\to\mathbb{S}^q$ the projection.\\
 Let
 $$
 \dots\to{\K^0}(\mathbb{S}^q)\xrightarrow{f^*}\K^0(\BG)\xrightarrow{j^*}\K^0(\BG^{q-1})\to 0
 $$
 be the long exact sequence associated to the pair $(\BG,(\BG)^{q-1})$.\\
 Let us show that $e$ is in the kernel of $j^*$.\\
 Now, the tangent bundle $T\BG|_{\BG^{q-1}}$ restricted to $\BG^{q-1}$ has by dimension reasons zero Euler class, so there is a nowhere vanishing section $s:\BG^{q-1}\to T\BG$. The section $s$ is homotopic to the zero section $t:\BG^{q-1}\to T\BG$, because so is any section of the tangent bundle. We have
 $$
 j^*(e_{\K})=j^*(i^*(u))=(ij)^*(u)=t^*(u)=s^*(u).
 $$
 But $u\in\K^*(T\BG)$ is represented by a complex which is exact outside the zero section, so $s^*(u)$ is zero.

 \end{proof}
 We have shown that $e_{\K}=\sign\chi(\BG)[\BG]$. This implies that the order of the class of unit is at most $|\chi(\BG)|$. If it were less, there would be an element $x$ of $\K^0(\BG)$ which gets mapped by $-\cup e_{\K}$ to $z[\BG]$ for $1\le z\le |\chi(\BG)|-1$. But if we had $x\cup e =z[\BG]$, then there would be the equality $(\ch(x))_0 \chi(\BG)[\BG]=z[\BG]$ where $\ch(x)_0$ is the degree zero part of $\ch(x)$ which is always an integer. So there would be a contradiction. This concludes the proof of the theorem.\qed
 \begin{rmk}
 We can even make precise the sign of $\chi(\Gamma)=\chi(\BG)$ in case $\rk G=\rk K$. The proof contained in the preceding section shows that the Euler characteristic of the compact dual $G'/K'$ is positive, and that the curvature tensor $R$ of $G/K$ is minus $R'$, the one of $G'/K'$. So the Euler density $\Omega$ of $G/K$, being linear in $q/2$ powers of $R$, has sign $\sign$ times that of $G'/K'$, and $\chi(G'/K')\ge 1$. So we deduce that $\chi(\Gamma)$ has sign $\sign$.
 \end{rmk}

\clearpage{\pagestyle{empty}\cleardoublepage}
\chapter{Equivariant Cyclic Cohomology of the Furstenberg boundary}
\section{Nistor's construction from the topological point of view}
\subsection{Simplicial Methods, Cyclic Modules and Group Homology}\label{IntroCyclic}
We are now going to introduce shortly cyclic and Hochschild homology, and then give an introduction to Nistor's construction, which is a nice demonstration of the interplay of notions from group theory, homology, cyclic homology, homotopy theory and homological algebra encountered in the subject.\\
The \emph{simplicial category} $\Delta$ is the category whose objects are the natural numbers, where $n\in\N$ is viewed as the set $\{0,\dots,n\}$, with morphisms the non-decreasing maps. A \emph{simplicial object} in a category $A$ is a contravariant functor $\Delta \to A$.\\
The cyclic category $\Sigma$ is the category whose objects are the natural numbers, where $n\in\N$ is viewed as the circle $\mathbb{S}^1$ together with the distinguished subset of the $n$-th root of unity. Then the morphisms are homotopy classes of once differentiable maps $\mathbb{S}^1\to\mathbb{S}^1$ of degree $1$ with nowhere negative derivative which map the distinguished subset to the distinguished subset. A \emph{cyclic object} in a category $A$ is a contravariant functor $\Sigma\to A$.\\
There are canonical morphisms in $\Delta$, namely the \emph{face maps}. The $i$-th face $\delta_i:[n-1]\to [n]$ is defined as the injection $\{0,\dots,n-1\}\to\{0,\dots,n\}$ with $i$ not in the image. The associated morphism $M_{n}\to M_{n-1}$ in $A$ is called the \emph{$i$-th face map} and is denoted by $d_i$. Similarly, there are \emph{degeneracy maps} $\sigma_i:M_n\to M_{n+1}$ associated to the surjection $\{0,\dots,n+1\}\to\{0,\dots,n\}$ which sends both $i$ and $i+1$ to $i$. The face and degeneracy maps satisfy a certain set of equations such as $d_id_j=d_{j-1}d_i$, such that they yield a presentation of the categories $\Delta$ resp. $\Sigma$, so that the simplicial objects in a category are exactly those sequences of objects together with morphism $d_i$ and $\sigma_i$ that satisfy the latter.\\
There is a canonical inclusion functor $\Delta\to\Sigma$ giving rise to a forgetful functor from the cyclic objects in $A$ to the simplicial ones. A cyclic object is then a simplicial one together with the additional structure of \emph{cyclic operators} $t_n:M_n\to M_n$ corresponding to the rotation morphisms in the category $\Sigma$.\\
A simplicial set (module...) is a simplicial object in the category of sets (modules...).\\
A complex algebra $B$ gives rise to a cyclic vector space by virtue of the following definitions.
\begin{align*}
M_n&=B^{\otimes(n+1)}\\
d_i(b_0\otimes\dots\otimes b_n)&=(b_0\otimes\dots\otimes b_ib_{i+1}\otimes\dots\otimes b_n)\text{ if }0\le i\le n-1\\
d_n(b_0\otimes\dots\otimes b_n)&=(b_n b_0\otimes\dots\otimes b_{n-1})\\
\sigma_i(b_0\otimes\dots\otimes b_n)&=(b_0\otimes\dots\otimes b_j\otimes 1\otimes b_{j+1}\otimes\dots\otimes b_n)
\end{align*}
(It suffices to suppose that $B$ is unital since, similarly to $\K$-theory, the natural definition of cyclic homology has the property that on non-unital algebras it is the reduced one of the unitization, see \cite{Loday}, in other words cyclic homology is a theory with compact supports).\\
A simplicial module gives rise to a canonical chain complex in the following way. Let $d_i$ denote the $i$-th face operator. Then $d=\sum_{i=0}^n d_i$ makes the sequence of vector spaces a chain complex.\\
Consider the canonical functor from the category of the sets to the category of modules associating to the set $X$ the free module with base $X$. Composition with this functor defines a functor from the category of simplicial sets to the category of simplicial modules. The \emph{simplicial homology} of the simplicial set is defined as the homology of the chain complex associated to the free simplicial module on the simplicial set.\\
A cyclic module gives rise to several canonical chain complexes in the following ways.\\
First, there is the forgetful functor from the category of cyclic modules to the category of simplicial modules. Hence, there is the chain complex associated to the underlying simplicial module. Its homology is called the \emph{Hochschild} homology of the cyclic module.\\
Second the presence of the cyclic operators endows the Hochschild complex with another ascending differential
$$
B:=(1-t_{n+1})t_{n+1}\varsigma_nN
$$
where $N=\sum_{i=0}^n t^i$ is the so-called norm operator. One has the crucial identity $bB+Bb=0$.\\
The operators $b,t$ and $N$ can be used to define a bicomplex, called the \emph{cyclic} or \emph{Connes-Tsygan}-bicomplex and denoted $\CC(A)$:
$$
\xymatrix{
\dots\ar[d]&\dots\ar[d]&\dots\ar[d]&\\
A_2\ar[d]^b&A_2\ar[l]^{1-t}\ar[d]^{-b'}&A_2\ar[l]^N\ar[d]^b&\dots\ar[l]^{1-t}\\
A_1\ar[d]^b&A_1\ar[l]^{1-t}\ar[d]^{-b'}&A_1\ar[l]^N\ar[d]^b&\dots\ar[l]^{1-t}\\
A_0&A_0\ar[l]^{1-t}&A_0\ar[l]^N&\dots\ar[l]^{1-t}
}
$$
If the cyclic module comes from a unital algebra, then the odd columns are contractible, and therefore $\CC()A$ can be shrinked, leading to an example of
\begin{defn}
A \emph{mixed complex} $A$ is a sequence $A_n, n\ge 0$ of modules together with anticommuting differentials in both directions, i. e. with a map $b$ of degree $-1$ and a map $B$ of degree $+1$ such that
$$
b^2=B^2=bB+Bb=0.
$$
Any cyclic module gives rise to a mixed complex.\\
Let $(A,b,B)$ and $(B,b',B')$ be mixed complexes. Then the product mixed complex $A\otimes B$ is defined to be the mixed complex $(\Tot A\otimes B,1\otimes b'+b'\otimes 1,1\otimes B'+B'\otimes 1)$.
\end{defn}
 $A\otimes B$ is an example of a mixed complex that is not necessarily associated to a cyclic module since there is no cyclic operator on $A_n\otimes B_m$ if $n\ne m$.\\
Now there is a canonical functor associating to a cyclic module $A$ a bicomplex $\mathcal{B}(A)$, called the $B-b$-bicomplex of $A$, having the Hochschild complex as the columns and the $B$-complex as rows:
$$
\xymatrix{
\dots\ar[d]&\dots\ar[d]&\dots\ar[d]\\
A_2\ar[d]^b&A_1\ar[d]^b\ar[l]^B&A_0\ar[l]^B\\
A_1\ar[d]^b&A_0\ar[l]^B&\\
A_0&&\\
}
$$
The fact that the $B-b$-bicomplex arises by eliminating the contractible subcomplex of the cyclic bicomplex leads to the statement that there is a quasi-isomorphism
$$
\Tot\mathcal{B}(A)\to\Tot\CC(A).
$$
The homology of these complexes is called the \emph{cyclic homology} of the cyclic module resp. mixed complex.\\
Any unital complex algebra gives rise to a cyclic module in a canonical way. This defines Hochschild and cyclic homology of algebras.\\
There are the following two archetypal examples.
\begin{exples}
\mbox{ }
\begin{enumerate}

\item
Let $V$ be a complex nonsingular affine algebraic variety. Then the Hochschild homology of $\mathcal{O}[V]$ is the complex of algebraic differential forms on $V$. The differential $B$, anticommuting with $b$, descends on Hochschild homology and induces the Cartan differential $d$ on the complex of differential forms. So it follows from an inspection of the spectral sequence associated to the $B-b$-bicomplex that the periodic cyclic homology of $\mathcal{O}[V]$ is the (algebraic) de-Rham cohomology of $V$.
\item
Let now $\Gamma$ be a discrete group, and consider the algebra
$\C\Gamma$. There is a natural notion of direct sum of cyclic modules, and the cyclic module of $\C\Gamma$ splits naturally as a direct sum . It is indexed by the conjugacy classes of $\Gamma$: $Z(\C\Gamma)_{<z>}$ is the submodule generated by those $(\gamma_0,\dots,\gamma_n)$ such that the product $\gamma_0\dots\gamma_n$ is conjugate to $z$. The summand of $<1>$ is called the \emph{homogeneous} part.
\end{enumerate}
\end{exples}
The functor ``free module'' from sets to modules gives rise to a functor from simplicial sets to simplicial modules. We are now going to look at the simplicial set underlying the cyclic module given by the homogeneous part of $Z(\C\Gamma)$. The construction is a beautiful reformulation and generalization of the process of manufacturing from a given discrete group a cell complex with the given group as fundamental group.\\
There is a functor from the category of small categories (i. e. the categories such that the collection of objects is a set in the set-theoretic sense) to the category of simplicial sets, called the \emph{nerve} of the category and defined as follows. The set of $0$-dimensional simplices is the set of objects. The set of $1$-dimensional simplices is the set of morphisms. The set of $n$-simplices is the set of $n$-fold compositions $h_0\circ\dots\circ h_{n-1}$ of composable morphisms $h_i$ in a category $C$. The faces of $h_0\circ\dots\circ h_{n-1}$ are $h_1\circ\dots\circ h_{n-1}$, $h_0\circ\dots\circ h_{n-2}$ and $h_1\circ\dots\circ h_{i-2}\circ (h_{i-1}\circ h_{i})\circ h_{i+1}\circ\dots\circ h_n$ for $i=1,\dots,n-1$. The degeneracies of $h_0\circ\dots\circ h_{n-1}$ are $1_{\text{target }h_0}\circ h_0\circ\dots\circ h_{n-1}$ through $h_0\circ\dots\circ h_{n-1}\circ 1_{\text{source }h_{n-1}}$. The definition can be restated in a more elegant way as follows. Let $\mathcal{C}$ be the category of small categories. There is a functor $i:\Delta\to\mathcal{C}$, and the functor Nerve from $\mathcal{C}$ to the category of simplicial sets is given by dualizing the composition
$$
\Delta^{\text{op}}\times\mathcal{C}\xrightarrow{i\times 1}\mathcal{C}^{\text{op}}\times\mathcal{C}\xrightarrow{\Hom}\text{Sets}.
$$
Denote the nerve of the category $\Gamma$, viewed as a category with one object and the group elements as morphisms, by $B.\Gamma$. Thus, $B.\Gamma^0$ is a point, and the set of $n$-simplices is the set of $n$-element ordered sets $(\gamma_0,\dots,\gamma_{n-1})$ of arbitrary group elements. Now there is the beautiful and easily verified fact that $B.\Gamma$ is isomorphic to the simplicial set underlying the simplicial module underlying the cyclic module given by the homogeneous part of the cyclic module associated to the algebra $\C\Gamma$. In fact, the bijection
$$
(\gamma_0,\dots,\gamma_{n-1})\mapsto (\gamma_{n-1}^{-1}\dots \gamma_0^{-1},\gamma_0,\dots,\gamma_{n-1})
$$
is compatible with the face and degeneracy maps in such a way that deleting the first entry corresponds to multiplying the first with the second and so on. This fact already determines the homogeneous Hochschild and cyclic homology of $\C\Gamma$. Namely, consider the functor ``geometric realization'' from simplicial sets to (compactly generated Hausdorff) topological spaces. Then it belongs to the coherence of the definitions that the singular homology of the geometric realization of a simplicial set is its simplicial homology, in other words the Hochschild homology of the free simplicial module on the simplicial set. Denote the (homotopy type of) the geometric realization of $B.\Gamma$ by $\BG$. Thus, the homogeneous Hochschild homology $\HH(\C\Gamma)_{<1>}$ of $\C\Gamma$ is equal to the singular homology $\Homol_*(\BG;\C)$ of the geometric realization of the simplicial space given by the nerve of the category $\Gamma$. A possible definition of \emph{group homology} of $\Gamma$ with coefficients in the constant $\Gamma$-module $\C$ is $\Homol_*(BG;\C)$. It is not difficult to show that, more generally, the Hochschild homology is group homology with coefficients in the adjoint representation $\C\Gamma\to\End_\C(\C\Gamma)$ defined by $\gamma\mapsto(\gamma'\mapsto \gamma\gamma'\gamma^{-1})$. (Group homology with non-trivial coefficients may also be defined topologically using coefficients in locally constant sheaves.) $\C\Gamma$ contains the constant submodule $\C$ corresponding to the neutral element. The homogeneous part of the Hochschild homology of $\C\Gamma$ is then the subspace corresponding to this submodule.\\
Just as Hochschild homology of a simplicial set is the singular homology of its geometric realization, the cyclic homology of a cyclic set can be read off from its geometric realization given the additional data of an $\mathbb{S}^1$-action that is defined on it, e. g., by \cite{Loday}, section 7.1, based on the action of the finite cyclic group on the $n$-simplices. One can show that the singular $\mathbb{S}^1$-equivariant homology of the geometric realization of a cyclic set, i. e. of its underlying simplicial set together with the $\mathbb{S}^1$-action, is the cyclic homology of the free cyclic module on the cyclic set. $\mathbb{S}^1$-equivariant homology of an $\mathbb{S}^1$-space $X$ means the ordinary homology of the homotopy quotient $\classb \mathbb{S}^1\times_{\mathbb{S}^1} X$. This fact determines the homogeneous cyclic homology $\HC_*(\C\Gamma)_{<1>}$ of $\C\Gamma$. Namely, the obtained action of $\mathbb{S}^1$ on $\BG$ is then trivial, and thus $\HC_*(\C\Gamma)_{<1>}=\Homol_*(\BG\times\classs\mathbb{S}^1)=\Homol_*(\BG)\otimes\Homol_*(\classs\mathbb{S}^1)=\Homol_*(\BG)\otimes\Homol_*(\mathbb{CP}^\infty)=\Homol_*(\BG)\otimes \C[[\varsigma]]$ where $\C[[\varsigma]]$ is the algebra of formal power series in the indeterminate $X$. $\varsigma$ corresponds to the periodicity $S$-operator.\\
However, there is a another means of computing $\HC(\C\Gamma)_{<1>}$ which fits much better with the generalization to crossed products. In fact, stretching the topological image even further, there is another canonical simplicial space associated to $\Gamma$, namely the universal covering simplicial space, denoted $E.\Gamma$. It is the nerve of the category which has the elements of $\Gamma$ as objects and exactly one morphism for any ordered pair of objects. So $E.\Gamma$ has as $n$-simplices ordered $n$-tuples $(\gamma_0,\dots,\gamma_n)$ of group elements, with structure maps given by
\begin{align*}
d_i(\gamma_0,\dots,\gamma_n)&=(\gamma_0,\dots,\widehat{\gamma_i},\dots,\gamma_n)\\
s_i(\gamma_0,\dots,\gamma_n)&=(\gamma_0,\dots,\gamma_i,1,\gamma_{i+1},\dots,\gamma_n)
\end{align*}
The notations $E.\Gamma$ and $B.\Gamma$ are now justified by the following clues. First, $\EG$, the geometric realization is \emph{contractible} because $E.\Gamma$ has an initial object. Second, the formula $(\gamma_0,\dots,\gamma_n)\gamma=(\gamma^{-1}\gamma_0,\dots,\gamma^{-1}\gamma_n)$ defines a simplicial right action, i. e. a right action compatible with all structure maps, of $\Gamma$ on $E.\Gamma$ and thus also on $\EG$. The group $\Gamma$ acts freely on $E.\Gamma$, and therefore so does it on $\EG$. Furthermore, any free simplicial action on geometrical realizations is proper. So the name $\classb.\Gamma$ makes sense. Third, the quotient $E.\Gamma /\Gamma$ of the action of $\Gamma$ on $E.\Gamma$ is isomorphic by the map
$$
(\gamma_0,\dots,\gamma_n)\mapsto (\gamma_n\inv\gamma_0,\gamma_o\inv\gamma_1,\dots,\gamma_{n-1}\inv\gamma_n)
$$
to the simplicial set underlying $Z(\C\Gamma)_{<1>}$, the homogeneous part of the cyclic module associated to the algebra $\C\Gamma$ having as $n$-simplices expressions $(\gamma_0,\dots,\allowbreak\gamma_n)$ such that the product equals $1$. The latter simplicial set in turn is isomorphic to $B.\Gamma$, as we have checked. The name $B.\Gamma$ is therefore also justified because geometric realization commutes with taking quotients.\\
The Hochschild complex $(\C[E.\Gamma],b)$ associated to the simplicial
module $E.\Gamma$ plays quite a central r\^ole in what follows. It is defined by
\begin{align*}
\C[E.\Gamma]_n&=(\C\Gamma)^{\otimes n+1}, n\ge 0\\
d(\gamma_0,\dots,\gamma_n)&=\sum_{j=0}^n(-1)^j(\gamma_0,\dots,\widehat{\gamma_i},\dots,\gamma_n).
\end{align*}
We already know that its homology is that of a point, and in fact the formula $\delta(\gamma_0,\dots,\gamma_n)=(\gamma_0,\dots,\gamma_n,1)$ defines a contracting homotopy. As is $E.\Gamma$, the complex $\C[E.\Gamma]$ is equipped with a $\Gamma$-right action, and the important feature is that this action is \emph{free}. The map $\epsilon:\C[E.\Gamma]\to\C,\sum a_\gamma\gamma\mapsto\sum a_\gamma$ defines a free resolution in the category of $\C\Gamma$-modules of the constant module $\C$. This resolution is called in the literature the ``bar resolution'' (because this meaningless word is familiar in similar contexts) and denoted $C^{\text{bar}}_\bullet\Gamma$. By the definition of group homology proper, the complex $\C[E.\Gamma]\otimes_{\C\Gamma}\C$ of $\Gamma$-coinvariants of $\C[E.\Gamma]$ computes group homology $\Homol_*(\Gamma;\C)$. We have again computed the homogeneous Hochschild homology of $\C\Gamma$ because taking the chain complex of a simplicial module commutes with taking quotients in the indicated sense.\\
Now the computation of the homogeneous cyclic homology of $\C\Gamma$ that needs to take into a account that we are dealing with cyclic, not merely simplicial sets, is hardly more complicated, if the notion of hyperhomology from homological algebra is used. In order to investigate $\mathcal{B}(Z(\C\Gamma)_{<1>})$, it is convenient to replace it by the $\Gamma$-equivariantly quasi-isomorphic Connes-Tsygan-bicomplex $\CC(Z(\C\Gamma)_{<1>})$, see \cite{Loday}. By definition, $\CC(Z)$ for a cyclic module $Z$ is defined as the periodic bicomplex
$$
(Z,b)\xleftarrow{1-t}(Z,-b')\xleftarrow{N}(Z,b)\xleftarrow{1-t}\dots,
$$
where $(Z,b)$ is the Hochschild complex with $b=\sum_{i=0}^n(-1)^id_i$, $(Z,-b')$ is the complex with same entries but $b'=\sum_{j=0}^{n-1}(-1)^jd_j$, and $1-t$ and $N$ are the maps we have already introduced. $(Z,b')$ is contractible for $Z=Z(A)$ of a unital algebra $A$. So $\CC(Z(\C\Gamma)_{<1>})$ is $\Gamma$-equivariantly quasi-isomorphic to the $\Gamma$-constant periodic complex in the first quadrant with $0$ in all rows except the lowest, and
$$
\C\xleftarrow{} 0\xleftarrow{} \C \xleftarrow{} \dots
$$
in the lowest. In other words, the bicomplex $\CC(Z(\C\Gamma)_{<1>})$ is a resolution of the latter complex. Furthermore, $\CC(Z(\C\Gamma)_{<1>})=\CC(B.\Gamma)=\CC(E.\Gamma/\Gamma)=\CC(E.\Gamma)\otimes_{\C\Gamma}$, and $\Tot \CC(Z(\C\Gamma)_{<1>})=\CC(E.\Gamma)\otimes_{\C\Gamma}\C$ computes by definition the group hyperhomology 
$$
\mathbb{H}(\Gamma;\C\leftarrow 0\dots)=\Homol_*(\Gamma;\C)\oplus\Homol_{*-2}(\Gamma;\C)\oplus\dots.
$$
(Of course, here one can read off the homology of the quotient directly from the quasi-isomorphism.) So $\HC_*(\C\Gamma)_{<1>}\cong\bigoplus_{j=0}^\infty\Homol_{*-2j}(\Gamma;\C)$. Dually, $\HC^*(\C\Gamma)_{<1>}\cong\Pi_{j=0}^\infty\Homol^{*-2j}(\Gamma;\C)$. In view of the Baum-Connes-conjecture, this is what we would have expected (maybe up to the fact that this is already exhausted by the homogeneous part).\\
We will now carry out the whole discussion, following Nistor, in detail for the more general case of crossed products $A\rtimes\Gamma$. The only new feature is that we need to make use of the simplicial Eilenberg-Zilber theorem.
\subsection{Crossed Products}
\noindent Let $A$ be a complex algebra with unit together with a $\Gamma$-action.\\
The cyclic vector space (it has become natural to say just ``cyclic module'') $Z(A\alcross\Gamma)$ of the algebraic crossed product has underlying vector space
$$
(A\alcross\Gamma)^{\otimes(n+1)}=A^{\otimes(n+1)}\otimes (\C\Gamma)^{\otimes(n+1)}=A^{\otimes(n+1)}\otimes\C[\Gamma^{\times(n+1)}].
$$
Let $<z>$ be a conjugacy class of $\Gamma$. Following Nistor, we denote by $L(A,\Gamma, z)$ the subspace of the cyclic module $Z(A\alcross\Gamma)$ of the crossed product that is generated by those $(a_0\otimes\dots\otimes a_n\otimes\gamma_0\dots\otimes\gamma_n)$ such that the product $\gamma_0\gamma_1\dots\gamma_n$ is conjugate to $z$. Since the conjugacy class of a product is invariant under cyclic permutations, the maps $d_i,s_i$ all preserve $L(A,\Gamma,z)$, and so there is a natural direct sum decomposition on the level of cyclic modules
$$
Z(A\alcross\Gamma)\cong\bigoplus_{<z>\in<\Gamma>}L(A,\Gamma,z),
$$
where the sum is over the conjugacy classes of $\Gamma$. The summand corresponding to the neutral element is called the \emph{homogeneous} part, the complement is called the ``inhomogeneous'' one. (This decomposition is in analogy with the splitting of the group ring of a finite group into parts corresponding to the conjugacy classes, described by Wedderburn theory). As a consequence, the cyclic homology and cohomology split naturally over the conjugacy classes as well, and both parts become orthogonal with respect to the pairing. We are only interested in a topological description of the homogeneous part, since the class of the trace $\tau$ on $\C\Gamma$ is concentrated in the homogeneous part.\\
\subsection{The Algebraic Case}
In order to separate the algebra from the analysis as far as possible, let us first treat a simplified case, namely that of algebraic cyclic homology associated to affine algebraic varieties. In that case, there need neither be taken projective tensor products for the definition of the cyclic modules, nor need the functionals be continuous in any sense for the cohomology. So we work in a purely algebraic context.
\begin{defn}
The group hyperhomology $\mathbb{H}(\Gamma,A_\bullet)$ of a discrete group $\Gamma$ with coefficients in a complex $A_\bullet$ of $\C\Gamma$-modules is defined as the hypertor functor $\mathbb{T}\mathrm{or}_*^{\C\Gamma}(\C,A_\bullet)$. The group hypercohomology with coefficients in a $\C\Gamma$-cochain complex $B^\bullet$ is defined as the hyperext functor $\mathbb{E}\mathrm{xt}^*_{\C\Gamma}(\C,B^\bullet)$.
\end{defn}
In the following theorem, $\Omega^\bullet V$ resp. $\Omega_\bullet V$ means the complex of \emph{algebraic} differential forms resp. its dual, the complex of algebraic de-Rham currents on a complex nonsingular affine algebraic variety $V$. Normally, $\Omega^\bullet V$ is a cochain complex in the sense that its differential has degree $+1$. However, here it has to be transformed into a chain complex by first viewing it as a $\Z/2\Z$-graded complex $\Omega^\ev V\leftrightarrows \Omega^\odd V$ and then interpreting the latter as a periodic chain complex, where $\Omega^\ev$ sits in even degree. A similar remark applies to the dual complex.
\begin{thm}
Let $V$ be a complex nonsingular affine algebraic variety, $\mathcal{O}[V]$ its coordinate ring, and $\Gamma$ be a discrete group acting on $V$ by diffeomorphisms.
\begin{enumerate}
\item
The homogeneous part of the algebraic periodic cyclic homology
$$
\HP_*\bigl(\mathcal{O}[V]\alcross\Gamma\bigr)_{<1>},\ *=0,1
$$
of the crossed product of the coordinate ring of $V$ by $\Gamma$ is naturally isomorphic (in both $V$ and $\Gamma$) to
$$
\mathbb{H}_{*}(\Gamma;\ \Omega^\ev V\leftrightarrows\Omega^\odd V),
$$
the group hyperhomology of $\Gamma$ with coefficients in the complex of differential forms on V.
\item
The homogeneous part of the cyclic cohomology of the crossed product is the group hypercohomology with coefficients in the de-Rham currents 
$$
\HP^*\bigl(\mathcal{O}[V]\alcross\Gamma\bigr)_{<1>}\cong\mathbb{H}^{*}(\Gamma;\ (\Omega^\ev V\leftrightarrows\Omega^\odd V)^*).
$$
\end{enumerate}
\end{thm}
In order to identify these hyperhomologies topologically, i. e. to find spaces whose homologies they are in a suitable functorial way,  $\Gamma$ and $V$ should enter with \emph{equal} functoriality. Later on, it will therefore be necessary to take Poincar\'e duality into account. However, it will be necessary to first work out the analytical analogue of the theorem.\\
We shall need some definitions to formulate the proof.
\begin{proof}
\begin{enumerate}
\item
Let $\Psi$ denote $L(A,\Gamma,1)$, the homogeneous part of the cyclic module as defined above.\\
Let $\widetilde{\Psi}$ be the ``covering space'' of $\Psi$ denoted $\widetilde{L}(A,\Gamma,1)$ in \cite{Nistor}, $\widetilde{\Psi}^n=(A\alcross\Gamma)^{\otimes(n+1)}$. $\widetilde{\Psi}$ is made a cyclic module by virtue of the following definitions:
\begin{align*}
d_i(a_0,\dots,a_n,\gamma_0,\dots,\gamma_n)=&(a_0,\dots,a_ia_{i+1},\dots,a_n,\gamma_0,\dots,\widehat{\gamma_i},\dots,\gamma_n)\\
&\mbox{for }0\le i\le n-1\\
d_n(a_0,\dots,a_n,\gamma_0,\dots,\gamma_n)=&(a_n a_0,a_1,\dots,a_{n-1},\gamma_0,\dots,\gamma_{n-1})\\
\sigma_i(a_0,\dots,a_n,\gamma_0,\dots,\gamma_n)=&(a_0,\dots,a_i,1,a_{i+1},\dots, a_n,\\
&\hspace{2em}\gamma_0,\dots,\gamma_i,\gamma_i,\gamma_{i+1},\dots,\gamma_n)\\
&\mbox{for }0\le i\le n\\
t(a_0,\dots,a_n,\gamma_0,\dots,\gamma_n)=&(-1)^n(a_n,a_1,\dots,a_{n-1},\gamma_n,\dots,\gamma_{n-1}).
\end{align*}
$\widetilde{\Psi}$ possesses a right $\Gamma$-action through cyclic module maps
$$
(a_0,\dots,a_n,\gamma_0,\dots,\gamma_n)\gamma=(\gamma\inv.a_0,\dots,\gamma\inv.a_n,\gamma\inv\gamma_0,\dots,\gamma\inv\gamma_n).
$$
One readily checks that the quotient $\widetilde{\Psi}/\Gamma$ identifies with $\Psi=L(A,\Gamma,1)$ by the $\Gamma$-invariant map
\begin{align*}
p:\widetilde{\Psi}\to\Psi:\quad&(a_0,\dots,a_n,\gamma_0,\dots,\gamma_n)\mapsto\\
&(\gamma_n\inv.a_0,\gamma_0\inv.a_1,\dots,\gamma_{n-1}\inv.a_n,\gamma_n\inv\gamma_0,\gamma_0\inv\gamma_1,\dots,\gamma_{n-1}\inv\gamma_n).
\end{align*}
Let $\bigl(Z(A),b_A,B_A\bigr)$ denote the mixed complex associated to $A$ (consisting of algebraic, not projective tensor products), and let $\bigl(\C[\Enerve],b,B\bigr)$ denote the mixed complex associated to the cyclic module $\C[\Enerve]$. Here, $\Enerve$ is the nerve of the category with objects the group elements as explained above. $\bigl(\C[\Enerve], b\bigr)$ is the contractible bar resolution of $\Gamma$, and $\Enerve/\Gamma=B.\Gamma,$ the classifying space of the category $\Gamma$.\\
Then the cyclic module $\widetilde{\Psi}$ is the product (in the sense of \cite{Loday}) of cyclic modules $Z(A)\times\C[\Enerve]$, i. e.
$$
\widetilde{\Psi}^k=Z^k(A)\otimes \C[\Enerve]^k=A^{\otimes(k+1)}\otimes(\C\Gamma)^{\otimes(k+1)}
$$
with cyclic structure defined by
$$
d_i=d^A_i\otimes d^\mathrm{Bar}_i,\ s_j=s_j^A\otimes s^\mathrm{Bar}_j,\ t_n=(-1)^nt^A_n\otimes t^\mathrm{Bar}_n.
$$
All tensor products are algebraic.\\
Let us first compute the ordinary cyclic homology and then show that $S$ stabilizes for large $n$.\\
\begin{align*}
\HC_n\bigl(A\alcross\Gamma\bigr)_{<1>}
&=\HC_n\bigl(\Psi\bigr)\\
&=\Homol_n\bigl(\Tot\mathcal{B}(\Psi)\bigr),
\intertext{$\mathcal{B}$ being, of course, the $(b, B)$-bicomplex of a mixed complex,}
&=\Homol_n\bigl(\Tot\mathcal{B}(\widetilde{\Psi}/\Gamma)\bigr)\\
&=\Homol_n\bigl(\Tot\mathcal{B}(\widetilde{\Psi})/\Gamma\bigr)\\
&=\Homol_n\bigl(\Tot\mathcal{B}(Z(A)\times \C[\Enerve])/\Gamma\bigr)
\intertext{$\Tot\mathcal{B}(Z(A)\times \C[\Enerve])$ is a complex of \emph{free} $\C\Gamma$-modules, so by \cite{Brown} the latter group is canonically isomorphic to group hyperhomology with coefficients in this complex.}
&=\mathbb{H}_n\bigl(\Gamma;\ \Tot\mathcal{B}\bigl(Z(A)\times \C[\Enerve]\bigr)\bigr)
\end{align*}
We now have to make use of a weakening of the notion of morphism of mixed complexes. Any mixed complex $C$ has associated $B-b$-bicomplex $\mathcal{B}C$ and a periodicity operator $S:\Tot\mathcal{B}C\to\Tot\mathcal{B}C[2]$ given by dividing by the subcomplex of the first line. An $S$-operator of two mixed complexes is a morphism of complexes $\Tot\mathcal{B}C\to\Tot\mathcal{B}C'$ that commutes with $S$ (\cite{Loday}, 2.5.14). Now there is an \emph{Eilenberg-Zilber theorem for cyclic homology} (see \cite{Loday}, 4.3.8), which states that there is a shuffle map
$$
\Tot \mathcal{B}(A\otimes B)\to\Tot\mathcal{B}(A\times B)
$$
which provides an $S$-morphism between the $\times$-product and the algebraic tensor product of two mixed complexes. This $S$-morphism restricts to the shuffle map on the Hochschild complex, where it induces an isomorphism with inverse the Alexander-Whitney-map, see \cite{Weibel}, page 277, and \cite{Loday}, 4.3.8. So this $S$-morphism induces an isomorphism on Hochschild homology. By a standard inductive argument involving the $SBI$-sequence, it induces also an isomorphism on the homology of the $\mathcal{B}$-complex. This $S$-morphism is a $\Gamma$-equivariant quasi-isomorphism, so (by \cite{Weibel}, 5.7.7/2) it induces also an isomorphism on group hyperhomology,
$$
\mathbb{H}_n\bigl(\Gamma;\ \Tot\mathcal{B}\bigl(Z(A)\otimes \C[\Enerve]\bigr)\bigr)\xrightarrow{\cong}\mathbb{H}_n\bigl(\Gamma;\ \Tot\mathcal{B}\bigl(Z(A)\times\C[\Enerve]\bigr)\bigr).
$$
Now there is the slant product $\backslash \epsilon$ (called cap product in \cite{Nistor}) with the generator
$$
\epsilon\colon  \C\Gamma\to\C,\ \sum a_\gamma\gamma\mapsto\sum a_\gamma
$$
of $\Homol_*(\Hom(\C[E.\Gamma],\C))$ (the augmentation of the bar resolution):
\begin{align*}
&\backslash \epsilon\colon \Tot Z(A)\otimes\C[\Enerve]\to Z(A),\\
&\backslash\epsilon\bigl((a_0,\dots,a_p)\otimes(\gamma_0,\dots,\gamma_q)\bigr)=0\ \mathrm{if}\ q\ge 1\\
&\backslash \epsilon\bigl((a_0,\dots,a_p)\otimes \gamma_0\bigr)=(a_0,\dots,a_n),
\end{align*}
which again is a $\Gamma-$equivariant quasi-isomorphism. So again by the $SBI$-se-quence it is also a quasi-isomorphism on the $\mathcal{B}$-bicomplex, and therefore also an isomorphism on hyperhomology:
$$
\mathbb{H}_n\bigl(\Gamma;\ \Tot\mathcal{B}\bigl(Z(A)\otimes\C[\Enerve]\bigr)\bigr)=\mathbb{H}_n\bigl(\Gamma;\ \Tot\mathcal{B}\bigl(Z(A)\bigr)\bigr)=\mathbb{H}_n\bigl(\Gamma; \Tot\mathcal{B}(A)\bigr).
$$
So the (homogeneous part) of the cyclic homology of the crossed product is simply the group hypercohomology with coefficients in the $B-b$- or cyclic bicomplex of $A$.\\
Note that it would be wrong to state an isomorphism of $\HC_*(\Gamma\allcross A)$ with the expression $\mathbb{H}_*(\Gamma;\HC_\bullet(A))$ although $\HC_\bullet(A)$ is, with the help of the axiom of choice, algebraically quasi-isomorphic to $\Tot\mathcal{B}(A)$. The reason is that the action of $\Gamma$ on $\HC_\bullet(A)$ may be trivial without the cyclic homology of the crossed product being trivial. The problem is that this quasi-isomorphism is not $\Gamma$-equivariant. In fact, one can obtain a counterexample whenever the hyperhomology spectral sequence does not degenerate at the $E_2$-term. For instance, take for $\Gamma$ the fundamental group of the Riemannian surface of genus 2 acting on the boundary $\mathbb{S}^1$ of the universal cover, as described in the K-theory part (see also \ref{caution}).\\
Now let us consider $\HP$ rather than $\HC$. Let us assume that the $S$-operator stabilizes for the algebra $A$. This means that there is an integer $n$ such that for all $m\ge n$ the operator $S:\HC_m(A)\to\HC_{m-2}(A)$ is an isomorphism, or equivalently that $A$ is of finite Hochschild dimension. This assumption will hold in all our applications. The $S$ operator $\Tot\mathcal{B}(A)\to\Tot\mathcal{B}(A)[2]$ is $\Gamma$-equivariant and induces therefore an operator $\mathbb{H}_n(\Gamma; \mathcal{B}(A))\to\mathbb{H}_n(\Gamma; \mathcal{B}(A)[2])$. One checks that all operations we performed on the bicomplex $\mathcal{B}(\Psi)$ commute with $S$. As a consequence, the operator between the group hyperhomology induced by the $S$-operator in the coefficients corresponds to the $S$-operator on cyclic homology under the isomorphism we have established. So the $S$-operator stabilizes also for the homogeneous part of the algebra $\Gamma\allcross A$. It follows that we obtain the homogeneous periodic cyclic homology as
$$
\HP_*(\Gamma\allcross A)_{<1>}\cong\mathbb{H}_*(\Gamma; \Totpr \CCP Z(A)),
$$
the group hyperhomology with coefficients in the completed periodic cyclic complex. $\Tot^\Pi$ means the product total complex of an unbounded bicomplex, and $\CCP$ means the periodic bicomplex, which is the complex $\CC$ continued infinitely to the left. In the definition of hyperhomology some care has to be taken: $\Totpr$ is an unbounded complex itself, and its resolution is supported in the first and second complex. We define the hyperhomology now using the direct sum total complex of the resolution.\\
However, the latter difficulty readily disappears when we come to consider the algebra $A=\mathcal{O}[V]$. Let us compute the periodic homology in this case again, without using the $S$-operator.\\
Let $\btdR(V)$ be the bicomplex of truncated de-Rham complexes of $V$:
$$
\xymatrix{
  {\dots}        \ar[d]     & {\dots}         \ar[d]                & {\dots}         \ar[d]  & {\dots}\\
  {\Omega^2V} \ar[d]^{0} & {\Omega^1V} \ar[l]^{d} \ar[d]^{0} & {\Omega^0V} \ar[l]^{d}\\
  {\Omega^1V} \ar[d]^{0} & {\Omega^0V} \ar[l]^{d}\\
  {\Omega^0V} 
}
$$
Note that $\btdR$ is a chain complex instead of a cochain complex.\\
Now let $A$ be the smooth commutative algebra $\mathcal{O}[V]$. There is a $\Gamma-$equivariant quasi-isomorphism $\mathcal{B}\bigl(Z(\mathcal{O}[V])\bigr)\cong \btdR(V)$ because we're in the algebraic framework on an affine algebraic variety. This quasi-isomorphism passes to the completed periodic bicomplex of $A$ to a quasi-isomorphism with a ``completed periodic de-Rham bicomplex''.  But $\Omega^n V=0$ for $n>\dim V$, so the completion process does nothing, and we are left with the periodic complex or bicomplex $\Omega^\ev V\leftrightarrows\Omega^\odd V$. This proves the assertion that
$$
\HP_*\bigl(\mathcal{O}[V]\alcross\Gamma\bigr)_{<1>}\cong\mathbb{H}_{*}(\Gamma;\ \Omega^\ev V\leftrightarrows\Omega^\odd V).
$$
This is defined for $*\in\Z$ but it is $\Z/2\Z$-periodic.\\
We have reached a point where a certain confusion that will later on tend to arise takes its origin: The complex $\Omega^\ev V\leftrightarrows\Omega^\odd V$ has to be understood as a \emph{chain complex}, even though it is normal and will be necessary to view it as a \emph{cochain complex}.
\item
For the computation of the cyclic and periodic cyclic cohomology, we could in the algebraic framework simply refer to the universal coefficient theorem, which states that the cohomology vector spaces are but the algebraic duals of the homology vector spaces. Nevertheless, we shall carry out an explicit computation, because later on there will be no such theorem for the theory that takes into account topologies on the crossed products.
\begin{align*}
\HC^n\bigl(\mathcal{O}[V]\alcross\Gamma\bigr)
= & \Homol_n\bigl(\Tot \Hom(\mathcal{B}(\Psi),\C)\bigr)\\
= & \Homol_n\bigl(\Tot \Hom(\mathcal{B}(\widetilde{\Psi})/\Gamma,\C)\bigr)\\
= & \Homol_n\bigl(\Tot \Hom_{\C\Gamma}(\mathcal{B}(\widetilde{\Psi}),\C)\bigr)\\
= & \Homol_n\bigl((\Tot \Hom(\mathcal{B}(\widetilde{\Psi}),\C))^\Gamma\bigr)\\
= & \mathbb{H}^n\bigl(\Gamma;\ \Tot\Hom(\mathcal{B}(\widetilde{\Psi}),\C)\bigr),
\end{align*}
where $\Gamma$ acts on a morphism as $\gamma.f=\gamma.f(\gamma\inv.-)$. The transpose of a $\Gamma$-equivariant quasi-isomorphism is also a $\Gamma$-equivariant quasi-isomorphism. We may therefore apply exactly the same steps as in the computation of the homology, replacing group hyperhomology with hypercohomology and the de-Rham complex with its (algebraic) dual. This yields
$$
\HP^*(\Gamma\allcross A)_{<1>}\cong\mathbb{H}^*(\Gamma; \Totsum \Hom (\CCP Z(A);\C)),
$$
if the $S$-operator stabilizes for $A$, and in the case of $A=\mathcal{O}[V]$ we obtain
$$\HP^*\bigl(\mathcal{O}[V]\alcross\Gamma\bigr)_{<1>}\ \cong\
\mathbb{H}^{*}\bigl(\Gamma;\ (\Omega_\ev\leftrightarrows\Omega_\odd V) \bigr).$$
\end{enumerate}
\end{proof}
\begin{cor}{(of the proof.)}\label{spectralsequence}
There are spectral sequences computing the (graded pieces of) cyclic homology in high degree:
$$
E^2_{pq}=\Homol_p\bigl(\Gamma;\ \HC_q(A)\bigr)\implies \Gr_p\ \HC_{p+q}\bigl(A\alcross \Gamma\bigr)_{<1>},
$$
and cohomology
$$
E_2^{pq}=\Homol^p\bigl(\Gamma;\ \HC^q(A)\bigr)\implies\Gr_p\ \HC^{p+q}\bigl(\Gamma\alcross A\bigr)_{<1>}.
$$
There is a canonical choice for them, and these choices coincide with the spectral sequences in \cite{Nistor}*{theorem 2.6}.
\end{cor}
\begin{proof}
These are the hyperhomology spectral sequences for
$$
\mathbb{H}_{p+q}(\Gamma; \Tot\mathcal{B}(A))
$$
resp.
$$
\mathbb{H}^{p+q}(\Gamma; \Tot \Hom(\mathcal{B}(A);\C)).
$$
They coincide with Nistor's since ours as well as his are the spectral sequences associated to the bicomplex $\C[E.\Gamma]\otimes_{\C\Gamma}\Tot\mathcal{B}(A)$ resp. its dual.

\end{proof}
\subsection{Bornological Algebras}
We now want to replace the variety $V$ with a smooth compact manifold, and the algebra $\mathcal{O}[V]$ with the locally convex $\CINF(V)$. If we tried to calculate the algebraic cyclic homology or cohomology of the crossed product $\Gamma\allcross\CINF(V)$, we would have the same problem as the one appearing in the fact that the algebraic periodic cyclic homology of $\CINF(V)$ is \emph{not} equal to $V$'s de Rham-cohomology. But when we consider in the complex that defines $\HP^*(V)$ only those functionals $(\CINF(V))^{\otimes n}\to\C$ which are continuous with respect to the topology on $(\CINF(V))^{\otimes n}$ given as a subspace of $(\CINF(V))^{\protens n}=\CINF(V^{\times n})$, we obtain a similar cyclic complex whose homology is commonly also called cyclic homology (this theory is \emph{not} to be confused with analytic cyclic homology), and which does satisfy $\HP_*(\CINF(V))\cong\oplus_i\HomoldR^i(V)$. The algebraic periodic cyclic cohomology would be much too large. Analogously, one might hope for a reasonable locally convex topology on the crossed product $\Gamma\allcross\CINF(V)$, which would, then, have to be completed. However, this is not immediate, at least if one wants to avoid to use a word metric on $\Gamma$. It turns out that it is a far better idea to equip $\Gamma\allcross\CINF(V)$ with a bornology instead of a topology. A bornology is a collection of subsets, called \emph{small} sets, that forms a covering and which is stable with respect to finite unions, containment and the algebraic operations. These axioms are satisfied by the collection of bounded subsets of a locally convex topological vector space, or by the system of relatively compact sets, or the system of equicontinuous sets, etc. See \cite{MEYERKTHEORY} and \cite{Meyerthesis} for good introductions to bornological algebras, as well as \cite{Hogbe} for a thorough introduction. A bounded map between bornological vector spaces is a linear maps which maps small sets into small sets. Thus there is the category of bornological vector spaces. A bornological vector space is said to be \emph{complete} if each small set is contained in another small set which is the unit ball of a Banach norm on its linear hull. A bilinear map $A\times B\to C$ from the product of two bornological spaces $A$ and $B$ to a bornological space $C$ is called bounded if it maps sets of the form $X\times Y$ to small sets, where $X$ is small in $A$ and $Y$ is small in $B$. A bornological algebra is a complete bornological vector space which at the same time is an algebra such that the multiplication $A\times A\to A$ is a bounded bilinear map. We have thus defined the category of complete bornological algebras,and this is the category we will work with for the rest of this thesis.\\
There is the notion of completed bornological tensor product on the category of complete bornological algebras (\cite{Meyerthesis}*{2.2.3}). It has the remarkable property that it reduces to the algebraic tensor product when one of the algebras has the fine bornology, and that it reduces to the projective tensor product $\protens$ of Fr\'echet algebras, when both algebras are Fr\'echet algebras endowed with the precompact bornology (\cite{Meyerthesis}*{theorem 2.29}).\\
In the following, we shall denote by $\altens$ the algebraic tensor product, which is also the complete bornological one when one of the two algebras has the fine bornology.\\
We will make use of two cyclic theories on the category of complete bornological algebras. The first one is ordinary cyclic homology associated to bornological algebras, the second one is analytic cyclic homology, as developped in the body of \cite{Meyerthesis}. In this chapter, we only use the former. Cyclic cohomology is defined simply as the homology of the subcomplex of the ordinary cyclic or periodic cyclic cocomplex consisting of \emph{bounded} functionals. Cyclic homology is defined as the homology of the cyclic complex where we use completed bornological tensor products $A^{\protens n}$ instead of algebraic ones $A^{\altens n}$. A similar definition is made for Hochschild homology and cohomology. In the appendix of \cite{Meyerthesis} can be found a proof of the $SBI$-sequence in this framework.\\
The category of complete bornological vector spaces has the advantage that even algebras without any additional structure are contained as a full subcategory. In fact, any algebra is a complete bornological algebra if all subsets contained in a bounded subset of a finite-dimensional subspace are taken as small sets. This bornology is called the \emph{fine} bornology. This bornology is complete because any finite dimensional normed vector space is a Banach space. Furthermore, the category of Fr\'echet algebras is also contained as a full subalgebra because the collection of all precompact sets of a Frechet algebra can be taken as the collection of small sets.\\
Let $V$ be a smooth manifold, not necessarily compact or orientable, and $\Gamma$ be a discrete group acting on $V$ by diffeomorphisms. Endow the crossed product algebra with a bornology by decreeing that sets that correspond under the vector space isomorphism $\Gamma\allcross A\cong\C\Gamma\altens A$ with the complete bornological tensor product of the fine algebra $\C\Gamma$ with the Frechet algebra $\CCINF(V)$ endowed with the precompact bornology (which, here, is the same as the bounded bornology) to small sets.\\
Another way of defining the bornology on the crossed product would be to view $\Gamma\times V$ as a smooth groupoid $\Gamma\ltimes V$. The convolution product on $\Gamma\allcross\CCINF(V)=\CCINF(\Gamma\ltimes V)$ is then separately continuous, so it is a complete bornological algebra when endowed with the precompact (which, here, is equal to the bounded) bornology (see \cite{Meyerthesis}*{4.4}). The convolution algebra is an example of an \emph{LF-algebra}, see \cite{Meyerthesis}*{example 2.7} for a definition and explanation of the notion of LF-agebra.\\
Let now $\Omega^\bullet V$ resp. $\Omega_\bullet V$ stand for the complexified $\mathcal{C}^\infty$-de Rham complex on $V$ resp. its continuous complex valued dual.
\begin{thm}
There are isomorphisms which are natural in $\Gamma$ and $V$ for the homogeneous parts of periodic cyclic homology and cohomology
\begin{align*}
\HP_*\bigl(\CINF(V)\alcross\Gamma\bigr)_{<1>}\ &\cong\ \mathbb{H}_{*}(\Gamma;\ (\Omega^\ev V\leftrightarrows\Omega^\odd V)),\quad *=0,1\\
\HP^*\bigl(\CINF(V)\alcross\Gamma\bigr)_{<1>}\ &\cong\ \mathbb{H}^{*}(\Gamma;\ (\Omega_\ev V\leftrightarrows \Omega_\odd V)),\quad *=0,1.
\end{align*}
\end{thm}
\begin{proof}
It is advantageous to follow the lines of the proof in the preceding section in the following way. All manipulations that led there to the result $\HC_*(\Gamma\allcross A)\cong\mathbb{H}_*(\Gamma;\HC_\bullet(A))$ came from explicit manipulations of the complexes calculating these groups. Therefore the proof will remain valid if we keep track of the differences between the complexes used in the preceding section and the ones used in this section. This means that we shall not define bornological cyclic modules, but work directly with the associated complexes.\\

Consider the typical entry $(\Gamma\allcross\CCINF(V))^{\protens n+1}$ in the cyclic bicmplex of the algebra $\Gamma\allcross\CCINF(V)$ as defined above. We have by definition
\begin{align*}
(\Gamma\allcross\CCINF(V))^{\protens(n+1)}=(\C\Gamma\otimes\CCINF(V))^{\protens (n+1)} 
\intertext{as a bornological vector space. Here the ordinary tensor product sign is at the same time the algebraic and the completed bornological tensor product, because $\C\Gamma$ is endowed with the fine bornology. Look at the case $n=1$, the higher ones following in the same way. By the associativity (which, as Meyer points out, is not trivially true) and commutativity of the completed bornological tensor product, we have} 
(\CINF(V)\alcross\C\Gamma)^{\protens 2} & =(\CINF(V)\otimes\C\Gamma)\protens(\CINF(V)\otimes\C\Gamma)\\
& =(\CINF(V)\protens\CINF(V))\otimes\C\Gamma\otimes\C\Gamma\\
& =\CINF(V\times V)\otimes(\C\Gamma)^{\otimes 2}.
\end{align*}
Now we see that all changes that have to be performed on the complexes, are that all tensor products between factors $\CCINF(V)$ are projective tensor products. All other tensor products remain algebraic.\\
This applies to the cyclic modules $L^n(A,\Gamma,z)$, and the decomposition into homogeneous and inhomogeneous part remains valid in a completely analogous way. So we also obtain a complex that computes group hyperhomology. In this context, we also dispose of the Eilenberg-Zilber theorem, because in the proof all tensor products can be replaced by bornological tensor products. The existence of the $SBI$-sequence remains valid, see \cite{Meyerthesis}*{p. 103}. So all we have to check is that the shuffle maps $sh$ and $sh'$, which we have mentioned earlier and which are defined precisely in \cite{Loday}*{4.2.8 and 4.3.2} are bounded. But this is clear from the definitions of these maps, which are finite linear combinations of permutations.\\
So we obtain an isomorphism $\HC_*(\Gamma\allcross \CCINF(V))\cong\mathbb{H}_*(\Gamma; \mathcal{B}(\CCINF(V)))$ where $\mathcal{B}(\CCINF(V)))$ now is the complex used in the literature for the cyclic homology of Fr\'echet algebras. So the first part of the claim is achieved by the fact that the latter complex is $\Gamma$-equivariantly quasi-isomorphic to the complex $\Omega^\bullet V$ of $\mathcal{C}^\infty$-de Rham-forms on $V$.\\
In the part concerning cyclic cohomology, we have to take \emph{bounded} duals everywhere. One checks therefore the following steps. Since we want to avoid to define bornological cyclic modules we work directly with the bornological complexes and their bounded duals.\\
Let us go through each of the important steps. The one involving
$\mathcal{B}(\Psi)=\mathcal{B}(\widetilde{\Psi}/\Gamma)$ needs no
checking since we are dealing with an isomorphism of complexes, both having $(\C\Gamma)^{\times (n+1)}\otimes A^{\protens(n+1)}$ as entries. In regard of bounded functionals the choice of tensor product is irrelevant. The same remark applies to the step $\mathcal{B}(\widetilde{\Psi}/\Gamma)=\mathcal{B}(\widetilde{\Psi})/\Gamma$. On the other hand, the step which identifies functionals on the quotient $\mathcal{B}(\widetilde{\Psi})/\Gamma$ with $\Gamma$-invariant functionals on $\mathcal{B}(\widetilde{\Psi})$ makes use of the trivial fact that composing a bounded homomorphism $A\otimes_{\C\Gamma}\C\to\C$ on the right with the (bounded) quotient map $A=A\otimes_\C\C\to A\otimes_{\C\Gamma}\C$ gives a bounded map, and of the easily checked fact that an invariant bounded functional $A\to\C$ factors via the two bounded maps $A\to A\otimes_{\C\Gamma}\C\to\C$ in our case. This follows from the fact that the quotient map $A\to A\otimes_{\C\Gamma}\C$ has a bounded section here.\\
A bounded homomorphism of bornological vector spaces or complexes has a transpose between the bounded duals, and if the former is a quasi-isomorphism, then so is the latter. So we end up with the continuous dual $\Omega_\bullet V$ of the de-Rham complex.
\end{proof}
It is more convenient to deal with the hyper(co)homology with coefficients in the de-Rham complex or its dual than with the coefficients $\Omega^\ev V\leftrightarrows \Omega^\odd V$. So we observe that the latter complex splits canonically into an infinite direct sum of \emph{cochain} complexes, namely $\dots\to\Omega^\ev \to\Omega^\odd\to\Omega^\ev\to\dots\cong\bigoplus_j (\Omega^j V)[2j]$. We will use this splitting in the sequel.
In the following corollary, we shall mean by $\Homol$ cohomology with arbitrary supports, as opposed to the one with compact supports. This remark only applies when $\BG$ is non-compact.
\begin{cor}[Connes's twisted cohomology, \cite{Co2}, III.2.$\delta$.]
There is a natural isomorphism in both $\Gamma$ and $V$ between the bounded periodic cyclic cohomology of the crossed product $\CINF(V)\alcross\Gamma$ and the cohomology of the homotopy quotient $\EG\times_\Gamma V$, twisted by the vector bundle $\tau$ on the homotopy quotient that is induced from the tangent bundle on $V$:
$$
\HP^n\bigl(\CINF(V)\alcross\Gamma\bigr)\cong\bigoplus_j\Homol_\tau^{n+2j}(\EG\times_\Gamma V).
$$
\end{cor}
Note that the important feature of the right hand side is that it is the cohomology of a topological space, but possesses the correct functoriality - contravariant in $\Gamma$, but covariant in $V$. This is because there is Poincar\'e-Lefschetz duality $\Homol^n_{TV}(V)=\Homol^n(DV, SV)\cong\Homol_{\dim V-n}(V)$ ($TV,DV,SV$ being, of course the tangent, disc and sphere bundle respectively of $V$) that does \emph{not} require $V$ to be oriented (more precisely, it is duality for the orientable  manifold with boundary $DV$, which is always orientable because the total space of the tangent bundle possesses an almost complex structure). This duality is $\Gamma$-equivariant and gives therefore rise to the covariant functoriality in $V$.
\begin{proof}[Proof (of the corollary)]
By \cite{Co2}, III.2.$\delta$, Proposition $11$, the twisted cohomology is computed by the bicomplex $C^\bullet$ whose $C^{pq}$ is the space of totally antisymmetric equivariant maps $\Gamma^{p+1}\to\Omega_{-q}V$, in other words by
$$
C^\bullet=\Hom(\Lambda_\bullet \C\Gamma,\Omega_{-\bullet}V)
$$
with boundary induced by
$$
d(\gamma_0,\dots,\gamma_{n+1})=\sum_{j=0}^{n+1}(-1)^j \gamma_0\wedge\dots\wedge\overset{\vee}{\gamma_j}\wedge\dots\wedge\gamma_{n+1}
$$
and by the de-Rham differential. One checks that $\Lambda_\bullet \C\Gamma$ is a free $\C\Gamma$-resolution of the constant module $\C$, because this is exactly the Koszul complex associated to the vector space $\C\Gamma$ together with the augmentation $\epsilon\colon \C\Gamma\to\C$. So we get back the hypercohomology of $\Gamma$ with coefficients in $\Omega_{\bullet}V$.
\end{proof}
\subsection{Equivariant homology, shriek maps and the Leray-Serre spectral sequence}
\subsubsection{Equivariant homology}
We are now in a position to apply Poincar\'e duality, thereby establishing an isomorphism with \emph{equivariant} de-Rham cohomology. We therefore make the following
\begin{defn}
Let $V$ be a smooth manifold, not necessarily compact or orientable, acted on smoothly by the group $\Gamma$.\\
The \emph{equivariant de-Rham homology} of $V$ is defined as $$\Homol^\Gamma_*(V)\defeq\mathbb{H}_*(\Gamma;\ \Omega_\ev V\lra\Omega_\odd V),$$ the group hyperhomology of $\Gamma$ with coefficients in the complex of complex-valued de-Rham currents of $V$.\\
Similarly, we define \emph{equivariant de-Rham cohomology} as
$$\Homol_\Gamma^*(V)\defeq\mathbb{H}^*(\Gamma; \ \Omega^\ev V\lra\Omega^\odd V),$$
the group hypercohomology of $\Gamma$ with coefficients in the complexified de-Rham complex of $V$.
\end{defn}
Examples of these groups will be given by identifying them with the ordinary homology resp. cohomology groups of the homotopy quotient of $V$ by $\Gamma$ in lemma \ref{examplesforequivarianthomology}(1).
\subsubsection{Poincar\'e duality}
Next we want to observe that these groups are in fact isomorphic to the hyperhomology groups we have met by calculating the cyclic homology. The interesting feature is only the functoriality.\\
\begin{lemmaanddef}\label{shriek}
Suppose $V$ is compact, and $\Gamma$ acts by orientation-preserving diffeomorphisms on $V$. Then there are isomorphisms which are natural in $\Gamma$
\begin{align*}
\mathbb{H}_*(\Gamma; \Omega^\ev V\lra\Omega^\odd V)&\ \cong\ \Homol_{2\Z+*+\dim V}^\Gamma(V)\\
\mathbb{H}^*(\Gamma; \Omega_\ev V\lra\Omega_\odd V)&\ \cong\ \Homol^{2\Z+*+\dim V}_\Gamma(V).
\end{align*}
These groups are indexed by a $\Z/2\Z$-grading.\\
By definition, the functoriality in $V$ is by fiberwise shriek maps, i. e. if $f\colon V\to W$ is a smooth $\Gamma$-equivariant orientation-preserving map of smooth compact orientable $\Gamma$-manifolds, denote the resulting map $\Homol^*_\Gamma(V)\to\Homol^*_\Gamma(W)$ by $f_!$. So $f_!$ is, by definition, the map which makes the following diagram commutative:
$$
\xymatrix{
  {\mathbb{H}^*(\Gamma;\ \Omega_\ev\lra\Omega_\odd V)} \ar[d]^{f_*}  & {\Homol_\Gamma^{2\Z+*+\dim V}(V)} \ar[l]_-{\cong}\ar[d]^{f_!} \\
  {\mathbb{H}^*(\Gamma;\ \Omega_\ev\lra\Omega_\odd W)}               & {\Homol_\Gamma^{2\Z+*+\dim W}(W)}\ar[l]_-{\cong},           \\ 
}
$$
The map $f_!$ is called the shriek map.\\
Similar statement and definition can be made for homology.
\end{lemmaanddef}
\begin{proof}
One easily checks that the map of complexes $\Omega^\bullet V\to\Omega_\bullet V$, given by
$$
\omega\mapsto\int_V-\wedge\omega
$$
is a $\Gamma$-equivariant morphism of complexes, because $\Gamma$ preserves the orientation. By Poincar\'e duality it induces an isomorphism on homology, so it induces an isomorphism on group hyperhomology and cohomology.
\end{proof}
\subsubsection{The Leray-Serre spectral sequence}
We now want to exhibit  a topological space, namely the homotopy quotient (also known as the Borel construction) whose cohomology are the equivariant groups above. In view of computability, this is the most useful formulation because it allows topology to enter the picture. The  spectral sequence we found then becomes the well-known Leray-Serre spectral sequence computing the cohomology of the total space of a fibration from that of base and fiber. The interesting feature about all that is only the shriek functoriality under $\Gamma$-equivariant diffeomorphisms. We want to avoid Poincar\'e duality in the base $\BG$, in other words we need fiberwise Poincar\'e duality for the fibration $\pi\colon \EG\times_\Gamma V\to\BG$. We therefore make the following digression on the Leray-Serre spectral sequence before stating the fact on the homotopy quotient.\\
Let $\pi\colon E\to B$ be a fibration over $B$. So there is a commutative diagram
$$
\xymatrix{
{E}\ar[rr]^{\pi}\ar[dr]&&{B}\ar[dl]^{p_*}\\
&{\pt}&
}
$$
which induces a similar diagram on the categories of sheaves with values in vector spaces of the respective spaces. Recall that the category of sheaves over a given space with values in an abelian category is again an abelian category.\\
Denote by $\underline{A}^\bullet_E$ the complex of sheaves of differential forms on $E$. The general theory of the Grothendieck spectral sequence associated to a composition of functors implies (\cite{Weibel}, section 5.7) that the cohomology of $E$ is given by the expression
$$\Homol^*(E;\C)=(\mathbb{R}p_*)(\pi_*\underline{A}^\bullet_E),$$
where $\mathbb{R}$ is the total derived functor. This is because there is an isomorphism
$$(\mathbb{R}p_*)(\pi_*\underline{A}^\bullet_E)=R(p_*\pi_*)(\C_E).$$
The last space is the cohomology of $E$ because $R(p_*\pi_*)(\C_E)=\Homol^*(S(\underline{A}^\bullet_E))=\Homol^*(\Omega^\bullet E)=\Homol^*(E)$, where $S$ is the global sections functor.. The Leray-Serre spectral sequence then is the spectral sequence computing this hyperderived functor.\\
End of digression.
\begin{lemma}\label{examplesforequivarianthomology}\mbox{ }\\
\vspace{-1em}
\begin{enumerate}
\item
There are natural (in both $\Gamma$ and $V$) isomorphisms
\begin{align*}
\Homol_\Gamma^*(V)&\cong\Homol^*(\EG\times_\Gamma V)\\
\Homol^\Gamma_*(V)&\cong\Homol_*(\EG\times_\Gamma V).
\end{align*}
\item
The Leray-Serre spectral sequence
$$
E_2^{pq}=\Homol^p\bigl(\BG,\ \underline{\Homol}^q(V, \C))\implies \Gr^p\ \Homol^{p+q}(\EG\times_\Gamma V)
$$
coincides with the spectral sequence computing group hypercohomology
$$
E_2^{pq}=\Homol^p\bigl(\Gamma;\ \Homol^q(V))=\Homol^p\bigl(\BG;\ \underline{\Homol}^q(V)\implies \Gr^p\ \mathbb{H}^{p+q}(\Gamma;\ \Omega^\bullet V)
$$
in such a way that the representation of $\Gamma$ on $\Homol^\bullet V$ as the monodromy representation $\Gamma=\pi_1(\BG)\to \End(\Homol^\bullet(V))$, thereby establishing the coefficient system $\underline{\Homol}^\bullet(V)$.
\end{enumerate}
\end{lemma}
\begin{proof}
\begin{enumerate}
\item
In fact, by \cite{Weibel}, theorem 6.10.10, the group $\Homol^{p+q}((\EG\times V)\big/\Gamma;\ \C)$ is equal to the group hypercohomology $\mathbb{H}^{p+q}(\Gamma;\ \Hom_{\mathbf{Ab}}(S_*(\EG\times V),\C))$, where $S_*$ is the singular chain complex functor. But $\Hom_{\mathbf{Ab}}(S_*(\EG\times V),\C)$ is by the de-Rham theorem $\Gamma$-equivariantly chain homotopy equivalent  to $\Omega^*(\EG\times V)$, which, in turn, is chain homotopy equivalent to $\Omega^* V$ via the map $\Omega^*V\to\Omega^*\EG\times V$ by the Poincar\'e lemma. This map is $\Gamma$-equivariant. This proves the statement.\\
\item
The hypercohomology spectral sequence is the one associated to the double complex $(\Omega^\bullet\EG\protens\Omega^\bullet V)^\Gamma$ (more precisely, the second one of the two canonical ones), because $\Omega^\bullet\EG$ is $\Gamma$-equivariantly chain homotopy equivalent to an injective resolution of $\C$ in the category of $\C\Gamma$-modules (\cite{Puschnigg}). So it is the spectral sequence associated to the filtered complex $\Tot(\Omega^\bullet\EG\otimes\Omega^\bullet V)^\Gamma$, faltered by columns. On the other hand, the inclusion
$$
\Tot(\Omega^\bullet\EG\otimes\Omega^\bullet V)^\Gamma\hookrightarrow\Tot(\Omega^\bullet\EG\protens\Omega^\bullet V)^\Gamma=\Omega^\bullet (\EG\times_\Gamma V)
$$
into the completed projective tensor product with respect to a natural Frechet topology on the differential forms is a quasi-isomorphism, where the right hand side takes its filtration from the left hand. The inclusion is a quasi-isomorphism because it induces an isomorphism on the $E_2$ term.\\
After all, we have defined a filtration on $\Omega^\bullet (\EG\times_\Gamma V)$. This filtration is exactly the Cartan filtration of differential forms, defined in \cite{Bryl}, definition 1.6.8. Consequently (see \cite{Bryl}, theorem 1.6.10) the associated Cartan spectral sequence identifies with the Leray spectral sequence. In our framework, his argument is that $\underline{A}_{\BG}^\bullet\otimes\pi_*\underline{A}_{\EG\times_\Gamma V}^\bullet$ is a soft resolution of the sheaf complex $\pi_*\underline{A}_{\EG\times_\Gamma V}^\bullet$, so
\begin{multline*}
(\R p_*)(\pi_*\underline{A}_{\EG\times_\Gamma V}^\bullet)=\Homol^*(p_*(\pi_*\underline{A}_{\EG\times_\Gamma V}^\bullet\otimes\pi_*\underline{A}_{\BG}^\bullet))\\
=\Homol^*(S(\pi_*\underline{A}_{\EG\times_\Gamma V}^\bullet\otimes\underline{A}_{\BG}^\bullet)).
\end{multline*}
Now $S(\pi_*\underline{A}_{\EG\times_\Gamma V}^\bullet\otimes\underline{A}_{\BG}^\bullet)$ is a double complex which gives the Cartan filtration of differential forms. The associated spectral sequence has $E_2$ term
$$
E_2^{pq}=\Homol^p(\BG;\ \underline{\Homol}^q(V;\ \C))\implies \Gr^p\ \Homol^{p+q}(\EG\times_\Gamma V).
$$
\end{enumerate}
\end{proof}
\begin{rmk}\label{caution}
The action of $\Gamma$ on the local coefficient system $\underline{\Homol}^\bullet(V)$ does \emph{not} specify the spectral sequence. For example, let $\Gamma$ be the fundamental group of a Riemann surface of finite genus, and let $V$ be the boundary of the Poincar\'e plane. The group $\Gamma$ acts in an orientation-preserving way on $V$, so the monodromy representation of $\Gamma$ on $\Homol^*(V)$ is trivial, but $\EG\times_\Gamma V$ is the unit sphere bundle of the tangent bundle of the Riemann surface, and this is not a homologically trivial bundle, as it has non-zero Euler class. So $\Homol^*(\EG\times_\Gamma V)$ is not product cohomology.
\end{rmk}
We now want to make explicit what the aforementioned shriek maps become under these isomorphisms.

\begin{rmk}
This proof shows how we can integrate the shriek map into the abstract Leray-Serre picture. Namely, one can show that there is defined a morphism 
$$
\pi_*\underline{A}^\bullet_{\EG\times_\Gamma V}\to\pi_*\underline{A}_{\EG\times\Gamma W}^{\bullet-k+l}
$$
of complexes of sheaves on $\BG$ by the map
\begin{multline*}
S(\pi_*\underline{A}_{\EG\times_\Gamma V}^\bullet\otimes\underline{A}_{\BG})=(\Omega^\bullet\EG\protens\Omega^\bullet V)^\Gamma\\
\to \Omega^\bullet\EG\protens\Omega^{\bullet-k+l}W)^\Gamma=S(\pi_*\underline{A}_{\EG\times_\Gamma W}^{\bullet-k+l}\otimes\underline{A}_{\BG}^\bullet)
\end{multline*}
not only on the level of the global sections but on the level of the sheaves themselves.
\end{rmk}
We now discuss the functoriality of the cohomology of the homotopy quotients induced by the the map $V\to \pt$ of the coefficient manifold in more down-to-earth-language.
\begin{defn}
The map on differential forms
$$
\Omega^\bullet(\EG\times_\Gamma V)=\Tot (\Omega^\bullet\EG\protens\Omega^\bullet V)^\Gamma\to\Tot\ (\Omega^{\bullet-k}\EG)^\Gamma=\Omega^{\bullet-k} \BG
$$
is  called integration along the fiber and denoted $\pi_!$.
\end{defn}
\begin{lemma}\label{fiberintegration}
\mbox{ }
\begin{enumerate}

\item
The map $\pi_!$ commutes with $d$, so $\pi_!$ induces a map
$$
\Homol^{p+k}(\EG\times_\Gamma V)\to\Homol^p(B).
$$
on cohomology.
\item
This map satisfies
$$
\int_B\tau\wedge\pi_!\omega=\int_E\pi^*\tau\wedge\omega
$$
for any differential forms $\omega\in\Omega^{|\omega|}\EG\times_\Gamma V,\ \tau\in\Omega^{|\omega|+\dim V-\dim\BG}\BG$.
\item
The induced map on cohomology comes from the map
$$
E_2^{pk}=\Homol^p(\BG,\ \underline{\Homol}^k(V))\to\Homol^p(\BG,\ \underline{\C}),
$$
induced by the (equivariant) integration $\Homol^p(V)\to\C$ on the $E_2$-term of the Leray-Serre spectral sequence.
\end{enumerate}
\end{lemma}
The map $\Homol^{p+k}(\EG\times_\Gamma V)\to\Homol^p(B)$ is also called integration along the fiber. The equality in the second statement is known from the integration along the fiber-map in the theory of smooth vector bundles (\cite{BottTu}).
\begin{proof}
\begin{enumerate}
\item
The commutativity of the diagram
$$
\xymatrix{
{(\Omega^p\EG\protens\Omega^k V)^\Gamma}\ar[d]^{1\otimes\int}\ar[r]^-{d} & {(\Omega^{p+1}\EG\protens\Omega^kV)^\Gamma}\ar[d]^{1\otimes\int}\\
{(\Omega^p\EG)^\Gamma}\ar[r]^-{d}&{(\Omega^{p+1}\EG)^\Gamma}
}
$$
is immediate (and is, in our context nothing but a special case of the fact we already used that morphisms of coefficient complexes induce morphisms on hypercohomology complexes).
\item
Let $(\chi_i)$ be a partition of unity on $\BG$ subordinate to an open cover that trivializes the bundle $\EG\times_\Gamma V\to\BG$. It suffices to consider only finitely many functions because the support of $\omega$ is compact, i. e. it suffices that $\pi^{-1}$ of the support of the $(\chi_i)$ comprises the support of $\omega$. Then by the linearity of the equation in question it suffices to show it for $\chi_i\tau$, i. e. in the case of a trivial fibration $X\times V\to X$, where $X$ is smooth and compact. Let $\pi_2$ denote the projection $X\times V\to V$. Approximate $\omega\in\Omega^\bullet(X\times V)=\Omega^\bullet X\protens\Omega^\bullet V$ by a series $\sum\pi^*\omega_v^j\wedge\pi_2^*\omega_h^j, $ where $\omega_h^i\in\Omega^\bullet X,\omega_v^i\in\Omega^\bullet V$.\\
Then 
\begin{align*}
\int_X\tau\wedge\pi_!\omega&=\sum_j\int_X\tau\wedge\pi_!(\pi^*\omega_v^j\wedge\pi_2^*\omega_h^j)\\
&=\sum_j\int_X\tau\wedge\omega_v^j\int_V\omega_h^j\\
&=\sum_j\int_{X\times V}\pi^*(\tau\wedge\omega_v^j)\wedge\pi_2^*\omega_h^j\\
&=\sum_j\int_{X\times V}\pi^*\tau\wedge\pi^*\omega_v^j\wedge\pi_2^*\omega_h^j=\int_{X\times V}\pi^*\tau\wedge\omega,
\end{align*}
where the first equality holds by a fundamental property of the Lebesgue integral, the second one by the definition of fiber integration, and the third one by Fubini's theorem.

\item
Immediate from the observation we made that the Leray-Serre spectral sequence is the one associated to the bicomplex $(\Omega^\bullet\EG\otimes\Omega^\bullet V)^\Gamma$.
\end{enumerate}
\end{proof}
\begin{rmk}
In the lines of the preceding remark, the integration along the fiber comes from a morphism of complexes of sheaves on $\BG$
$$
\pi_*\underline{A}^\bullet_{\EG\times_\Gamma V}\to(\id_*)\underline{A}^\bullet_{\BG}[-n]=\underline{A}^\bullet_{\BG}[-n].
$$
where $[-n]$ denotes shifting a complex by $-n$ degrees (and in general one can show that $f_!$ induces a morphism $\underline{A}^\bullet_{\EG\times_\Gamma V}\to\underline{A}^{\bullet-k+l}_{\EG\times_\Gamma W}$.)\\
Applying the functor $\R p_*$ to the latter sequence, we obtain the integration over the fiber map
\begin{multline*}
\pi_!\colon \Homol^*(E;\ \C)=(\R p_*)(\pi_*\underline{A}^\bullet_{\EG \times_\Gamma V})\\
\to (\R p_*)(\underline{A}^\bullet_{\BG})[-k]=(Rp_*)(\underline{\C}_{\BG})[-k]=\Homol^{*-k}(\BG).
\end{multline*}
Let us look at the induced spectral sequence map on the $E_2$ level. As we have already noted, we get
$$
E_2^{pk}=(R^p)(\Homol^k(\pi_*(\underline{A}^\bullet_{EG\times_\Gamma V})))=\Homol^p(\BG,\ \underline{\Homol}^k(V))\to\Homol^p(\BG,\ \C),
$$
So passing to the abutment, this spectral sequence map induces
\begin{gather*}
\Homol^{p+k}(\EG\times_\Gamma V)=\bigoplus_{x+y=p+k}E_\infty^{xy}\to E_\infty^{pk}\\
\hookrightarrow E_2^{pk}=\Homol^p(\BG;\ \underline{\Homol}^k(V))\to\Homol^p(\BG;\ \C),
\end{gather*}
where the second map is the projection in the direct sum onto the ``upper edge'' terms, and where the third map is injective because there are only zero differentials arriving at these terms. So this is the way we can read off the integration along the fiber directly from the Leray-Serre spectral sequence in the form it occurs.
\end{rmk}
\subsubsection{The spectral sequence for cyclic cohomology}
We now reach our first goal, namely the identification of Nistor`s spectral sequence.
\begin{prop}
The cohomological spectral sequence of corollary \ref{spectralsequence} coincides with an infinite direct sum of Leray-Serre spectral sequences.
\end{prop}
\begin{proof}
Putting together what we have said about both spectral sequences, we see that the Leray-Serre sequence is the hypercohomology spectral sequence computing $\mathbb{H}^*(\Gamma; \Omega^*(V))$. In turn, the spectral sequence of \ref{spectralsequence}, which coincides with Nistor's, is the one computing $\mathbb{H}^*(\Gamma; \Omega_\ev V\leftrightarrows\Omega_\odd V)\cong\mathbb{H}^*(\Gamma; \Omega^\ev V\leftrightarrows\Omega^\odd V)$. The complex $\Omega^\ev V\leftrightarrows \Omega^\odd V$ is a direct sum of complexes $\Omega^* V$, and it is this in the category of complexes of $\C\Gamma$-modules as well.
\end{proof}
\section[Lattices, the Furstenberg boundary and the splitting principle]{Lattices in semisimple Lie groups, the Furstenberg boundary and the splitting principle}We shall apply the previous computation of cyclic cohomology to $\K$-theory and index theory, thereby approaching the idempotent theorem for the group ring of higher rank lattices. Let us begin with lattices in $\SL(q,\C)$.
\subsection{The case $G=\SL(n,\C)$}
\begin{prop}\label{splprplelementar}
Let $G$ be the simple Lie group $\SL(n,\C)$ for $n\ge 2$, $\Gamma$ a torsion-free (not necessarily cocompact, nor even of finite covolume) discrete subgroup of $G$, $K=\SU(n)$ the preferred maximal compact subgroup of $G$, $\EG=G/K$ the symmetric space with its left $\Gamma$-invariant Riemannian structure of non-positive curvature, (as explained in \ref{symmetric space}) and $B$ the minimal parabolic and maximal solvable subgroup of upper triangular matrices of $G$, $B\backslash G$ the corresponding flag manifold, seen as a right $\Gamma$-space. Consider the fibration
$$
\xymatrix{
{B\backslash G\times_\Gamma G/K}\ar[d]^{p}\ar@{=}[rr] && {(B\backslash G)\times_\Gamma\EG}\ar[d]\\
{\Gamma\backslash G/K}\ar@{=}[r] &\BG &{\pt\times_\Gamma\EG}\ar@{=}[l]
}
$$
induced by collapsing $B\backslash G$ to a point. It has the property that the map on cohomology
$$
p^*\colon \Homol^*(\BG)\to\Homol^*((B\backslash G)\times_\Gamma\EG)
$$
(with rational coefficients, without support condition) is injective. Consequently,
$$
p_*\colon \Homol_*((B\backslash G)\times_\Gamma\EG)\to\Homol_*(\BG)
$$
is surjective.
\end{prop}
Recall that there are torsion-free subgroups of finite index in $\SL(n,\Z)$ (see \cite{Brown}). They provide therefore canonical examples for $\Gamma$.
\begin{proof}
The ordinary splitting principle for complex vector bundles can most easily be formulated as follows (see for example \cite{BottTu}). There is a functor $\Fl$ from vector spaces to spaces, associating to a $n$-dimensional vector space $U$ its flag manifold $\Fl(U)$, i. e. the set of maximal flags
$$
0=U_0\subset U_1\subset \dots \subset U_{n-1} \subset U_n =U,
$$
$U_i$ being an $i$-dimensional subspace of $U$. This set is a $\SL(n,\C)$-homogeneous space and in particular a compact smooth manifold because $\SL(n,\C)$ acts on it transitively with isotropy subgroup the closed subgroup $B$. So $\Fl(U)\cong B\backslash G$. The functor $\Fl$ extends to a functor from the category of complex vector bundles over a given space to that of fiber bundles over it. The total space of the flag bundle $\Fl(V)$ over $X$ associated to a fixed real vector $V\to X$ is a splitting manifold for the vector bundle, i. e. the bundle $V$ pulled back to the total space of $\Fl(V)$ along the projection $\Fl(V)\to X$ splits into a direct sum of line bundles, and, important to us, the cohomology $\Homol^*(X)$ of the base injects by pulling back along the projection into the cohomology $\Homol^*(\Fl(V))$ of $\Fl(V)$ (see \cite{BottTu} for the splitting principle). Trivial vector bundles with, say, fiber $U$ are, of course, mapped to trivial fibrations with fiber $\Fl(U)$. In particular, for the trivial vector bundle of rank $n$ over $G/K$ we have
$$
\Fl(U\times G/K)=\Fl(U)\times G/K =B\backslash G\times G/K.
$$
Now the bundle $U\times G/K$ comes with a canonical $\Gamma$ action, namely the diagonal action induced from the action of $\Gamma$ as a subgroup of $\SL(n,\C)$ on $\C^n=U$ and the natural action on $\EG$. So there is a $\Gamma$-action on $\Fl(U\times G/K)$, because $\Fl$ is a functor. Let us now show that the following diagram commutes.
$$
\xymatrix@C=4.5em{
\Fl(U\times G/K)\ar[d]^{\Fl(\gamma\inv\times\gamma)}\ar@{=}[r]&\Fl(U)\times G/K\ar[d]^{\gamma\inv\times\gamma}\\
\Fl(U\times G/K)\ar@{=}[r]&\Fl(U)\times G/K.
}
$$
For the factor $G/K$ this is exactly the statement that $\Fl$ is a functor from the category of vector bundles to that of fiber bundles. For the factor $\Fl(U)=\Fl(\C^n)$ this follows from the fact that the homeomorphism $G/B\to \Fl(U), xB\mapsto xB \zeta$ where $\zeta$ is the standard flag of $\C^n$ is equivariant. So we may divide on both sides of the equation $\Fl(U\times G/K)=\Fl(U)\times G/K$ by $\Gamma$, because $\Gamma$ acts freely and properly o
n $\EG$, and therefore so does it on vector or fiber bundles over $\EG$.\\
Second one checks that the natural map from $\Fl(V)$ to $\Fl(V/\Gamma)$ identifies $\Fl(V)/\Gamma$ with $\Fl(V/\Gamma)$ when the action on the base is free and proper. So altogether we have shown that $\Fl(U\times_\Gamma G/K)=\Fl(U\times G/K)/\Gamma=B\backslash G\times_\Gamma G/K$. Therefore the space under consideration is a splitting manifold, and the projection maps coincide by definition, i. e. the map $p$ is in fact the map that provides the embedding of the cohomology of the base $\BG$ into that of the splitting manifold. This proves the first half of the theorem.\\
The second follows from the universal coefficient theorem. 
\end{proof}
\begin{rmk}
This proof shows that we can even tell the cohomology ring of $\EG\times_\Gamma G/B$ as a ring. Since the structure group $\Gamma$ of the complex vector bundle $\EG\times_\Gamma U$ is discrete, it is flat, hence it has no non-zero rational Chern classes. So the cohomology of the associated flag bundle is the product
$$
\Homol^*(\EG\times_\Gamma G/B)=\Homol^*(\Gamma)\otimes \Homol^*(G/B)=\Homol^*(\Gamma)[x_1,\dots,x_n]\bigg/\biggl(\thickspace\prod_{i=1}^n (1+x_i)=1\biggr).
$$
\end{rmk}
\subsection{Connected semisimple Lie groups}
Now we want to extend the proposition to the case of more general (complex connected) semisimple Lie groups. Let us first gather some basic facts on them. The following facts can be found in any standard reference on Lie groups, for instance in \cite{Serre}.
\begin{enumerate}
\item
Any complex semisimple Lie algebra $\mathfrak{g}$ has a unique (up to isomorphism) compact real form, i. e. a real subalgebra $\mathfrak{g}_0$ such that $\mathfrak{g}\cong\mathfrak{g}_0\otimes\C$ and $\mathfrak{g}_0$ is compact (it is the Lie algebra of a compact group).
\item
Any complex semisimple Lie group (algebra) has a subgroup (subalgebra) which is at the same time maximal solvable and cocompact (such subgroups are called Borel or minimal parabolic subgroups) .
\item
A maximal abelian real subalgebra of $\mathfrak{g}$ is a Cartan subalgebra (i. e. nilpotent and its own normalizer) iff any of its elements is semisimple (diagonalizable in the adjoint representation). $\dim_\R\mathfrak{g}$ is then called the rank $\rk \mathfrak{g}$ of $\mathfrak{g}$.\\
If $\mathfrak{g}$ is compact, then a subalgebra is Cartan iff it is maximal abelian.
\item
If $\mathfrak{g}_0$ is a real form of $\mathfrak{g}$, then $\mathfrak{g}_0$ is abelian (resp. semisimple) in the real sense iff $\mathfrak{g}$ is in the complex sense. Moreover, $\rk \mathfrak{g}=2\rk\mathfrak{g}_0$ because $\mathfrak{h}_0\subset\mathfrak{g}_0$ is Cartan iff $\mathfrak{h}_0\otimes\C\subset\mathfrak{g}$ is.
\item
Any compact Lie algebra is reductive (a direct sum of an abelian and a semisimple algebra).
\end{enumerate}
\begin{rmk}\label{condition}
 One can prove the following statement. Let $G$ be a complex semisimple Lie group. Let $K$ be the maximal compact subgroup corresponding to the compact real form of the Lie algebra. Let $T$ be the maximal torus of $K$. 
Then $T$ is also a maximal torus for a Borel subgroup $B$.
\end{rmk}
The conclusion obviously holds in the case of $G=\SL(n,\C)$ for all $n$.\\
Since the statement of the conclusion is the only condition on the Lie group we make use of in the entire chapter, the remark allows to replace $\SL(n,\C)$ everywhere in the chapter, in particular in theorem \ref{essence}, by a complex semisimple Lie group.
\begin{prop}
In the situation of the preceding proposition, let $\Gamma$ be a torsion-free lattice of $G$.
\begin{enumerate}
\item 
There is a homotopy commutative diagram as follows.
$$
\xymatrix{
{G\times_B\classb B}\ar[d]^{(8)}\ar@{=}[r]&{G\times_B\classb B}\ar[d]^{(5)}&{K\times_T\classb T}\ar[d]^{(7)}\ar[l]^{(6)}\\
{\EG\times_\Gamma G\times_B\classb B}\ar[r]\ar[d]&{\classs B}\ar[d]^{(4)}   &{\classs T}\ar[d]^{(3)}\ar[l]^{(1)}\\
\BG\ar[r]&{\classs G}           &{\classs K}\ar[l]^{(2)}
}
$$
where the maps $(1)$, $(2)$, $(3)$, $(4)$, $(6)$ and $(8)$ are induced by inclusions, all others by obvious projections.
\item
The columns are homotopy fibrations.
\item
The arrows $(1)$, $(2)$ and $(6)$ are homotopy equivalences.
\item
The arrows $(3)$ and $(4)$ are injective on rational cohomology.
\item
The arrow $(7)$ and hence the maps $(5)$ and $(8)$ are surjective on rational cohomology.\item
The left hand fibration $\EG\times_\Gamma G\times_B\classb B\to \BG$ is cohomologically trivial as well. In particular $\pi^*$ is injective.
\end{enumerate}
\end{prop}
\begin{proof}\mbox{ }\\
\vspace{-1.5em}
\begin{enumerate}
\item
The commutativity up to homotopy follows immediately from the definitions of the maps and the fact that $\classs$ is a functor, except for the lower left square. The latter is a special case of the following statement.\\
For any connected Lie group $G$ with the two closed subgroups $H$ and $H'$ the following diagram is a pull-back diagram.
$$
\xymatrix{
{\classb H\times_H G\times_{H'} \classb H'}\ar[r]\ar[d] & \classs H'\ar[d]\\
\classs H\ar[r]&\classs G,
}
$$
where the upper and left arrows are given by projecting onto the right resp. left factor, and the other two arrows are given by maps $H\to G$ and $H'\to G$. This is checked by verifying that both compositions classify the same $G$-principal bundle, namely $\classb H\times_H G\times G\times_{H'} \classb H'\to \classb H\times_H G\times_G G\times_{H'} \classb H'=\classb H\times_H G\times_{H'} \classb H'$.\\
This proves the commutativity of the lower left square since for $H=\Gamma$ and $H'=B$ this pullback $\EG\times_\Gamma G\times_B \classb B$ is homotopic to $\EG\times_\Gamma G/B$ (because it projects onto the first two factors with contractible fiber $\classb B$).\\
(Note however that the map $\BG\to \classs G$ is in cohomology the projection onto the zeroth component degree by the Chern-Weil theory of characteristic classes.)\\
(Second note: A principal bundle is by definition a space equipped with a \emph{right} group action satisfying some properties. If the formula $b.(x,y)=(xb\inv,by)$ is used to define the homotopy quotient, then this is a left action. However, changing to $(xb,b\inv y)$ does not change the quotient).
\item
For the first column, this is immediate. For the two others, this follows from the special case $H=\pt$ of the statement made on the preceding item. In fact, it shows that the compositions of the vertical maps in question factor via $\classs\pt=\pt$, so is homotopically trivial. Second, $\classs H'\to \classs G$ is surjective, since $\classb G$ is also a $\classb H'$, thus $\classs G=\classb G/G=\classb G/H'/G=\classb H'/H'/G=\classs H'/G$.
\item
Since $B/T$ and $G/K$ are contractible, the statement follows from the same reasoning as in the preceding item.
\item
The statement for the map (3) is part of A. Borel's homotopical formulation of the splitting principle, see \cite{Borel}, p. 67. In fact, he proves that the cohomology of $\classs K$ is isomorphic by this map to the invariants under the action of the Weyl group of $T$ in $K$. The statement for the map $(4)$ follows from that for $(3)$ given that $(1)$ and $(2)$ are homotopy equivalences.
\item
In terms of the Leray-Serre spectral sequence for a fibration $F\xrightarrow{i} E\to B$, the map $i^*$ is given by
$$
\Homol^n(E)\to E_\infty^{n0}\hookrightarrow E_2^{0n}=\Homol^n(F).
$$
The second map is injective because on the left hand edge only arrive zero differentials. If this map is also surjective, as is proved in the preceding items, all other differential vanish as well. So the spectral sequence degenerates at the $E_2$ term and gives product cohomology.
\end{enumerate}
\end{proof}
\noindent
So the fibration $G/B\to \EG\times_\Gamma\to\BG$ splits in rational cohomology.
\subsection{$\SL(n,\C)$ in $\K$-theory}
In a third presentation of the splitting principle, we now want to extend the results of theorem \ref{splprplelementar} from the case of de-Rham cohomology to $\K$-theory. This means that now we want to take torsion phenomena into account. So we set out to prove the \emph{split injectivity} of the map $\pi^*\colon \K^*(\BG)\to\K^*(\EG\times_\Gamma G/B)$ if $G=\SL(n,\C)$. This gives another proof of theorem \ref{splprplelementar} from a different (though related) and more general viewpoint because the Chern character is a rational isomorphism.\\
Consider the following sequence of ever larger subgroups of $G=\SL(n,\C)$.
$$
\begin{array}{ccccc}
G_1&=&B&=&\left\{\text{all matrices of the form}
{\tiny
\left(
\begin{array}{cccccc}
* & * & * &\ldots & *& *\\
0 & * & * &\ldots & *& *\\
0 & 0 & * &\ldots & *& *\\
\vdots & \vdots & \vdots & \ddots &&\vdots\\
0 & 0 &0  &   \dots&* &*\\
0 & 0 & 0 & \cdots & 0& *
\end{array}
\right)
}
\right\}\cap G
\\
\\
G_2&=&&&\left\{\text{all matrices of the form}
{\tiny
\left(
\begin{array}{cccccc}
* & * & * &\ldots & *& *\\
0 & * & * &\ldots & *& *\\
0 & 0 & * &\ldots & *& *\\
\vdots & \vdots & \vdots & \ddots &&\vdots\\
0 & 0 &0  &   \dots&* &*\\
0 & 0 & 0 & \cdots & *& *
\end{array}
\right)
}
\right\}
\cap G\\
\vdots&&&&\vdots\\
G_{n-1}&=&&&\left\{\text{all matrices of the form}
{\tiny
\left(
\begin{array}{cccccc}
* & * & * &\ldots & *& *\\
0 & * & * &\ldots & *& *\\
0 & * & * &\ldots & *& *\\
\vdots & \vdots & \vdots & \ddots &&\vdots\\
0 & * & *  &   \dots&* &*\\
0 & * & * & \cdots & *& *
\end{array}
\right)
}
\right\}
\cap G\\
G_n&=&G,&&
\end{array}
$$
i. e. at every step one more column is filled.\\
The fibration $\classs B\to \classs G$ thus factors as a sequence of fibrations
$$
\classs B\to\classs G_2\to\dots\to\classs G.
$$
We now claim that each single step is cohomology-injective.
\begin{lemma}
The fiber of the fibration $\classs G_{i-1}\to \classs G_i$ is (homotopic to) $\cpspace^{n-i}$.
\end{lemma}
\begin{proof}
We have to identify the homogeneous space $G_i/G_{i-1}$. $G_i$ acts on $\C^{i}$ by restricting to the right lower corner. The action passes to a transitive action on $\cpspace^{i-1}$. $G_{i-1}$ is the stabilizer of the base-point of $\cpspace^{i-1}$ corresponding to the first basis vector.\\
\end{proof}
\begin{lemma}
Each fibration step $\classs G_i\to\classs G_{i+1}$ can be identified with a universal fibration, i. e. with a fibration $\cpspace^{i-1}\hookrightarrow \proj(\vare_i)\to \classs \U(i)$ where $\proj(\vare_i)$ is the projectivization of the universal bundle $\vare_i\to\classs\U(i)$.
\end{lemma}
\begin{proof}
By looking at the foregoing proof we see that we may reduce the problem to the maximal compact subgroups (which may be done for purposes of homotopy). The $i$-step fibration $\classs G_i\to\classs G_{i+1}$ is homotopic to
$$
\xymatrix{
\U(i)/\U(1)\times\U(i-1)\ar[r] &\classs(\U(1)\times\U(i-1))\ar[d]\\
&\classs \U(i),
}
$$ which is exactly the universal one.
\end{proof}
\begin{lemma}
The universal fibration $\cpspace^{i-1}\hookrightarrow \proj(\vare_i)\to \classs \U(i)$ has the property that the restriction to the fiber $\K^*(\proj(\vare_i))\to\K^*(\cpspace^{i-1})$ is surjective.
\end{lemma}
\begin{proof}
This is a corollary of Grothendieck's introduction of the Chern classes.\\
In fact, let $\mathcal{O}(1)\to\cpspace^n$ be the dual of the tautological bundle. There is the fact that $\K^*(\cpspace^n)$ is $\Z[X]/(X^{n+1})$ where $X$ is a generator of degree $0$ which is given by $X=[\mathcal{O}(1)]-[1]$ (In particular, $\K^1(\cpspace^n)$ vanishes).\\
Now we generalize this to the general bundle. Denote by $\mathcal{O}_\vare(1)\to \proj(\vare_i)$ the dual of the tautological bundle over $\proj(\vare_i)$. Let $i\colon \cpspace^n\to\proj(\vare_i)$ be the inclusion. Now we have that $i^*(\mathcal{O}_\vare(1))=\mathcal{O}(1)$. It follows that $i^*$ is surjective in $\K$-theory.\\
Using the Atiyah-Hirzebruch spectral sequence, one could also make use of a similar reasoning to cohomology, where one may use the Serre spectral sequence. In fact, the cohomology of $\classs\U(n)$ is concentrated in odd degree, as is that of $\cpspace^{i-1}$. So there is no room for non-trivial differentials and the edge homomorphisms must be injective resp. surjective, which proves the same.
\end{proof}
Since any vector bundle $E\to Y$ is the pull-back of the universal bundle, we can conclude that the particular nature of the base space does not play any role. The following proposition is standard and proved in the same way as the above lemma.
\begin{prop}
$\K^*(\proj(E))$ is a \emph{free} $\K^*(Y)$-module with basis $[1],X_E,\dots,(X_E)^n$ for an \emph{arbitrary} vector bundle $E\to Y$, where $X_E$ is the class of the tautological line bundle over $\proj(E)$ minus $[1]$. Here, we have to understand $\K$-theory with possibly non-compact support, i. e. the homotopy invariant version.\\
In particular $\pi^*$ is split-injective.
\end{prop}
\begin{cor}
$\pi^*\colon \K^*(\BG)\to\K^*(\EG\times_\Gamma G/B)$ is split injective for $\K$-theory with and without support condition.
\end{cor}
\begin{proof}
We have shown this for each of the single composing steps of the fibration $\classs B\to\classs G$. So the split injectivity holds for $\classs B\to\classs G$ and by the same argument as in the preceding sections also for $\EG\times_\Gamma G/B\to\BG$. The statement for compact support $\K$-theory (which interests us the most) is the content of the following lemma \ref{lfinjective}(2).
\end{proof}
\begin{lemma}\label{lfinjective}\mbox{ }\\
\vspace{-1.5em}
\begin{enumerate}
\item
If $\pi^*\colon \Homol^*(B)\to\Homol^*(E)$ for an orientable fibration $F\hookrightarrow E\to B$ of smooth manifolds with compact fiber but not necessarily compact base space is injective, then so is the map $\pi^*\colon \Homolc^*(B)\to\Homolc^*(E)$ between the cohomologies with compact support. Dually, if $\pi_*$ is surjective between ordinary homologies, then so is $\pi_*$ between the locally finite homologies.
\item
A $\K$-theoretical analogue to the cohomological part of the preceding statement holds. We may then also replace ``split injective'' for ``injective'' everywhere.
\end{enumerate}
\end{lemma}
\begin{proof}\mbox{ }\\
\vspace{-1.5em}
\begin{enumerate}
\item Choose an exhaustion $\oslash\subset K_0\subset K_1\subset\dots\subset B$ of the base by compact subsets which are smooth connected manifolds with boundary.\\
Now the point is that we know the splitting principle for the fibrations $\pi^{-1}K_i\to K_i$ and $\pi^{-1}(B-K_i)\to B-K_i$ for all $i$ because all we have used in the proof is that the bundle is the flag bundle of a complex vector bundle. There have not entered any particular properties of the base space.\\
Denote by $E_i=\pi\inv(K_i)$. So $\Homol^*(E_i)$ is a free module over $\Homol^*(K_i)$, and  $\Homol^*(E-E_i)$ is a free module over $\Homol^*(B-K_i)$. Moreover, the maps induced by the inclusions $K_i\to K_{i+1}$ and $E_i\to E_{i+1}$ map the canonical basis into the canonical basis. By the excision theorem, $\Homol^*(B,B-K_{i})\cong\Homol^*(K_{i+1},K_{i+1}-K_i)$ and $\Homol^*(E,E-E_{i})\cong\Homol^*(E_{i+1},E_{i+1}-E_i)$ because we may assume that the closure of $B-K_{i+1}$ is contained in the interior of $B-K_i$. Consider the long exact cohomology sequence at the $i$-th step.
$$
\normalsize{
\xymatrix{\dots\ar[r]&\Homol^{k-1}\ar[d]\ar[r]&\Homol^k(E_{i+1},E_{i+1}-E_i)\ar[d]\ar[r]&\Homol^k(E_{i+1})\ar[d]\ar[r]&\Homol^k(E_{i+1}-E_{i})\ar[d]\ar[r]&\dots\\
\dots\ar[r]&\Homol^{k-1}\ar[r]&\Homol^k(K_{i+1},K_{i+1}-K_i)\ar[r]&\Homol^k(K_{i+1})\ar[r]&\Homol^k(K_{i+1}-K_i)\ar[r]&\dots
}
}
$$
Let $\xi^j$ denote the canonical $\Homol^*(X)$-base of the non-relative groups $\Homol^*(\pi\inv X)$ for $X=K_i$ or $X=K_{i+1}$. Clearly these bases correspond to each other under the map $i$ induced by the inclusion so we will not distinguish them by notation.\\
Furthermore, the relative group is also a module under $\Homol^*(E)$, and so there is a multiplication by $\xi^j$ in $\Homol^*(K_{i+1},K_{i+1}-K_i)$.\\
It is a general fact then that the boundary operator acts as a derivation: 
$$d(z\xi)=d(z)\xi+(-1)^{|z|}zd\xi.$$
Since $\xi$ is, as we have said, in the image of the preceding group, $d\xi=0$. Thus all maps of the long exact sequence are linear with respect to multiplication by $\xi$. There is a map $(\Homol^*(K_{i+1},K_{i+1}-K_i))^{\oplus n}\to \Homol^*(E_{i+1},E_{i+1}-E_i))$ defined by sending a family $(x_j)$ to $\sum \pi^*(x_j)\xi^j$. We want to prove that this map is an isomorphism. For surjectivity, let $\alpha\in\Homol^k(E_{i+1},E_{i+1}-E_i)$ be arbitrary. Consider its image $i(\alpha)\in\Homol^k(E_i+1)$ in the long exact sequence. Since $\Homol^*(E_{i+1})$ is a free module over $\Homol^*(K_{i+1})$, there is a unique expression $i(\alpha)=\sum_0^{n-1} \pi^*(x_j)\xi^j$ for $x_j\in\Homol^{k-2j}(K_{i+1})$. Since the sequence is exact, $0=p(i(\alpha))=\sum_0^{n-1} p(x_j)\xi^j\in\Homol^k(E_{i+1}-E_i)$ and we see that all $p(x_j)=0$. So all $x_j$'s are in the image of $i$, say $i(y_j)=x_j$ for $y_j\in\Homol^k(K_{i+1},K_{i+1}-K_i)$, and $\alpha-\sum_0^{n-1} y_j\xi^j\in\Ker\colon  \Homol^k(E_{i+1}),E_{i+1}-E_i)\to\Homol^k(E_{i+1})$. This means that $\alpha-\sum y_j\xi^j=\delta(\beta)$ for $\beta\in\Homol^{k-1}(E_{i+1}-E_i)$. Now we may express $\beta$ in the basis as $\beta=\sum \pi^*(z_j)\xi^j$ for $z_j\in\Homol^{k-1-2j}(K_{i+1}-K_i)$. In other words, $\alpha=\sum \pi^*(y_j+\delta(z_j))\xi^j$ which proves that there is an expression for $\alpha$ in terms of $(\xi^j)$ and the map is surjective.\\
To prove uniqueness, let $\sum y_j\xi^j=0$. So $\sum i(y_j)\xi^j=0$ which implies that $i(y_j)=0$. So $y_j=d(z_j)$ for some $z_j$'s. Then $\sum d(z_j)\xi^j=0=d(\sum z_j\xi^j)$ which means that $\sum z_j\xi^j=p(\sum q_j\xi^j)=\sum p(q_j)\xi^j$. So $z_j=p(q_j)$ and $y_j=d(p(q_j))=0$ by the exactness.\\
So we conclude that the relative group $\Homol^*(E,E-E_i)$ is a direct
sum of $\Homol^*(B,B-K_i)$. Moreover, the canonical splitting map $\Homol^*(B,B-K_i)\to \Homol^*(E,E-E_i)$ is given by $x\mapsto \pi^*(x)\xi^0$, in other words by the map $\pi^*$. By the naturality of the proof, this fact passes to the direct limit $\Homol_c^*(E)=\lim_i \Homol^*(E,E-E_i)$.\\
The statement for homology follows from the duality of $\pi^*$ and $\pi_*$.\\
Remark: An alternative proof could be given as follows. By explicit calculations with differential forms, one can see that the map which makes the following diagram commute
$$
\xymatrix{
\Homolc^*(E)\ar[r]^-\cong & \Homol_{\dim E-*}(E)\\
\Homolc^*(B)\ar[r]^-\cong \ar[u]^{\pi^*} & \Homol_{\dim B-*}(B)\ar@{.>}[u]
}
$$
can be read off explicitely in the Leray-Serre spectral sequence. Namely, it is the map which sends the homology to the limit of the uppermost row in the $E^2$-term of the spectral sequence for ordinary homology. This map is injective since all differentials vanish. However, we wrote down a purely topological proof to make it valid in K-theory as well.
\item
The above proof applies verbatim to the $\K-$theory case.
\end{enumerate}
\end{proof}
\subsection{The consequence for cyclic cohomology}
In the following theorem, we reach the goal of the chapters two and three in view of the sequel of the thesis.
\begin{thm}\label{essence}
Let $\Gamma$ be a torsion free, not necessarily cocompact discrete subgroup of a connected semisimple Lie group $G$ which has a Borel subgroup $B$ which satisfies the condition appearing in the conclusion of remark \ref{condition}. For instance, let $G=\SL(n,\C)$.\\
Then the map
$$
\HP^0\bigl(\CINF(B\backslash G)\alcross\Gamma\bigr)_{<1>}\longrightarrow\HP^0(\C\Gamma)_{<1>}
$$
induced from shrinking $B\backslash G$ to a point is surjective. In particular, the class
$$
[\tau]\in\HP^0(\C\Gamma)_{<1>}
$$
of the canonical trace
$$
\tau\colon \C\Gamma\to\C,\ \sum a_\gamma\gamma\mapsto a_e
$$
is in the image.
\end{thm}
The point about this theorem is that it shows that the trace pairing $<e,\tau>$, where $e$ is an idempotent in $\C\Gamma$, depends only on the value of $e$ in $\K_0(\Gamma\allcross\CINF(G/B))$ resp. $\HP_0(\Gamma\allcross\CINF(G/B))$. From that fact, and the knowledge of the K-theory $\K_0(\Gamma\rlcross C(G/B))$ we will develop the proof of the idempotent theorem.
\begin{proof} Let us apply all we have done so far to the case $V=B\backslash G,\ W=\pt,\ p\colon V\to W$.
By making explicit the functoriality for the map $p$ of the cyclic cohomology under all the isomorphism we have established so far, namely
\begin{align*}
\HP^*(\CINF(V)\alcross\Gamma)_{<1>}&\cong\mathbb{H}^*(\Gamma; \Omega_\ev V\lra\Omega_\odd V)\cong\mathbb{H}^{*+\dim V}(\Gamma; \Omega^\ev V\lra\Omega^\odd V)\\
&\cong\bigoplus_j\Homol^{*+\dim V+2j}(\EG\times_\Gamma V)
\end{align*}
we see that the following diagram commutes:
$$
\xymatrix{
{ \Homol^{\ev}(\BG)}\ar[r]^{\cong}&{\HP^0(\C\Gamma)_{<1>}}\\
{\Homol^{\ev+ \dim B\backslash G}((B\backslash G)\times_\Gamma\EG)}\ar[u]^{p_!}\ar[r]^{\cong}&{\HP^0(\CINF(B\backslash G)\alcross\Gamma)_{<1>},}\ar[u]_{p_*}
}
$$
where the left hand vertical arrow is the shriek map from definition \ref{shriek}.
Let $\pi\colon E\to B$ be an orientable fibration (i. e base, fiber and total space are coherently oriented). We want to relate the fiber integration $\pi_!$, which lowers the degree of a cohomology class by the fiber dimension, with the ordinary degree-preserving $\pi^*$. We therefore show that $\pi_!$ deserves its name, namely that it equals the composition $\Homol^*(E)\cong\Homollf_{\dim E-*(E)}\to\Homollf_{\dim E-*}(B)\cong\Homol^{*-\dim*(B)-\dim E}(B)$ which involves Poincar\'e duality for the whole space. Here, $\Homollf$ stands for locally finite homology (\cite{BottTu}). Note that $\Gamma$ acts by orientation preserving diffeomorphisms on $\EG=G/K$ and on $B\backslash G$, so the orientation condition is satisfied.\\
The map $\pi_!$ is induced by the integration over the fiber-map, which we shall also denote by
$$\pi_!\colon \Omega^i E\to \Omega^{i-k}B,\ \omega\mapsto\int_{\pi^{-1}x}\omega,\ x\in B$$
where the fiber is compact orientable and of dimension $k$. If we compose this map with Poincar\'e duality in the base $B$, we obtain the de-Rham current
$$
\tau\mapsto\int_B\tau\wedge \pi_! \omega
$$
on $B$, where $\tau$ has to have compact support. By lemma \ref{fiberintegration}, this is the same as 
$$
\biggl(\tau\mapsto\int_E \pi^*\tau\wedge \omega\biggr)=\pi_*\biggl(\zeta\mapsto \int_E\zeta\wedge\omega\biggr).
$$
This shows that the following diagram commutes:
$$
\xymatrix@C=5em{
{\Homol^*(E)}\ar[r]^{\text{PD}}\ar[d]^{p_!}&{\Homollf_{\dim E-*}(E)}\ar[d]^{p_*}\\
{\Homol^{*+\dim B-\dim E}(B)}\ar[r]^{\text{PD}}&\Homollf_{\dim E-*}(B),
}$$
in other words, our map $p_!$ is the shriek map for the projection of the fibration. So if the $p_*$ map of the fibration between the locally finite homologies is surjective - as assured by lemma \ref{lfinjective}, then so is the $p_!$.\\
Note that in the proof we \emph{cannot} dispense with Poincar\'e duality for the base.
\end{proof}
The following corollary will not be used in the sequel, but it illustrates the way the theorem can prove useful, and it sets the theorem a little closer to the framework of the first chapter.
\begin{cor}
The class of unit $[1]\in\K_0(\Gamma\allcross\CINF(G/B))$ is not torsion.
\end{cor}
\begin{proof}
Let $\tau'$ be a cocyle for the preimage of the trace class, i. e. the restriction of $\tau'$ to $\C\Gamma$ is cohomologous to $\tau$. Then $<\tau',[1]>=<\tau,1>=1$, so there is a functional on $\K_0(\Gamma\allcross\CINF(G/B))$ that detects $[1]$.
\end{proof}
\clearpage{\pagestyle{empty}\cleardoublepage}
\chapter{The Dirac-dual Dirac method, equivariant Bott periodicity and index theory}
We now set out to obtain a higher index formula in our setting, i. e. we wish to dispose of a formula for the pairing
$$
\K_0\bigl(\Gamma\allcross\CINF(G/B)\bigr)\otimes\HP^0\bigl(\Gamma\allcross\CINF(G/B)\bigr)_{<1>}\longrightarrow\C.
$$
where we are interested only in the value on the constructed cyclic cocycle $\tau'\in\HP^0\bigl(\Gamma\allcross\allowbreak\CINF(G/B)\bigr)_{<1>}$ whose restriction $i^*\tau'$ to $\C\Gamma$ is cohomologous to the trace: $i^*\tau'=\tau$. Our goal is to calculate $<e,\tau>=<i_*e,\tau'>$.\\
Concretely, we will \emph{compute} $\HA_0(\Gamma\allcross\CINF(G/B))$ by a Dirac-dual Dirac method. See the synopsis for a short explanation why this is needed. We would have liked to compute $\K_0$ instead of $\HA_0$. However that turned out to be too difficult. Even the result for $\HA$ contains the far most involved arguments of the whole paper, and a computation of $\HA$ fortunately suffices for the proof of the idempotent theorem.\\
In the result it is visible which classes can come from $\HA_0(\mathcal{R}\Gamma)$.\\
See also the appendix \ref{foliation} for an overview of how and under what conditions Connes' foliation index theorem would yield the result.


\section{The $C^*$-algebraic Dirac-dual Dirac construction}
In this section, we shall explain how to compute $\K_*(\Gamma\rlcross C(G/B))$. Later on, we will mimic the procedure on smooth subalgebras in order to obtain a calculation of $\HA_*(\Gamma\allcross \CINF(G/B))$.\\
Let $q=\dim B/T=\dim G/K=\dim EG$, as throughout the thesis.\\
In the entire sequel we will assume that there exists a $\Gamma$-equivariant Spin$^c$-structure on $\EG=G/K$.
\begin{rmk}\label{SPINN2}
As we restrict to the case $G=\SL(n,\C)$, the Spin$^c$ assumption we made during the construction of the assembly map in that form is automatic, i. e. $\EG$ is always equipped with a $G$-equivariant Spin$^c$ structure.\\
More generally, whenever $G$ is simply connected, the condition is satisfied.
\end{rmk}
\begin{proof}
One proofs inductively that $\SU(n)$ is simply connected as follows. $\SU(2)$ is the three-sphere. $\SU(n)/\SU(n-1)$ is $\C\mathbb{P}^{n-1}$, as one can see by making $\SU(n)$ act transitively on the lines in $\C^n$. So the induction step is completed by a look at the long exact homotopy sequence associated to that fibration.\\
Now the only obstruction against the existence of a $G$-equivariant Spin-structure on a finite dimensional representation space of a compact Lie group $G$ is two-torsion in the fundamental group of $G$ (see \cite{Botttwotorsion}).
\end{proof}
\begin{prop}\label{Spinc}\label{SPINN1}
There is an invertible bivariant $G$-equivariant Kasparov element
$$
\beta\in\KK_q^G(C(G/B),C_0(G/T))
$$
such that the restriction $\beta'\defeq\res_G^1(\beta)\in\KK_q(C(G/B),C_0(G/T))$ is the usual Thom element for the bundle $G/T=KAN/T=K\times_T A\times N\to K/T=KAN/TAN=G/B$, where $G=KAN, B=TAN$ is the Iwasawa resp. Langlands decomposition of $G$ resp. $B$. $G$ acts by left multiplication on both $G/B$ and $G/T$.
\end{prop}
In case $G=\SL(n,\C)$, the subgroup $A$ of $G$ consists of the diagonal matrices with real positive entries and determinant 1, $N$ consists of the unipotent upper triangular matrices, i. e. of the matrices of the form $1+$ strictly upper triangular matrix, $B$ of the upper triangular matrices with determinant 1, and $T$ of the unitary diagonal matrices with determinant 1. $B$ is an extension
$$
1\to NA\to B\to T\to 1
$$
and the torus acts on $NA$ by conjugation (trivially on A), and $B$ is a semidirect product $(NA)\rtimes T$ in the sense of Lie groups.
\begin{proof}
The point is to show that the standard Thom element $\beta'$ is an equivariant cycle, i. e. it is in the image of the restriction homomorphism $\res_G^e$.\\
Let $X$ be a space gifted with a left $G$- and a commuting right $H$ action, $H$ being any locally compact group, such that the $H$-action is free and proper. Denote by $\Morita_{X,H}\in\KK_0^G(C_0(X)\rtimes H,C_0(X/H))$ ($G$ acts trivially on $H$) the class of the Kasparov triple given by the $C_0(X/H)$-Hilbert-$C^*$-module
$$
\Hilm^+=\biggl\lbrace f\in L^2(X):\int_H |f(-.h)|^2 dh\in C(X)^H=C(X/H)\biggr\rbrace,\ \Hilm^-=0
$$
i. e. the sections of the Hilbert space bundle whose Hilbert space over a point in $X/H$ is exactly the corresponding $H$-orbit on $X$ the obvious representation of $G$ and of $C_0(X)\rtimes H$ on $\Hilm^+$, and the zero operator. Then $\res_G^1(\Morita_{X,H})\in\KK_0(C_0(X)\rtimes H,C_0(X/H))$ is the standard invertible Kasparov triple providing the strong (Rieffel-)Morita equivalence between $C_0(X)\rtimes H$ and $C_0(X/H)$. The proof that $\Morita_{X,H}$ is invertible goes through as in the non-equivariant case (see\cite{Kas}, Theorem 2.18). Denote by $\Morita^{-1}_{X, H}$ its inverse.\\
The existence of a $T$-invariant Spin$^c$-structure for $B/T$ is automatic for the following reasons. It is also a consequence of our standing assumption that $\EG=G/K$ has a $K$-invariant Spin$^c$-structure. In fact, the existence of a $G$-invariant Spin$^c$-structure on $G/K$ is equivalent to saying that the adjoint representation $K\to\SO(\mathfrak{g/k})$ lifts to $K\to$Spin$^c\to\SO$. Now the representation of $T$ on $\mathfrak{b/t=g/k}$ is a restriction of the representation of $K$ and thus also lifts. Recall that our standing assumption is always satisfied if $G=\SL(n,\C)$ because then $K$ is simply connected. In the case of $B$ acting on $B/T$, one may also argue that representations of the torus are always Spin, because it is simply connected, (see remark \ref{SPINN1}) or also because representation of a torus my be diagonalized into complex representations on $\R^2=\C$.\\
Consider now the dual Dirac-element $\eta_{B/T}\in \KK^B_q(\C,\Cliff (T^*(B/T)\otimes\C))$ constructed by Kasparov. $\Cliff$ stands for the sections vanishing at infinity of the complex Clifford-algebra associated to a complex vector bundle. The $\Spin$ structure allows to view $\eta_{B/T}$ as an element in $\KK^B_q(\C,C_0(B/T))$. Let $G$ act trivially on $C_0(B/T)$. Then $\eta_{B/T}$ defines an element in the group $\KK_q^{G\times B}(\C,C_0(B/T))$. Denote by $\tau_A$ the homomorphism $\KK^G(C,D)\to\KK^G(C\otimes A,D\otimes A)$ that tensors the Hilbert module, representation and the operator with the $G$-algebra $A$. Finally, denote by $j^H$ the descent homomorphism $\KK^{G\times H}(C,D)\to\KK^G(C\rcross H,D\rcross H)$, see \cite{Kas}, 3.11.\\
Then we are in a position to define $\beta$. Namely, we let $\beta$ be the Kasparov product
\begin{multline*}
\Morita_{G,B}^{-1}\ \otimes_{C_0(G)\rtimes B} j^B(\tau_{C_0(G)}(\eta_{B/T}))\ \otimes\\ \otimes_{C_0(G\times B/T)\rtimes B}\Morita_{G\times B/T,B}\in \KK_q^G(C(G/B),C_0(G/T)).
\end{multline*}
So we view the spaces $G/B$ resp. $G/T$ as the crossed-product $C^*$-algebras $C_0(G)\rtimes B$ resp. $C_0(G)\rtimes T$ and then use the Kasparov element (which is also the one appearing in the Connes-Kasparov-Rosenberg-conjecture) to reduce the crossed product with $B$ by a crossed product with its maximal compact subgroup $T$. This description already explains that $\beta'=\res_G^1(\beta)$ is invertible, since for solvable Lie groups such as $B$, the dual Dirac element $\eta_{B/T}$ is invertible with inverse the Dirac element. So, $\beta'=\res_G^1(\Morita)\otimes j^B(\tau_{C_0(G)}(\eta)\otimes\res_G^1(\Morita^{-1})$ is invertible as well. $\beta'$ is exactly the invertible element that induces the K-theoretic Thom isomorphism.
\end{proof}
As our goal is to represent $\beta$ by a homomorphism (such that $\beta$ restricts to a map between the smooth subalgebras), we need an element in $\KK_0$ rather than in $\KK_q$. We solve this by simply adding a one-dimensional trivial bundle with trivial actions.\\
More precisely, in case $q$ is odd, there is an invertible element
$$
\beta\in \KK_0^G(C(G/B),C_0(G/T\times \R),
$$
$\R$ being given the trivial $G$-action. In fact, there is a canonical inver-tible element in $\KK^G_1(\C,C_0(\R))$, namely the generator of $\K^1(\mathbb{S}^1)$. We take the outer product of this cycle with the $\beta$ constructed before. So let $V$ be the manifold $\R$ if $q$ is odd, the one-point-space otherwise, equipped with trivial $G$-action. Then consider $\beta\in\KK_0^G(C(G/B),C_0(G/T\times V))$.\\
\section{The smooth Dirac-dual Dirac method}
In this section, we will describe how to find an analogue in $\HA$ of the procedure described in the preceding section.
\subsection{Split-exactness and Quasi-Homomorphisms}
In this subsection, we introduce the concept of quasi-homomorphism due to Cuntz. In (\cite{Cuntznewlook}), he showed that $\KK_0(A, B)$ is equal to the set of homotopy classes of homomorphisms from $\q A$ to $B\otimes\Ko$. The algebra $\q A$ is an algebra associated functorially to $A$ in such a way that any quasi-homomorphism, being a seemingly only slightly more general notion of homomorphism, from $A$ to $C$ is the same as a homomorphism from $\q A$ to $C$. So every KK-element can be represented by a quasi-homomorphism from $A$ to $B\otimes\Ko$. Our goal is to represent the Dirac and dual Dirac-elements as quasi-homomorphisms, because this allows to check if the smooth sub-algebras are preserved.\\
This means that we now want to change back from the $C^*$-algebraic context to the far more algebraic one of bornological algebras, i. e. we will again consider $\K_0$ and $\HP$ as functors on the categories of bornological algebras, where $\K_0$, though, simply forgets the bornology and looks at stable equivalence classes of idempotents. Algebraic $\K$-theory is in some aspects a much more subtle theory than the $C^*$-algebraic one. For example, it is \emph{not} invariant under tensoring with the compact operators of a separable Hilbert space, nor is it homotopy-invariant, nor even invariant under polynomial homotopies (\cite{Rosenberg}, exercise 3. 2. 24).\\
However, it is a fundamental theorem of algebraic $\K$-theory that it is a \emph{split-exact} functor. In other words, each split-exact sequence $0\to J\to A\to A/J\to 0$ of algebras induces an exact sequence $0\to K_0(J)\to \K_0(A)\to \K_0(A/J)\to 0$ (which therefore also splits). The $C^*$-algebraic proof relies on excision and the 6-term sequence, and the algebraic proof is similar, so it can not be trivial since it involves algebraic $\K_1$, the algebraic boundary map and algebraic excision.\\
However, split-exactness suffices to define a functoriality under certain smooth $\KK$-elements. More precisely, we shall define quasi-homomorphisms of algebras and we shall show how they induce $\K$-theory maps.\\
What we will say will $\emph{a fortiori}$ apply to the split-exact functors $\HP$ and $\HA$ as well.\\
The content of the following lemmas is taken from \cite{MEYERKTHEORY}.
\begin{defn}
Let $A, B, D, D'$ be bornological algebras. Let $0\to B\to D'\to A\to A$ be a sequence of bornological algebras, i. e. the maps are bounded algebra homomorphisms. It is said to be an extension if it is exact as vector spaces, $B$ carries the subspace bornology of $D$, and the induced map $D'/B\to A$ is a bornological isomorphism.\\
A \emph{special quasi-homomorphism} from $A$ to $B$ is an extension $0\to B\to D'\to A$ together with two splittings  $f_+, f_-:A\to D'$ (i. e. bounded algebra homomorphisms such that $pf_+=pf_-=1_A$) such that the map $f_+-f_-$ is also bounded and its image is contained in that of $B$.\\
A \emph{quasi-homomorphism} is an algebra $D$ together with a closed ideal $B$ (i. e. $B$ has the subspace bornology) and two bounded algebra homomorphisms $f_\pm:A\to D$ such that $f_+-f_-$ is bounded and has image in $B$. One no longer requires that $D/B$ be isomorphic to $A$.\\
Sometimes one writes a quasi-homomorphism as $f_\pm:A\rightrightarrows D\rhd B$.\\
A special quasi-homomorphism $f_\pm:A\rightrightarrows D'\rhd B$ induces a map $\K_0(A)\to \K_0(B)$ and a map $\HA_0(A)\to \HA_0(B)$ as follows.\\
We have $p_*((f_+)_*-(f_-)_*)=p_*(f_+)_*-p_*(f_-)_*=\id_*-\id_*=0$, so by the split-exactness $\Image ((f_+)_*-(f_-)_*)\subset \K_0(B)\subset \K_0(D').$ This means that one may define $(f_\pm)_*$ as $(f_+)_*-(f_-)_*$.\\
Let now $f_\pm:A\rightrightarrows D\rhd B$ not be special. Then $D':=A\oplus B$ with the multiplication
$$
(a, b)(a', b')=(a'+f_+(a)b'+bf_+(a'), bb')
$$
defines a bornological algebra $D'$ (one may write $D'=A\ltimes B$ with the action of $A$ via $f_+$ on $B$ as an ideal inside $D$ is understood) together with two maps $f_+',f_-':A\to D'$, defined as
$$
f_+'(a)=(0,a),\ f_-'(a)=(f_+(a)-f_-(a),a)
$$
which are sections. Then one defines $(f_\pm)_*=(f_\pm')_*$.\\
Since $\HP_0$ is also split-exact in each variable (now the splitting is required to be bounded since we consider cyclic cohomology associated to bornological algebras), the definition makes perfectly sense for cyclic homology. We have thus defined natural transformations
$$
(f_\pm)^*:\HP(B,-)\to\HP(A,-)
$$
and
$$
(f_\pm)_*:\HP(-,A)\to\HP(-,B).
$$
Moreover, these come from a bivariant element $f_\pm\in\HP(A,B)$ since $\HP(A,D')=\HP(A,B)\oplus \HP(A,A)$ and $f_+-f_-\in\HP(A,B)\subset \HP(A,D')$.\\
Furthermore, the fact that $(f_+)^*$ and $(f_-)^*$ in cyclic cohomology are dual to $(f_+)_*$ and $(f_-)_*$ in $\K$-theory implies that $(f_\pm)^*$ is dual to $(f_\pm)_*$.
\end{defn}
We now want to describe ways of producing quasi-homomorphisms out of quasi-homo-morphisms. First, one may tensor a special quasi-homomorphism $f_\pm:A\rightrightarrows D\rhd B$ with a bornological algebra $F$ (it need not be special if the algebra $F$ is nuclear, anyway) to obtain the special quasi-homomorphism $f_\pm\otimes 1:A\otimes F\rightrightarrows D\otimes F\rhd B\otimes F$.\\
Second we want to describe the operation of taking crossed products of quasi-homomor-phisms with a group $G$. We will discuss first the case where the quasi-homomorphism $f_\pm:A\rightrightarrows D\rhd B$ is made of algebras of the form (finite matrices over) $\CCINF(M)$ where $M$ is equipped with a smooth action of the connected Lie group $G$. First we define this operation on the level of bornological vector spaces as $\CCINF(G)\protens -$, the completed bornological tensor product with the algebra $\CCINF(G)$, endowed with the precompact (=bounded) bornology. Then since $\CCINF(G)\protens\CCINF(M)\cong\CCINF(G\times M)$, which is the vector space of the smooth convolution algebra of the groupoid $G\ltimes M$, the split short exact sequence $0\to A\to D'\to B\to 0$ is in fact mapped to a split short exact sequence of algebras $0\to G\ltimes A\to G\ltimes D'\to G\ltimes B\to 0$ by the distributivity of $\protens$ and $\oplus$. This proves
\begin{lemma}
A special quasi-homomorphism $f_\pm:A\rightrightarrows D\rhd B$ of finite matrix algebras of smooth function algebras of smooth $G$-manifolds induces a special quasi-homomorphism $G\ltimes f_\pm:G\ltimes A\rightrightarrows G\ltimes D'\rhd G\ltimes B$.
\end{lemma}
The condition of this lemma, of course, does not seem very beautiful. However, we wish to reach our goal without getting involved in too many functional-analytic complications.
\begin{rmk}
Forgetting bornologies, one may in a straightforward manner associate to the $G$-equivariant and not necessarily special quasi-homomorphism $A\rightrightarrows D\rhd B$ of $G$-algebras, $G$ being a group, a quasi-homomorphism $G\ltimes A\rightrightarrows G\ltimes D\rhd G\ltimes B$. Algebraically, the associated special quasi-homomorphism $G\ltimes A\rightrightarrows (G\ltimes D)'\rhd G\ltimes B$ is then equal to $G\ltimes A\rightrightarrows G\ltimes D'\rhd G\ltimes B$. This justifies that we only take bornological crossed products of special quasi-homomorphisms.
\end{rmk}
The second way of taking crossed products is with discrete groups $\Gamma$. Here the situation is even simpler because on the bornological vector space level we merely need to tensor with $\C\Gamma$ endowed with the fine bornology, yielding the algebraic tensor product. We can thus associate to a not necessarily special equivariant quasi-homomorphism $A\rightrightarrows D\rhd B$ of $\Gamma$-algebras, where $\Gamma$ acts by bounded automorphisms, a quasi-homomorphism $\Gamma\ltimes A\rightrightarrows \Gamma\ltimes D\rhd \Gamma\ltimes B$.
\subsection{Some preliminaries on Cuntz' picture of $\KK$-elements}
\begin{defn}
In order to deal with homotopies of quasi-homomorphisms, we need to introduce the $\q$-picture of quasi-homomorphisms due to Cuntz (\cite{Cuntznewlook}).\\
Let $A$ be a complex algebra. $\Q A$ is defined as the free product $A*A$, i. e. the free algebra on the elements of two copies of $A$, divided by the relations coming from the multiplication in the two copies. As $\Q A$ satisfies the universal property that the data of an algebra map $\Q A\to B$ is the same as the data of two homomorphisms $A\to B$, there is a map $\ Q A\to A$ coming from twice the identity $A\to A$. Let $\q A$ denote the kernel.\\
We put bornologies on these algebras as follows: Recall that the algebraic $QA\cong A*A$ is isomorphic to $A\oplus A\otimes A'\oplus A'\otimes A\oplus A\otimes A'\otimes A\oplus \dots$, a direct sum of alternating tensor products of the algebra $A$ and a copy $A'$. Consider the bornology generated by all sets of the form $S^1\times\dots\times S^n$ in an $n$-fold tensor product, where all $S^i$ are small in $A$ or $A'$. Its bornological completion is then $A\oplus A\protens A'\oplus\dots$, a direct sum of completed tensor products. So the tensors are completed, the direct sum isn't. This gives a complete algebra, because all $S^1\times\dots\times S^n$ are completant.
\end{defn}
In the sequel, we shall make use of the fact that $\HP$ and $\HA$ satisfy excision for extensions with a bounded linear section (for $\HP$, this is shown in \cite{excision}, for $\HA$ this is proved in \cite{Meyerthesis}).
\begin{prop}\mbox{ }\\
\vspace{-1.5em}
\begin{enumerate}
\item
There are inclusions $\iota_1,\iota_2:A\to \Q A$ such that given arbitrary bounded algebra homomorphisms $f_1,f_2:A\to D$ there is a unique bounded algebra homomorphism $f_1*f_2:\Q A\to D$ that makes the diagram
$$
\xymatrix{&A\ar[rd]_{\iota_1}\ar@<0.5 ex>[rd]^{\iota_2}\ar[ld]_{f_2}\ar@<0.5 ex>[ld]^{f_1}&\\
D&&\Q A\ar[ll]^{f_1 * f_2}
}
$$
commute.
\item
$\q A$ is functorially $\HA$- and $\HP$-equivalent to $A$.
\end{enumerate}
\end{prop}
The first part of the proposition says that $\Q A$ has in the complete bornological category the universal property that arbitrary pairs of algebra maps $f_1,f_2:A\to D$ correspond bijectively to bounded  algebra maps $f_1 * f_2:\Q A\to D$.
\begin{proof}\mbox{ }\\
\vspace{-1.5em}
\begin{enumerate}
\item
The map $$A^{\otimes j}\to D,a_1\otimes \dots \otimes a_j\mapsto f_1(a_1)f_2(a_2)f_1(a_3)\dots f_{1+(j\mod 2)}(a_j)
$$
is the composition of
$f_1\otimes\dots\otimes f_{1+(j\mod 2)}$ with algebra multiplication, thus bounded. So in fact there is $\Q A\to D$ with the required properties, which clearly is unique, too
\item
The proof is completely analogous to the proof that for a $C^*$-algebra $B$ there is an invertible element in $\KK(\q B,B)$. In fact, the homotopy used in \cite{Cuntznewlook} only consists in multiplying functions such as $\sin t$ or $\cos t$ on the algebra. They therefore merely require the use of an algebraic tensor product $A\to A\otimes\CINF[0,1]$, not of a completion. So we dispose of the required homotopy invariance in $\HA$. In short, the smoothness requirement for homotopy invariance is largely satisfied.\\
Let us write down the argument in our situation. Define a map from $\Q A$ to $A\oplus A$ by $\iota_1*\iota_2$ where $\iota_1$ and $\iota_2$ are the inclusions $A\to A\oplus A$, and a map from $A\oplus A\cong\diag(A,A)$ to $M_2(\Q A)$ by the inclusion, using the two embeddings $A\to \Q A$. Then the composition $\Q A\to A\oplus A\to M_2(\Q A)$ is homotopic to the upper-left-corner embedding $\Q A\to M_2(\Q A)$. In fact, take the path from $A\to \diag(0,A)$ to $A\to \diag(A,0)$. This gives a path between the two maps $\Q A\to M_2\Q A:\diag(\id *0,0*\id)$ and $\diag(\id*\id,0)$, the latter map being the upper-left-corner embedding. This proves that the map $A\oplus A\to \Q A$ is surjective in $\HA$ and in $\HP$. Second, by decomposing $M_2(A\oplus A)$ as $M_2(A)\oplus M_2(A)$ and rotating the latter entry to the first coordinate, one also sees that this map is injective.\\
By excision, this implies that $\q A$ is $\HA$-equivalent and $\HP$-equivalentto $A$. The existence of a bounded linear section is automatic.\\
Moreover, an equivalence is given by the map $\q A\to A$ defined by restriction of $\id*0$ to $\q A\subset \Q A$.
\end{enumerate}
\end{proof}
In the following lemma, the word ``functorial'' shall refer to functoriality in algebra maps $A\to B$.
\begin{lemma}\label{qpullinto}
Let $A$ be a complete bornological algebra endowed with a bounded action of a discrete group $\Gamma$. Then $\Q A$ and $\q A$ come with functorial actions as well. Furthermore, there is a functorial homomorphism $\theta: \q(A\alcross \Gamma)\to(\q A)\alcross\Gamma$ which is an $\HP$- and $\HA$-equivalence.
\end{lemma}
\begin{proof}
The first statement is clear from the functoriality of the construction of $\Q A$ and $\q A$.\\
The two canonical homomorphisms $A\to \Q A$ are $\Gamma$-equivariant, so they pass to two homomorphisms $A\alcross\Gamma \to (\Q A)\alcross\Gamma$. By the universal property of $\Q A$, this gives rise to a homomorphism $\widetilde{\vartheta}:\Q (A\alcross\Gamma)\to (\Q A)\alcross\Gamma$. Define $\vartheta$ as its restriction to $\q (A\alcross\Gamma)$.\\
Consider the morphism of short exact sequences
$$
\xymatrix{
 0 \ar[r]         & \q(A\alcross\Gamma)\ar[r]\ar[d]^{\vartheta}        & \Q(A\alcross\Gamma)\ar[r]\ar[d]^{\widetilde{\vartheta}}         & A\alcross\Gamma \ar[r]\ar@{=}[d]            & 0\\
 0 \ar[r]         & (\q A)\alcross\Gamma\ar[r]                       & (\Q A)\alcross\Gamma\ar[r]               & A\alcross\Gamma\ar[r]                                   & 0.
}
$$
By excision, it suffices to prove the statement that $\widetilde{\vartheta}$ is an $\HA$- and $\HP$-equivalence. Again, the existence of a bounded linear section is automatic.\\
It is a classical fact that the obvious homomorphism $A\oplus A\to \Q A$ is an $\HA$- and $\HP$-equivalence. The proof only uses functorially defined homotopies, so it passes to an analogous statement for the homomorphism $(A\oplus A)\alcross\Gamma\to(\Q A)\alcross\Gamma$. Consider the isomorphism $(A\alcross\Gamma)\oplus (A\alcross\Gamma)\to (A\oplus A)\alcross\Gamma$. The following diagram commutes
$$
  \xymatrix{
  \Q(A\alcross\Gamma)\ar[d]^{\widetilde{\vartheta}}   &   (A\alcross\Gamma)\oplus (A\alcross\Gamma)\ar[l]\ar[d]  \\
  (\Q A)\alcross\Gamma & (A\oplus A)\alcross\Gamma.\ar[l]
}
$$
So $\widetilde{\vartheta}$ has the desired property because all other three maps have it.
\end{proof}
\subsection{Algebraic and analytic Morita equivalences} \label{Morita}
We now want to state two lemmas concerning compatibilities of various isomorphisms induced by Morita equivalences. These frequently arise throughout this thesis, as it has already become clear in the parts concerning KK-theory, and they will become even more important in the sequel. The content of this subsection is basically a reformulation in our context of the notion of Morita equivalence developed in \cite{CuntzMorita}, as well as of the results contained therein and in \cite{SOMEWHERECUNTZ}.\\
The first lemma is about the invertible elements in $\HP$ and $\HA$ induced by the Morita equivalences, the second one assures that in the case the algebras are strongly equivalent $C^*$-algebras the isomorphism in analytic cyclic homology coming from the first lemma and the fact that $\HA=\HPL$ on $C^*$-algebras is compatible with the one coming from the invertible $\KK$-element as described by Kasparov, under the Chern character.\\
As always, we work in the category of complete bornological algebras.\\
Let us first recall how algebraic Morita equivalences are defined and how they induce isomorphisms in cyclic homology. Goodwillie's theorem states that, if
$$
0\to I\to A\to B\to 0
$$
is an extension with linear section, then
$$
\widehat{A}:=\varprojlim A/I^n\to A/I=B
$$
induces invertible elements in $\HP_*$ and $\HA$ (the original reference is \cite{Goodwillie}, in our context, a proof for $\HP$ is given in \cite{SOMEWHERECUNTZ}, for $\HA$ in \cite{Meyerthesis}*{3.1.7}). A Morita context for two algebras $A$ and $B$ is an algebra of the form 
$D=\left(\begin{smallmatrix}A&E\\F&B\end{smallmatrix}\right)$, i. e. an algebra $D$ with a splitting into four linear subspaces such that if the elements are written as matrices multiplication is compatible with matrix multiplication. A Morita context $D$ induces a Morita equivalence between $A$ and $B$ if $A\cong D(FE)$ i. e. if in $A$ is dense the linear subspace inside $D$ generated by all products $ef, e\in E, f\in F$ and similarly $B\cong D(EF)$.\\
For example there is a Morita context between $A$ and $A\otimes \ell^p$, where $\ell^p=\ell^p(\mathcal{H})$ is the algebra of compact operators on a Hilbert space $\mathcal{H}$ with a $p$-summable sequence of characteristic values, given by $A\otimes \ell^p(\mathcal{H}\oplus\C)$. The other and more important Morita equivalences that we will use throughout this thesis will all be based on the following form. Let $X$ be a possibly noncompact Riemannian manifold, and $G$ a Lie group right-acting on $X$ in a smooth, free and proper, but not necessarily cocompact way. Since we do not suppose $G$ to be connected, this also includes the case of a discrete group acting on $X$. Then the algebras $A=\CCINF(X)\alcross G$ and $B=\CCINF(X/G)$ fit into a Morita context $D=\left(\begin{smallmatrix}A&E\\F&B\end{smallmatrix}\right)$ where $E=F=\CCINF(X)$. The multiplication formulas are explicitly given in \cite{Puschnigg}, but there is a more conceptual way to see that $D$ is a Morita context. Recall that two $C^*$-algebras $N$ and $H$ are said to be (strongly) Morita-equivalent if there is a full Hilbert $H$-module $\Hilm$ such that $N$ is isomorphic to the $C^*$-algebras of compact operators $\Ko(\Hilm)$ on $\Hilm$ (\cite{Blackadar}*{13.7.1(b)}). This is in fact an equivalence relation. It follows that $H$ is Morita equivalent to the compact operators on its Hilbert module $\Hilm\oplus H$, the orthogonal direct sum of $\Hilm$ and $H$ viewed as a Hilbert module over itself. $\Ko(\Hilm\oplus H)$ decomposes as a direct sum
$$
\Ko(\Hilm\oplus H)\cong\Ko(\Hilm)\oplus \Ko(\Hilm,H)\oplus\Ko(H,\Hilm)\oplus \Ko(H)=\left(\begin{smallmatrix}N&\Ko(H,\Hilm) \\\Ko(\Hilm, H)&H\end{smallmatrix}\right).
$$
In the case where $H$ is the $C^*$-completion $C_0(X/G)$ of $\CCINF(X/G)$ and $N$ is the $C^*$-completion $C_0(X)\rcross G$ of $\CCINF(X)\alcross G$, the Hilbert module $\Hilm$ over $H=C_0(X/G)$ is the Hilbert module of continuous sections of the following Hilbert space bundle over $X/G$. The fiber over $xG\in X/G$ is the Hilbert space $L^2(xG)$ on the fiber $xG\subset X$ of $X\to X/G$ with respect to the measure given by the Riemannian metric. The action of $H$ on $\Hilm$. So $\Hilm$ is a set of (equivalence classes of) functions on $X$. One then has in fact that the compact operators on this Hilbert module is the set of continuous sections of the $C^*$-algebra bundle over $X/G$ whose fiber over $xG$ is the $C^*$-algebra of compact operators on the Hilbert space $L^2(xG)$ on the fiber. This $C^*$-algebra of sections readily identifies with the $C^*$-algebra $G=C_0(X)\rcross G$ by indicating covariant representations of the pair $(C_0(X),G)$ on each fiber in the straightforward way that $C_0(X)$ acts by pointwise muliplication on $\Hilm$, which is, as we have said, a set of funtions on $X$, and letting $G$ act on the functions on $X$ by right translations. (Caution: This is a representations by unitary operators in each fiber, not by compacts, but the image of the representation of the crossed product induced by the covariant pair is the compact operators in each fiber.) So we have found a Hilbert module $\Hilm$ over $C_0(X/G)$ such that $C_0(X)\rcross G\cong\Ko(\Hilm)$. Now define a bounded injective linear map
$$
D=\left(\begin{smallmatrix}A&E\\F&B\end{smallmatrix}\right)\to \Ko(\Hilm\oplus H)
$$
as follows. The maps $A\to G$ and $B\to H$ are the inclusions of the dense subalgebras into their $C^*$-completions. The map $X\to \Ko(B,\Hilm)$ is given by
$$
X\times B\ni(x,b)\mapsto xb\in\Hilm,
$$
the map $X\to \Ko(\Hilm,B)$ is given by
$$
X\times\Hilm\ni(x,\xi)\mapsto <x,\xi>\in B.
$$
After checking that this is an injection of $D$ into $\Ko(\Hilm\oplus H)$ and that the image is a subalgebra, this defines an algebra structure on $D$, and one readily checks that $D$ is a Morita context.\\
There are certain ways of manufacturing from a Morita context another one. If $D$ is a Morita context between $A$ and $B$, and $J$ is another complete bornological algebra, then $J\protens D$ is a Morita context between $J\protens A$ and $J\protens B$. Second, suppose in the example above that there is given a smooth left action of a discrete group $\Gamma$ on $X$ that commutes with the right action of $G$. Then one verifies without difficulties that $\Gamma\allcross D$ is a Morita context between $\Gamma\allcross A$ and $\Gamma\allcross D$. More generally, whenever $A$ and $B$ are Morita-equivalent $\Gamma$-algebras via a Morita context $D$ which is also a $\Gamma$-algebra such that the inclusions of $A$ and $B$ into $D$ are equivariant, then $\Gamma\allcross D$ is a Morita context between $\Gamma\allcross A$ and $\Gamma\allcross B$.\\
This finishes the explanation of the Morita contexts that induce the Morita equivalences used in this thesis.\\
Let us now turn attention to how and why a Morita context induces invertibles in $\HP$ and $\HA$. Apply Goodwillie's theorem to the extension
$$
0\to\Ker(T(E\otimes F))\to A)\to T(E\otimes F)\to A\to 0
$$
where $TK$ is the tensor algebra of a vector space $K$, and $E\otimes F$ is the algebraic tensor product over $\C$ of $E$ and $F$. The map $T(E\otimes F)$ comes from the obvious linear map $E\otimes F\to A$ by the universal property of the tensor algebra. Denote by $\widehat{T}(E\otimes F)$ the adic completion $\varprojlim T(E\otimes F)/(\Ker(T(E\otimes F)\to A))^n$ of $T(E\otimes F)$ with respect to the adic filtration defined by the kernel of this extension. Then $\HP$ and $\HA$ of $\widehat{T}(E\otimes F)$ are those of $A$ by Goodwillie's theorem.\\
Consider the operator $\phi$ of cyclic permutation of tensor products
$$
\phi:T(E\otimes F)\to T(F\otimes E): e_1\otimes f_1\dots e_k\otimes f_k\mapsto f_k\otimes e_1\otimes f_1\otimes\dots\otimes e_k.
$$
The next lemma involves the $X$-complex, as defined, for instance, in \cite{SOMEWHERECUNTZ}.
\begin{lemma} The map $\phi$ passes to an isomorphism $\widehat{\phi}:\widehat{T}(E\otimes F)\to\widehat{T}(F\otimes E)$, and furthermore to a homotopy equivalence $X(\widehat{T}(E\otimes F))\to X(\widehat{T}(F\otimes E))$. Consequently, there is an invertible element in $\HP_0(A, B)$.\\
Second, in the case of the Morita contexts given above, the same is valid in analytic cyclic homology.
\end{lemma}
\begin{proof}
  The first part of the lemma is proved in \cite{CuntzMorita}. The second part of the lemma is proved in \cite{Puschnigg} in the case of the Morita equivalence between $\CCINF(X)\alcross G$ and $\CCINF(X/G)$ and the equivalences derived from these by one of the operations described above. The Morita equivalence between $A$ and $A\protens\ell^p$ is the stability statement of $\HA$ and proved in \cite{Meyerthesis}*{3.2.3}.
\end{proof}
The next lemma involves the bivariant \emph{local cyclic homology} $\HPL_*$ of $C^*$-algebras, as defined by Michael Puschnigg (see \cite{Pulocal} and the introduction). The facts that $\HPL$ has a bivariant Chern-Connes character and that $\HPL_*\cong\HA_*$ on $C^*$-algebras give a very useful way to compute analytic cyclic homology of $C^*$-homology (and in many cases, it seems that this is the only one, since $\HA$ does probably not possess a bivariant Chern-Connes-character). It implies that $\HA_*$ coincides on strongly Morita equivalent $C^*$-algebras. We now want to show that this Morita equivalence isomorphism is compatible with the one given by the preceding lemma.
\begin{lemma}
Let $X$ be a smooth Riemannian manifold with a smooth, free and proper action of the not necessarily connected Lie group $G$. Then the algebraic Morita-equivalence $\CCINF(X)\alcross G\sim\CCINF(X/ G)$ is compatible in $\K$-theory with the $C^*$-algebraic strong Morita equivalence $C_0(X)\rtimes\Gamma\sim C_0(X/\Gamma)$ in the sense that the following diagram commutes.
$$
\xymatrix{
\HA_0(\CCINF(X)\alcross G)\ar[d]_\Morita^\cong\ar[r]^{i_*}&\HA_0(C_0(X)\rtimes G)\ar[d]_\Morita^{\cong}\\
\HA_0(\CCINF(X/ G))\ar[r]^{i_*}&\HA_0(C_0(X/ G)),
}
$$
where the right arrow is induced from the isomorphism with local cyclic homology and the invertible $\KK$-Morita-element as explained above.
\end{lemma}
\begin{proof}
We shall use the nomenclature and maps described in the introduction of this seubsection. By construction, there is a commutative diagram of inclusions
$$
\xymatrix{
N\ar[r]&\Ko(\Hilm\oplus H)&H\ar[l]\\
A\ar[r]\ar[u]&D\ar[u]&B\ar[l]\ar[u]
}
$$
So the claim is proved if we can show that the two arrows in the first line are invertible in $\KK$, that their composition induces the invertible Morita element in $\KK(N,H)$, and analogously in the lower line that both arrows are invertible in $\HA$ and that their composition induces the invertible Morita element in $\HA(A,B)$.\\
Let us abbreviate $\Ko(\Hilm,H)$ by $\Ko$. This $\Ko$ is itself Morita-equivalent to $H$, and the inclusion $H\to\Ko$ is inverse to the Morita element in $\KK(\Ko,H)$. This follows from the fact that the Morita element is represented by the $\Ko-H$-bimodule $\Hilm\oplus H$ and the zero operator, and the inclusion is represented by the $H-\Ko$-bimodule $\Ko$. So the composition over $\Ko$ is represented by the $H-H$-bimodule $\Ko\otimes_{\Ko}(\Hilm\oplus H)\cong \Hilm\oplus H$, where both $H$ act by zero on $\Hilm$. This describes the sum of $1\in\KK(H,H)$ with a degenrate cycle, in other words $1\in\KK(H,H)$. Since we already know that the Morita element is invertible, the inclusion must be its inverse. A similar calculation applies to the other arrow of the first line, and so their composition describes the Morita element of $\KK(A,B)$ by the transitivity property of the construction.\\
The statement for the second line follows from the similar fact that $D$ is Morita equivalent to $A$ because there can be written down an obvious $3\times 3$-Morita context linking them, and the equivalence is induced by the corner embeddings.
\end{proof}

\subsection{The smooth Dirac-dual Dirac elements}
Let W be a finite-dimensional representation space of a compact torus $T$. We shall now show that there is also a smooth version of the Dirac and dual Dirac-elements in $\KK^T_{\dim W}(\C,C_0(W))$ resp. $\KK^T_{\dim W}(C_0(W),\C)$, represented by quasi-homomorphisms

\begin{align*}
&\beta :\C\rightrightarrows\CINF(W)\otimes M_2\rhd\CCINF(W)\otimes M_2,\\
&\alpha :\CCINF(W)\rightrightarrows \mathcal{L}(L^2 W)\rhd\ell^p(L^2 W)\otimes M_2
\end{align*}
For $W$ we will take $B/T\times V$ (recall that $B/T$ is $T$-equivariantly diffeomorphic to $\mathfrak{b/t}$, see \cite{Kas}), where $V$ is the manifold $\R$ if $q=\dim B-\dim T$ is odd, the one-point-space otherwise. $V$ is equipped with trivial $G$- and $B$-actions. We introduce this $V$ in order to assure that we are dealing with an \emph{even} KK-cycle.\\
The actions of $B$ on $B/T$ used to define the crossed products are \emph{not} the ``right'' ones. In fact they are defined in such a way that they factor through $T$.

\subsubsection{The $T$-equivariant dual-Dirac element for $B$}
Consider the classical equivariant Bott element for the even-dimensional vector space $B/T\times V$ associated to the chosen $T$-equivariant Spin$^c$-structure as defined by Atiyah (\cite{Atiyah Bottpe}). It is an element in $\K_T^0(B/T\times V)$, the equivariant $\K$-theory with compact support of $B/T\times V$ and can thus be represented as (and is in fact most easily defined as) a $T$-equivariant complex $(E^\bullet, \partial^\bullet)$ of hermitian vector bundles.
We may assume that the complex is of length one (\cite{ATIYAHKTHEORY})
$$
\xymatrix{0\ar[r]& E^0\ar[r]^{\phi}& E^1\ar[r]&0}.
$$
More precisely, consider the equivariant even and odd spinors such that the bundle of Clifford algebras of $T^*(B/T\times V)$, the cotangent vector bundle, is the endomorphism bundle of the spinors. We may choose $E^0$ to be the even equivariant spinors, $E^1$ as the odd ones and the operator between both to be multiplication with the radial vector field. This means consider the distance function $d(0,-)$ associated to the ($T$-invariant) Riemannian metric. Then the radial vector field is the differential of the function $\sqrt{1+d(0,-)^2}$, and it can be viewed as a section of the Clifford algebra bundle, acting on the spinors by Clifford multiplication.\\
We now want to translate this description of the Bott element, following Atiyah, into the idempotent picture of $\K$-theory. To this end, consider the hermitian vector bundle $E^0\oplus E^1$ over $B/T\times V$ and the subbundle given by the graph $\{(v\oplus v'):\pi(v)=\pi(v'), \phi(v)=v'\}$. There is in the endomorphism bundle of $E^0\oplus E^1$ a canonical section given in each fiber by orthogonal projection onto the graph. It is a family $e_x, x\in B/T\times V$ of idempotents $e_x\in M_n\CINF(B/T\times V)$, called the graph projection. This is now a $T$-equivariant element, and one of the elements defining the same $\K$-theory element as the equivariant Bott element. It remains to subtract the formal dimension of the idempotent. In order to give the difference compact supports, we proceed as follows. The procedure is described in \cite{Puasymptotic}, and we will give a short synopsis here. Denote by $e_\infty$ the constant family given by the idempotent $\lim_{x\to\infty} e_x$ of the family $(e_x)$ at infinity. We want the difference $(e_x)-e_\infty$ to have compact support. Since this only holds approximately, we have to cut down $(e_x)-e_\infty$ to zero by a function $\chi $ taking the value 1 on the compact set and zero outside a larger compact set. The resulting matrix-valued function is an idempotent inside the small compact set and outside the large one. However, the difference to an idempotent in the annulus between the two is small. This problem is solved by application of functional calculus with the function
$$
f:x\mapsto x+\sum_{k=0}^\infty \bigl(x-\frac{1}{2}\bigr)\bigl(x-x^2\bigr)^k
$$
known from the formula for the Chern character with values in the $X$-complex, applied to the function $\chi \times (e_x)$. In fact, the smallness of the difference between the two idempotents implies that the power series converges. If an idempotent is plugged in, the function is the identity. This implies that $f(\chi (e_x))$ is an idempotent everywhere. So the function $(e_x)-e_\infty$ has compact support. Moreover, the $T$-equivariance is preserved by elementary reasons by functional calculus.\\
So define 
$$
\beta^T_+:\C\to M_n(\CINF(B/T\times V))
$$
by the family $f(\chi(e_x))$ of idempotents, and similarly $\beta^T_-$ by the constant family $(e_\infty)$.
Thus the map $\C\to M_n\CINF(B/T\times V)$ defined by $\beta^T_+-\beta^T_-$ takes its values in $M_n\CCINF(B/T\times V)$. So by construction the beta quasi-homomorphism
$$
\beta^T:\C\to M_n\CINF(B/T\times V)\rhd M_n\CCINF(B/T\times V)$$
is defined and $T$-equivariant.
\subsubsection{The $B$-equivariant Dirac element for $B$}
We are first going to define a smooth version of the invertible element in $\KK_0^B(C_0(B/T\times V),\C)$.\\
In contrast to the dual-Dirac element, the Dirac operator is equivariant with respect to the entire non-compact group $B$.\\
Recall that $B/T\times V$ is endowed with a $T$-invariant Riemannian metric, and that we may assume that there exists a $T$-equivariant Spin$^c$-structure on $B/T\times V$. Denote by $\Spinor^+$ and $\Spinor^-$ the bundles of positive resp. negative spinors, and by $\partialdirac$ the Dirac operator from sections of $\Spinor^+$ to $\Spinor^-$. Consider the operator $\dirac=\left(\begin{smallmatrix} 0&\partialdirac^*\\ \partialdirac&0\end{smallmatrix}\right)$, being an endomorphism of the sections of the direct sum $\Spinor^\pm$ of positive and negative spinors.\\
The spinors are a hermitian vector bundle, so we may consider the Hilbert space $\Hil=L^2(B/T\times V, \Spinor^\pm)$ of square-integrable sections of $\Spinor^\pm$. $\Hil$ comes with a canonical representation $\rho$ of $C_0(B/T\times V)$ by pointwise multiplication. $\Spinor^\pm$ is a $B$-equivariant vector bundle since $B$ acts from the left on $B/T$, and thus on $B/T\times V$, by translations. Thus $B$ also acts on $\Hil$. Furthermore, denote by $\epsilon:\Hil\to\Hil$ the grading operator coming from the grading of the spinors. The expression
$$
F=\dfrac{\dirac+\epsilon\cdot 1}{\sqrt{(1+\dirac^2)}},
$$
is well-defined as the composition of (the commuting operators) $\dirac+\epsilon\cdot 1$ and the result of functional calculus of $\frac{1}{\sqrt{(1+x^2)}}$ with the essentially self-adjoint unbounded operator $\left(\begin{smallmatrix}0&\partialdirac^*\\ \partialdirac&0\end{smallmatrix}\right)$. Using that $\epsilon$ and $\dirac$ anticommute, one easily checks that $F^2=1$. Consider the two expressions
$$
\begin{array}{cccc}
&C_0(B/T\times V)&\rightarrow&\mathcal{L}(\Hil)\\[0.2cm]
\rho_+:&f&\mapsto&\frac{1+\epsilon}{2}\rho(f)\frac{1+\epsilon}{2}\\[0.2cm]
\rho_-:&f&\mapsto&F\frac{1-\epsilon}{2}\rho(f)\frac{1-\epsilon}{2}F.
\end{array}
$$
These are bounded algebra homomorphisms, given that $\rho$ is grading-preserving and therefore commutes with the projections. The following proposition is based on a classical result of Weyl.
\begin{lemma}
For every real $p>\dim B/T$ such that the difference $\rho_+-\rho_-:C_0(B/T\times V)\to\Ko(\Hil)$ restricted to the dense subalgebra $\CCINF(B/T\times V)$ takes its values in the Schatten ideal $\ell^p(\Hil)$ of compact operators with $p$-summable characteristic values.
\end{lemma}
\begin{proof} $F$ is a compact perturbation of an odd operator $\frac{\dirac}{\sqrt{1+\dirac^2}}$ (see \cite{CoNCDG}).
\end{proof}
This is very useful to us, because the periodic cyclic cohomology of $\ell^p$ (see \cite{SOMEWHERECUNTZ}) is not degenerate, in contrast to that of $\Ko$.\\
The lemma implies that the pair of morphisms $\rho_+$ and $\rho_-$ defines a quasi-homomor-phism $$\alpha^B=\rho_\pm:\CCINF(B/T\times V)\rightrightarrows\mathcal{L}(\Hil)\rhd\ell^p(\Hil),$$ which by the homogeneity of the construction is $B$-equivariant. ($B$ acts on $\Hil=L^2(B/T\times V, \Spinor^\pm)$ by conjugation)
\subsection{The proof of equivariant Bott periodicity in analytic cyclic homology}
We now want to prove that $\alpha$ and $\beta$ are $T$-equivariantly inverse to each other. More precisely , we want to homotop between $\alpha\circ\beta$ and the identity resp. $\beta\circ\alpha$ and the identity by a path of \emph{equivariant} homomorphisms. A second and quite different proof, relying on the realization of the Dirac morphism as a deformation, is given in the appendix.
\begin{thm}[Equivariant Bott periodicity, smooth version]\label{abhomotopy}\label{XXX}
Let $W$ be an even-dimensional Euclidean vector space. Consider the composition
$$
M_n\alpha\circ\q\beta:\q^2\C\to\q(M_n(\CCINF(W)))\to\ell^p W\otimes M_n
$$
of $M_n\alpha$ with $\q\beta$.\\
Consider furthermore the canonical homomorphisms $\q^2\C\to\q\C\to\C\to\ell^p$ coming twice from $(\id*0:\Q\C\to\C)|_{\q\C}$ (which induce, as we have shown, isomorphisms) and the canonical projector onto a $L^2$-solution of the Dirac operator.\\
Then the following diagram commutes up to smooth homotopy.
$$
\xymatrix{
\q^2\C\ar[rr]^{\alpha\circ\q\beta}\ar[rd]&&\ell^p\\
&\C\ar[ur]
}
$$
Furthermore, suppose $W$ is equipped with an isometric action of the compact Lie group $T$ such that our standing assumption that the representation $T\to\SO(W)$ lift to $\Spin^c$ is satisfied. Let $T$ act on $L^2=L^2 W$ as left translation operators and on $\ell^p=\ell^p L^2 W$ by conjugation.\\
Suppose the map $\C\to\ell^p$ is given by the projection onto a $T$-invariant solution of the Dirac operator. Then the diagram commutes up to smooth homotopy through a path of $T$-invariant homomorphisms.
\end{thm}
To prepare for the proof, we state the following
\begin{lemma}\label{yyy}
Let $e$ and $f$ be projections in $\mathcal{L(H)}$, such that $e-f$ is compact. Then $e$ is smoothly homotopic to a projection $e'$, such that the intersection of $e'\mathcal{H}$ and $f\mathcal{H}$ has finite codimension in $e'\mathcal{H}$ as well as in $f\mathcal{H}$. If $e$ and $f$ are invariant under a linear action of a group on $\mathcal{H}$, then so are $e'$ and the homotopy joining them. If $e-f\in\ell^p$, then so is the path between $e-f$ and $e'-f$.
\end{lemma}
\begin{proof} Let $v=(2e-1)(2f-1)-1$. Then $v$ is compact and normal. Let $w'=2+v$. Let finally $p$ be a spectral projection of finite rank to $v$, which is chosen in such a way that the norm of $2 - w(t)$, $w(t) = 2 + ((1-t)+tv(1-p))$, $t\in [0,1]$, is smaller than 2. Then $w(t)$ is invertible. We have $w(0)=3$ and the rank of $w(1)-w'$ is finite. By construction, $(2e-1)w' = w'(2f-1)$, so the rank of $ew(1)-w(1)f$ is finite. Conjugation with the family of invertibles $w(t)$, $t \in [0,1]$, gives a homotopy of idempotents between $e$ and an idempotent $e(1)$, such that $e(1)-f$ has finite rank. Possibly after application of functional calculus, one may assume that $e(t)$, $t \in [0,1]$, are projections.\\
If $e$ and $f$ are invariant under a group action, then so are all operators thus constructed.\\
Let finally $\mathcal{H}' = \Ker (e(1)-f)$. This is a subspace of finite codimension in $\mathcal{H}$. One readily checks that $e(1)\mathcal{H}' = f\mathcal{H}' = \mathcal{H}''$ is a subspace of finite codimension in $e(1)\mathcal{H}$ as well as in $f\mathcal{H}$. For $e'$ one takes orthogonal projection onto $\mathcal{H}''$.\\
The last assertion follows from the fact that only finitely many characteristic values are concerned.
\end{proof}
The next lemma provides the crucial cancellations necessary to make the map $\q ^2\C\to \ell^p$ factor through $\C$.
\begin{lemma}
Let $f*f':A\to B\rhd C$ be a quasi-homomorphism, and let $\Psi:A\to B$ be a homomorphism such that $(\Image f)(\Image\Psi)=0$, i. e. such that both homomorphisms are orthogonal, and similar for $f'$ and $\Psi$.\\
Then $(f+\Psi)(f'+\Psi)$ coincide as homomorphisms $\q A\to C$ with $(f*f')$.\\
This holds in the category of abstract algebras as well as in that of complete bornological ones.
\end{lemma}
\begin{proof}
Since $(f+\Psi)*(f'+\Psi)=(f*f')+(\Psi *\Psi)$ on $\Q A$, this follows from $\Psi *\Psi=0$ on $\q A$. This in turn follows from the exactness of $0\to\q A\to\Q A\to A \to 0$ together with the fact that $\Psi *\Psi=\Psi\circ(\id *\id)$.
\end{proof}
\begin{cor}
Consider a map $\q^2\C\to A$ represented by four idempotents $(e_0*e_1*e_2*e_3)$. Upon composing with $A\to M_2 A$, this map is smoothly homotopic to the one given by the permuted set $(e_0*e_2*e_1*e_3)$, and analogously for the permuted set $e_3*e_1*e_2*e_0$. When all data are equivariant w.r.t. a group action, then so is the homotopy.
\end{cor}
\begin{proof}
The composition  with $A\to M_2 A$ has the effect that we may assume, by rotating $\left(\begin{smallmatrix}e_1&0\\ 0&0\end{smallmatrix}\right)$ to $\left(\begin{smallmatrix}0&0\\ 0&e_1\end{smallmatrix}\right)$, that both $e_0$ and $e_2$ are orthogonal to both $e_1$ and $e_3$, and similarly for $e_3$.\\
Now we apply the following steps to the quasi-homomorphism $(e_0*e_1)*(e_2*e_3)$. We first add $(e_1*0)+(0*e_2)$ to obtain $((e_0+e_1)*(e_1+e_2))*((e_2+e_1)*(e_3+e_2))$. Now we may on the left hand side subtract $(e_1*e_1)$, and on the right hand side $(e_2*e_2)$, using again the lemma.
\end{proof}
Let us now shift interest back to the dual Dirac morphism $\q \CCINF(\R^2)\to\ell^p(L^2\R^2)\otimes M_2$.
\begin{proof}[Proof (of the theorem).]
We want to show that the composition induces an invertible element in $\HA$, and we want the proof to be sufficiently general to apply also when all algebras are tensored with another one, which also has a torus action, and when we take crossed products with $T$.\\
First we observe that the Dirac morphism extends to the unitalization of $\CCINF(\R^2)$ to give a map $\q (\CCINF(\R^2))^+\to\ell^p(L^2\R^2)\otimes M_2$ and that we may assume that the two vector bundles give idempotents on the sphere which are both constant in a neighborhood of the north pole $\infty$ such that they yield a map $\Q\C\to(\CCINF\R^2)^+$. The resulting map $\Q\C\to\mathcal{L}\otimes M_2$ is then given by four idempotents $(e_0*e_1*e_2*e_3)$ where $e_0=p_1Vp_1,e_1=p_1V'p_1,e_2=Fp_2Vp_2F, e_3=Fp_2V'p_2F$. $V$ and $V'$ denote the projections onto the two vector bundles, and $p_1$ and $p_2$ the two projections onto the graded pieces of the Hilbert space.\\
Denote for simplicity $\mathcal{L}'=\mathcal{L}\otimes M_2$.\\
The process of assuring that $(e_0,e_1)\perp (e_2,e_3)$, as required by the foregoing corollary has the effect that we have to deal with two homomorphisms $\Q ^2\C\to\mathcal{L}'\otimes M_2$, represented by a quadruple of idempotents which we shall, by abuse of notation, also denote by $(e_0*e_1*e_2*e_3)$. If we restrict to $\q ^2\C\to\mathcal{L}'\otimes M_2$, the corollary says that this map coincides with the one represented by $(e_0*e_2*e_1*e_3)$. This is the same as the restriction of a map $\Q ^2\C\to \mathcal{L}'\otimes M_2$, given by the same quadruple. The latter homomorphism, in turn, is homotopic to a homomorphism $\Q ^2\C\to\mathcal{L}'\to\mathcal{L}'\otimes M_2$, i. e. a left-upper-corner embedding. So the restriction to $\q ^2 \C$ is a quasi-homomorphism $\q\C\rightrightarrows \mathcal{L}'\rhd\mathcal{L}'$.\\
The two homomorphisms $\q \C\to\mathcal{L}'$ now are in fact homomorphisms $\q\C\to\ell^p\otimes M_2$ because they are quasi-homomorphisms $\C\rightrightarrows\mathcal{L}'\rhd\ell^p\otimes M_2$ by Weyl's theorem.\\
In other words, we have a pair of homomorphisms $\q\C\to\ell^p\otimes M_2$, which is a homomorphism $\Q\q\C\to \ell^p\otimes M_2$. With the help of lemma \ref{yyy} and the cancellation lemma, this homomorphism even factors through $M_\infty\otimes M_2\subset\ell^p\otimes M_2$.\\
In the following, we shall make use of the easily checked fact that two irreducible subrepresentations (so that they are both 2-dimensional) of a finite dimensional real representation of $T$ can be homotopped equivariantly until they coincide, if and only if they have the same weights $(\chi_1,\dots,\chi_n)\in\Z^n$ if $T=(\mathbb{S}^1)^{\times n}$. This is because the representation space splits into a sum over all weights, and in each weight space, all irreducible invariant subspaces have same weight.\\
By repeated applications of the two lemmas, we may thus replace our quadruple by a quadruple $(e*f*0*0)$, where for each weight $\chi\in\Z ^n$, it either does not occur as a subrepresentation of the image of $e$ or $f$, or it occurs in $e$, or in $f$, but not in both. We set $n(\chi)=0$ in the first case, $n(\chi)=$ the complex dimension of the subrepresentation of $e$ in the second, and $n(\chi)=$ minus the complex dimension of the subrepresentation of $f$ in the third.\\
It remains to identify these integers. This can be done by passage to the $C^*$-completion.\\
In fact, observe that $(e_1*\dots*e_3)$ defines an element in $\KK^T(\C,\C)$, and by construction of the Dirac-dual Dirac element, this element is the unit of the representation ring $R(T)=\KK^T_0(\C,\C)$. This follows because $\alpha$ and $\beta$ defined exactly the invertible elements in $\KK^T(\C,C_0(\R^2))$ and $\KK^T(C_0(\R^2),\C)$ given by equivariant Bott periodicity after passage to the $C^*$-completion. The fact that $R(T)$ is generated freely by all such $\chi$ allows to deduce from $\beta\otimes\alpha=1\in\KK^T(\C,\C)$ that $n(0)=1$ and $n(\chi)=0$ if $\chi\ne 0$. On the other hand $[e]-[f]$ defines the same element in $R(T)$. This proves that $e$ is a one-dimensional projector onto an representation of weight 0, and $f=0$.\\
The map $e*0*0*0$ is exactly the description of the map $\q ^2\C\to\ell^p\otimes M_n$ that factors through $\C$.\\
Thus this map does factor through $\C$, as was the claim.\\
It remains to generalize from the representation $\R^2$ to $W$. In the entire sequel, the following argument allows to reduce to $\R^2$. The representation of $T$ on $W=B/T\times V$ splits into a direct sum of irreducibles, and the desired element in $\HA_0((A\protens\CCINF(W))\alcross T,A\alcross T)$ is a composition of elements induced from quasi-homomorphisms, one for each irreducible one. The construction then also allows to take crossed products with $B$, as long as the action factors through $T$.\\
\end{proof}
\section{The proof that $\alpha$ and $\beta$ are inverse isomorphisms on HA of the crossed product}
We have shown that  $\alpha\circ\beta:\HA(A\alcross B)\to\HA((A\protens\CCINF(W))\alcross B)\to\HA(A\alcross B)$ is the identity whenever the representation of $B$ on $W$ factors through $T$. In this section, we set the goal to get rid of that condition.\\
As a consequence of the theorem, we may tensor the commuting triangle diagram (bornologically) with $\CCINF(G)$, take $B$ crossed products, and we may take $\Gamma$-crossed products to obtain a diagram which commutes up to smooth homotopy.\\
For the following corollary, remember that $B$ always acts on $B/T\times V$ via $B\to\alTor\to T$.
\begin{cor}
The elements
$$
\beta\in\HA_0(\Gamma\allcross\CCINF(G)\alcross B,\Gamma\allcross (\CCINF(G\times B/T\times V)\alcross B)
$$
and
\begin{align*}
\alpha\otimes&\Morita,\ \mathrm{where}\\
&\alpha\in\HA_0(\Gamma\allcross (\CCINF(G\times B/T\times V)\alcross B,\Gamma\allcross (\CCINF(G)\protens\ell^p)\alcross B)\ \mathrm{and}\\
&\Morita\in\HA_0(\Gamma\allcross (\CCINF(G)\protens\ell^p)\alcross B,\Gamma\allcross \CCINF(G)\alcross B)
\end{align*}
satisfy $\beta\otimes\alpha\otimes\Morita=1$. Furthermore, $\beta$ is injective, $\alpha$ is surjective in $\HA_*$.\\

\end{cor}
\begin{proof}
We apply the theorem to $W=B/T\times V$. So firstly, the composition
$$
\HA(\C,\CCINF(B/T\times V)\otimes \HA_0(\CCINF(B/T\times V,\ell^p)
$$
maps $\beta\otimes\alpha$ exactly to the element in $\HA(\C,\ell^p)$ coming from $\C\to\ell^p$ i. e. to $\Morita$.\\
Now we may take the induction from $B$ to $G$ and the $\Gamma$-descent of the diagram. More precisely, we have described how to obtain from a quasi-homomorphism a new one taking bornologically complete tensor products or descents. Thus we have described the operation of sending a quasi-homomorphism $x$ as in the diagram to $\Gamma\allcross (\CCINF(G)\protens x)\alcross B$.\\
Now the homotopy, if given as $\q^2\C\mapsto\ell^p\protens\CINF[0,1]$ can be taken the same process because it is given by a path of \emph{$T$-invariant homomorphisms}.\\
Finally, the algebra $\CINF[0,1]$ can be pulled out of the $B$ and $\Gamma$-crossed products because the action decomposes as a tensor product action such that the one on $\CINF[0,1$ is trivial. So $\CINF[0,1]$ can be pulled out on vector space level by the associativity of the completed bornological tensor product and this is an algebra map as well, thus yielding the requisite homotopy.\\
\end{proof}
\subsection{A two-step retraction from $B\to T$ through group homomorphisms}
The goal of this section is to introduce the group retraction that will allow the calculation of $\HA$ of the ``right'' crossed product where the action of $B$ does not factor.\\
We now have to restrict from the case of general semisimple Lie groups to the case $G=\SL(n,\C)$. So $B$ is the subgroup of upper triangular matrices.\\
Denote by $\alTor$ the algebraic non-compact torus of diagonal matrices. In the Iwasawa decomposition $\alTor=TA$. $A$ is the group of diagonal matrices with real positive entries whose product is $1$. $A$ is isomorphic to the abelian group $\R^{n-1}$. Since $\alTor$ is an abelian group, this decomposition is a direct product of abelian groups. So there are canonical projections $\alTor\to T$ and $\alTor\to A$ (and, of course, injections $T\to\alTor$ and $A\to\alTor$) which are group homomorphisms.\\
The following formula defines a $[0,1]$-family of group homomorphisms $B\to B$: For $t\in[0,1]$, set
$$
r_t:\left( \begin{array}{cccc}
z_{1} & * & * &\ldots \\
0 & z_2 & * &\ldots \\
0 & 0 & z_3 &\ldots \\
\vdots & \vdots & \vdots & \ddots
\end{array} \right)
\mapsto
\left( \begin{array}{cccc}
z_{1} & t* & t^2* &\ldots \\
0 & z_2 & t* &\ldots \\
0 & 0 & z_3 &\ldots \\
\vdots & \vdots & \vdots & \ddots
\end{array} \right),
$$
i. e. $r_t$ multiplies the $k$-th off-diagonal with $t^k$.\\
This gives a smooth path between $r_1=\id_B$ and the projection $B\to\alTor$ onto diagonal matrices through group homomorphisms.\\
Furthermore there is a smooth path $(s_t)_{s\in[0,1]}$ of homomorphisms $A=\R^{n-1}\to\R^{n-1}=A$ given by $s_t(x)=(1-t)x$, which as well is a path through group homomorphisms between $s_1=\id_{\R^{n-1}}$ and the zero homomorphism. Taking the direct product of this path with the compact torus $T$, we see that there is a path, which we shall also denote by $(s_t)$, between the identity $s_1$ on $\alTor$ and the projection $\alTor \to T$.\\
So if we define a $B$-action on source and range of $\beta^T$ by composing the $T$-actions with $s_0\circ r_0:B\to\alTor\to T$, the homomorphism $\beta^T$ is $B$-equivariant. The $B$-action however is not the right one. In fact, the action $B\to\alTor\to T\to \AutTop (B/T)$ defined via these two group retraction is the action of $T$ on $B/T$, i. e. under the exponential mapping exactly the adjoint representation $T\to \SO(\mathfrak{b}/\mathfrak{t})$.\\
The space $Y=G\times_B B/T\times V$ where $B$ acts on $B/T$ via $T$, where the action on $G$ is the standard one, is a homogeneous vector bundle over $G/B$, endowed with a $\Gamma$-action. However, the algebras $\Gamma\allcross\CCINF(Y)$ and $\Gamma\allcross\CCINF(G/B\times V)$ are \emph{not} homotopy equivalent in any obvious way.

\section{Computation of the Dirac map on the homogeneous part of cyclic cohomology}
In this section, we want to compute the effect of $\alpha$ on periodic cyclic homology. The results we will develop for $\HA$ do not immediately apply to that case too, because the homotopy invariance is too weak.
\begin{thm}\label{zzz}
The quasi-homomorphism $\alpha$ induces an isomorphism
\begin{align*}
\HP_*(\Gamma\allcross\CCINF(G/T\times V))_{<1>}&\cong \HP_*(\Gamma\allcross\CCINF(G\times B/T\times V))_{<1>}\\
\longrightarrow\HP_*(\Gamma\allcross\CCINF(G)\alcross B)_{<1>}&\cong\HP_*(\Gamma\allcross\CCINF(G/B))_{<1>}
\end{align*}
\end{thm}
\begin{proof}
It suffices to prove that $\alpha$ induces an isomorphism on the $E_2$-term of the spectral sequence. This map is
$$
\HP_*(\CCINF(G\times B/T\times V)\alcross B)\to\HP_*(\CCINF(G)\alcross B).
$$
This is a map from $\Homolc^{\ev/\odd}(G/T\times V)$ to $\Homolc^{\ev/\odd}(G/B)$. The fact that this map is the Chern character of the K-theory map induced by
$$
j^B(C_0(G)\otimes\alpha)\in\KK(C_0(G\times B/T\times V)\rtimes B,C_0(G)\rtimes B),
$$
where $\alpha\in\KK^B(C_0(B/T\times V),\C)$, is quite obvious, and proved formally by the commutativity of the following diagram
$$
\xymatrix{
\HP_*(\CCINF(G/T\times V))\ar[r]^-\alpha\ar[d]^\cong&\HP_*(\CINF(G/B))\ar[d]^\cong\\
\HA_*(\CCINF(G/T\times V))\ar[r]^-\alpha\ar[d]^\cong&\HA_*(\CINF(G/B))\ar[d]^\cong\\
\HA_*(C_0(G/T\times V))\ar[r]^-{\overline{\alpha}}&\HA_*(C(G/B))\\
\K_*(C_0(G/T\times V))\ar[u]_\ch\ar[r]^-{\overline{\alpha}}&\K_*(C(G/B))\ar[u]_\ch
}
$$
where $\overline{\alpha}$ denotes the $C^*$-completion of $\alpha$, and the first three lines are all canonically and in a consistent way isomorphic to a map $\Homolc^{\ev/\odd}(G/T\times V)\to\Homol^{\ev/\odd}(G/B)$. For the fourth line, this applies when it is tensored with $\C$, and then the map that arises is the Thom isomorphism. Or simply, the map is an isomorphism because $\overline{\alpha}$ is invertible.
\end{proof}

\section{Computation of analytic cyclic homology of the crossed products}
We may now state the most important technical theorem of this paper, namely the computation of the analytic cyclic (co-)homology of the crossed product $\HA_*(\Gamma\allcross\CINF(G/B))$.
\begin{thm}\label{thmwithdiagram}
The element $\alpha\in\HA_0(\Gamma\allcross\CINF(G/B),\Gamma\allcross\CCINF(G/T))$ has a left inverse as an $\HA$-element in the sense that that there is $\beta$ with $\beta\otimes\alpha=1$. So if we think of the induced maps in homology, $\alpha$ is surjective. %

\end{thm}
As a consequence, we obtain a surjective map on the cohomology of $\Gamma\backslash G/T$ by Morita equivalence.\\
The proof of the theorem will rely on the following preliminary definitions and lemma.\\
Consider the algebra $\CINF [0,1]\protens\ell^p$. It is a subalgebra of the $C^*$-completion $C[0,1]\otimes_{C^*}\Ko$. The latter algebra is the algebra of continuous maps from $[0,1$ to the compacts i. e. of the sections of the trivial bundle with fiber $\Ko$ over the interval. If we view the fiber over $t$ as $\Ko L^2 (B/T)$ i. e. as the completion of $\ell^p L^2 (B/T)_t$ where $(B/T)$ is the space $B/T$ endowed with the action $r_t$ of $B$ we see that this bundle is equipped with a $B$-action. So $B$ acts over $t$ via $r_t$. This $B$-action makes the algebra of sections a $B$-algebra and preserves the subalgebra $\CINF [0,1]\protens\ell^p$. We have thus defined the crossed product $\CINF [0,1]\protens \ell^p\alcross B$
\begin{lemma}[A partial homotopy result for the actions]
Consider the algebra
$$
\Gamma\allcross \CCINF(B/T\times V\times [0,1])\times B,
$$
where $B$ acts as explained above on the fiber over $t\in[0,1]$ via $r_t$.\\
Consider the evaluation at $0$ from this algebra to $\Gamma\allcross \CCINF(B/T\times V)\times B$ where $B$ acts via the algebraic torus $\alTor$. This map is invertible in $\HA$ and $\HP$.\\
More precisely, consider the extension
$$
0\to \Ker(\ev_0)\to \bigl(\CCINF(G\times [0,1])\protens\ell^p\bigr)\alcross B\to \bigl(\CCINF(G)\protens\ell^pL^2(B/T)_0\bigr)\alcross B\to 0.
$$
Then the ideal is $\HA$-contractible. This also holds for the crossed product with $\Gamma$.\\
A similar statement holds if we let $B$ act on the bundle via $s_t\circ(B\to\alTor)$ so that the evaluation goes to the algebra with $T$-action.\\

\end{lemma}
\begin{rmk}[Warning]
\emph{A priori}, there is absolutely no reason to assume that the evaluations in the endpoints $1$ of the intervals should be invertible. In fact, \emph{any} $\KK$-element from $\KK(A, C)$ can be represented with the tools of asymptotic morphisms as a correspondence $\K(A)\cong\K(B)\to\K(C)$. This is shorthand for morphisms induced by bivariant elements. The algebra $B$ is an algebra of sections, the first isomorphism is of similar nature as the ones presented in the lemma, namely evaluations at one of the endpoints of the interval, and the second one is evaluation at the other endpoint.\\
However, our proof of the theorem, involving the whole machinery of the Dirac operator and equivariant Bott periodicity, will show that in our particular case, the evaluations are invertible, though.
\end{rmk}
\begin{proof} (of the lemma).\\
The idea is that for every $t\ne 0$ the action $r_t$ is faithful. This implies that there is an isomorphism between the ideal and an algebra of sections of a constant bundle over $(0,1]$. More precisely, define the isomorphism
$$
\Ker(\ev_0)\to (\CCINF(G\times (0,1])\protens \ell^p L^2(B/T)_1)\allcross B
$$
by $f\mapsto (t\mapsto f(r_t^{-1}))$ is $f$ is a smooth function on $B$. Here for functions with compact support we have to understand functions with all derivatives vanishing at zero. The action of $B$ on $(0,1]$ is trivial and constant on $\CCINF(G)\protens \ell^p$, thus $\CCINF(0,1]$ can be pulled out of the crossed product yielding an isomorphism
$$
(\CCINF(G\times (0,1])\protens \ell^p L^2(B/T)_1)\allcross B\cong \bigl((\CCINF(G)\protens \ell^p L^2(B/T)_1)\allcross B\bigr)\protens \CCINF(0,1].
$$
The latter algebra is $HA$-contractible by an explicit formula for a contraction, see \cite{Loday}*{p. 120}.\\
Taking crossed products with $\Gamma$ does not alter the argument at all because all maps and homotopies are by construction $\Gamma$-equivariant.
\end{proof}
\begin{proof} (of the theorem)
Consider the following commuting diagram in the category of bornological algebras endowed with right $B$-action and bivariant $\HA$-elements.\\
The indices are supposed to denote the $B$-actions. Actions on algebras that may be viewed as algebras of smooth sections of some bundle over $[0,1]$ such that the action depends on the fiber $t$ are denoted with brackets in the index i. e. $[T,\alTor]$ means that the action on the fiber over $t$ is supposed to be the action of $B$ via the group homomorphism $r_t$. A similar remark holds for $[t,\alTor]$ endowed with the action $s_t\circ(B\to\alTor)$.
$$
{\normalsize
\xymatrix@C=6em{
   (\CCINF(B/T\times V))_T \ar[r]^{\alpha}  & (\ell^p)_T \ar[r]^{\Morita}       &       \C
\\
(\CCINF(B/T\times V\times[0,1]))_{[T,\alTor]}\ar[r]^-{\alpha\protens\CINF[0,1]}\ar[u]^{\ev_0}\ar[d]_{\ev_1}& (\ell^p\protens\CINF[0,1])_{[T,\alTor]}\ar[r]^{\Morita}\ar[u]^{\ev_0}\ar[d]_{\ev_1}&\CINF[0,1]\ar[u]^{\ev_0}\ar[d]_{\ev_1}
\\ 
(\CCINF(B/T\times V))_{\alTor}\ar[r]^{\alpha}&(\ell^p)_{\alTor}\ar[r]^{\Morita}&\C
\\
(\CCINF(B/T\times V\times[0,1]))_{[\alTor,B]}\ar[r]^-{\alpha\protens\CINF[0,1]}\ar[u]^{\ev_0}\ar[d]_{\ev_1}&(\ell^p\protens\CINF[0,1])_{[\alTor,B]}\ar[r]^{\Morita}\ar[u]^{\ev_0}\ar[d]_{\ev_1}&\CINF[0,1]\ar[u]^{\ev_0}\ar[d]_{\ev_1}\\
(\CCINF(B/T\times V))_B\ar[r]^{\alpha}&(\ell^p)_B\ar[r]^{\Morita}&\C
}
}
$$
The algebras on the right hand side column are endowed with the trivial action.\\
Forgetting the actions, it is clear that all maps are isomorphisms in $\HA$. Since for all $t$ $\alpha_t$ are equal, an inverse is given by $\beta_0$. However we want to show that the $\alpha$ on the last line is invertible in such a way that the proof also applies when the induction and descent process is applied.\\
In other words, we now take completed bornological tensor products with $\CCINF(G)$, take crossed products with $B$ from the right which is possible since all maps are $B$-equivariant and finally crossed products with $\Gamma$ from the left, which is also possible.\\
First we observe that the arrows in the right hand column are invertible by homotopy invariance proper. Second, the Morita arrows are invertible, so all four arrows in the middle column are invertible as well. Now we have shown in the preceding lemma that both arrows $\ev_0$ are invertible. Furthermore it was the content of the preceding section that $\alpha$ is surjective, in the sense of having a one-sided $\HA$-inverse, if we take the $T$-action. So the $\alpha$ in the first line is surjective. So the second, third, fourth-line $\alpha$'s are surjective, which makes the fifth-line $\alpha$ surjective as well.
\end{proof}
Concerning $\HP_{<1>}$, we have carried out the necessary computation in the preceding section. The preceding proof does not apply verbatim to $\HP_{<1>}$, because there, the homotopy invariance theorem is weaker in the sense that it only applies to the algebraic tensor product. Nevertheless, by using the spectral sequence, we didn't need to take crossed products with $\Gamma$. We knew that all arrows induce isomorphisms on the $E^2$-term and thus got the desired result, at least right right for $\alpha$, where we need it. (However, we expect the theorem to hold in $\HP$ as well, which would have the consequence of a vanishing inhomogeneous part).\\
Next we shall state the main theorem in the form in which we shall actually need it, i. e. we will show in which sense our computations of the crossed products in the various theories are compatible.
\begin{cor}\mbox{ }\\
\vspace{-1.5em}
\begin{enumerate}
\item
The canonical natural transformation
$$
\HA_*(\Gamma\allcross \CINF(G/B))\to\HP_*(\Gamma\allcross \CINF(G/B))_{<1>}
$$
induced from the inclusion of the total complex with growth condition into the product total complex, followed by projection onto the homogeneous part, is an isomorphism.
\item
The map
$$
\HA_*(\Gamma\allcross \CINF(G/B))\to\HA_*(\Gamma\rlcross C(G/B))
$$
is an isomorphism.
\end{enumerate}
\begin{proof}\mbox{ }\\
\vspace{-1.5em}
\begin{enumerate}
\item
Consider the commutative diagram
$$
\xymatrix{
\HP_*(\Gamma\allcross \CCINF(G/T\times V))_{<1>}\ar[r]^\alpha&\HP_*(\Gamma\allcross \CCINF(G/B))_{<1>}\\
\HA_*(\Gamma\allcross \CCINF(G/T\times V))\ar[r]^\alpha\ar[u]&\HA_*(\Gamma\allcross \CCINF(G/B)).\ar[u]
}
$$
By the Morita equivalence, the left vertical arrow is an isomorphism.
\item
For this we need to invoke quite an advanced tool, namely local cyclic homology as developed by Puschnigg (see \cite{Pulocal} and introduction). It is a fact proved in \cite{Pulocal} that for a large class of algebras, that $C^*$-algebras belong to, local cyclic homology coincides with analytic cyclic homology. However, in cohomology, there is no analogous statement available. Moreover, there is in full generality a \emph{bivariant Chern-Connes character} from $\KK$-theory to bivariant local cyclic homology with the desired property of being a multiplicative transformation. This may be called a generalized Grothendieck-Riemann-Roch theorem.. It follows that $\HA_*$ of $C^*$-algebras may be calculated by the Dirac-dual Dirac method. More precisely, the following diagram commutes, and once again, one can read off from the fact that some arrows are isomorphisms that the others are as well.
$$
\xymatrix{
\HA_*(\Gamma\allcross \CCINF(G/B))\ar[r]^\alpha\ar[d]^{i_*}&\HA_*(\Gamma\allcross \CCINF(G/T))\ar[d]^{i_*}\\
\HA_*(\Gamma\rlcross C(G/B))\ar[r]^\alpha&\HA_*(\Gamma\rlcross C_0(G/T))\\
\HPL_*(\Gamma\rlcross C(G/B))\ar[r]^\alpha\ar[u]_\cong&\HPL_*(\Gamma\rlcross C_0(G/T))\ar[u]_\cong\\
}
$$
\end{enumerate}

The proof shows (by the uniqueness of an inverse, but in the weak sense on the level of maps) that an inverse can be constructed as follows. It is the composition of the dual-Dirac element from $\C$ to $(\CCINF(B/T\times V))_T$ followed by the composition of the left hand side arrows of the diagram.
\end{proof}
\end{cor}
\section{The idempotent theorem}
The goal of this section is to show how to use the previous results to prove the statement of the idempotent theorem. The large diagram might convince the reader that it was necessary to compute $\HA_*(\Gamma\allcross\CINF(G/B))$. Moreover, it shows how to reduce our problem to the Connes-Moscocivi index theorem, avoiding Connes' foliation index theorem.\\
Let $X$ be a smooth $\Gamma$-manifold.\\
We constructed in section \ref{DdDG} a Dirac element $\alpha\in\KK^\Gamma(C_0(\EG),\allowbreak\C)$ and in the proof of the theorem \ref{mainKtheorem} of the K-theory part an ``assembly map'' between $\K$-groups of $C^*$-algebras
$$
\K^{*+q}(\EG\times_\Gamma X)\to\K_*(\Gamma\rlcross C_0(X))
$$
as right Kasparov product with the element $\Morita\otimes j(\tau_{C_0(X)}(\alpha^\Gamma))$ where $\Morita$ denotes the obvious Morita invertible element, $j$ the descent homomorphism and $\tau_A$ the operation of tensoring Kasparov-$\KK^\Gamma$-cycles with a fixed $\Gamma-C^*$-algebra $A$. The construction relied on existence of a $G$-equivariant Spin$^c$-structure on $\EG=G/K$. However, this is automatically satisfied for $G=\SL(n,\C)$, see remark \ref{SPINN2} and the proof of theorem \ref{SPINN1}. More generally, complex semisimple Lie groups are always simply connected, so that the existence is automatic for them as well.\\
The following lemma is in content very similar to theorem \ref{SPINN1}, yet appearing here from the viewpoint of the $\Gamma$-Dirac element, instead of that of the $B$-Dirac element.
\begin{lemma}
If $X=G/B$, the assembly map
$$
\K^{*+q}(\EG\times_\Gamma G/B)\to\K_*(\Gamma\rlcross C(G/B))
$$
is an isomorphism.
\end{lemma}
\begin{proof}
Let us show that, if we compose
$$
j^\Gamma(\alpha^\Gamma\otimes C(G/B))\in\KK_q(C_0(\EG\times_\Gamma G/B))
$$
with the invertible
$$
j^(\Gamma\times B)(C_0(G)\otimes \beta^B)\in KK_q(\Gamma\rlcross C(G/B),\Gamma\rlcross C_0(G/T)),
$$
then we obtain an invertible element. Again, we suppress in the notation $\KK$-equivalences coming from Morita equivalences. In fact, by a calculation of the form we have come to know already, the composition is equal to the product
$$
j^{\Gamma\times B}(C_0(\EG\times G)\otimes \beta^B)\otimes j^{\Gamma\times B}(\alpha^\Gamma\otimes C_0(G\times B/T)).
$$
The first factor is invertible because $\beta^B$ is, and the second is the invertible element in $\KK_0(C_0(\EG\times_\Gamma G/T),C_0(\Gamma\backslash G/T))$ coming from the Thom isomorphism for the vector bundle $\EG\times_\Gamma G/T\to \Gamma\backslash G/T$. This is a consequence of the construction of $\alpha$. Precisely, the statement follows from the principle that the assembly map for $\Gamma$ with coefficients in a \emph{proper} $\Gamma$-algebra is an isomorphism (see \cite{Valette}*{theorem 9.4}).\\
In other words, consider the following commuting diagram
$$
\xymatrix@C=145pt{
 {\K^q(\EG\times_\Gamma G/B)}
\ar[r]^{-\otimes\Morita\otimes j(\tau_{C(G/B)}(\alpha^\Gamma)}
\ar[d]^{-\otimes\Morita\otimes j(\tau_{C_0(\EG)}(\beta)}
&
{\K_0(\Gamma\rlcross G/B)}
\ar[d]^{-\otimes j(\beta)}\\
\K^q(\EG\times_\Gamma G/T\times V)\ar[r]^{-\otimes\Morita\otimes j(\tau_{C_0(G/T\times V)}(\alpha^\Gamma)} & {\K_0(\Gamma\ltimes C_0(G/T\times V))}
}
$$
In the upper row, $X=G/B$, in the lower one $X=G/T\times V$.
The diagram commutes by the associativity of the exterior Kasparov product. The lower arrow is an isomorphism because modulo Morita equivalence it is the Thom isomorphism for the bundle $\EG\times_\Gamma G/T\times V\to \Gamma\backslash G/T\times V$. In fact, the ordinary Thom isomorphism holds between the cohomology without support condition of the base space and the cohomology with compact supports of the vector bundle. \\
We already know that the vertical arrows are isomorphisms because $\beta$ is invertible. So the upper arrow is an isomorphism.
\end{proof}
Let now $e$ be an idempotent in $\C\Gamma$. We have
$<e,\tau>=<e,i^*\varphi>=<i_*e,\varphi>=<i_*e,\beta^*\widetilde{\varphi}>=<\beta(i_*e),\widetilde{\varphi}>=<\ch \beta(i_*e),\widetilde{\varphi}>.$
\begin{thm}
The image $\ch(e)\in\HP_0(\C\Gamma)_{<1>}$ of $e$ is also an image of a $[D]\in \K_0(\BG)$ under $\mu\ch$, if $\Gamma$ acts cocompactly.
\end{thm}
\begin{proof}
Consider the following diagram\label{Schlussdiagramm}
$$
\xymatrix@C=0em@R=1.8em{
{\K_0(\BG) \ar[rd]^{(1)} \ar[rr]^{(2)}\ar[d]^{(3)}}  &         &   {\mathbb{H}_\ev(\Gamma)\ar[rd]^{(4)}}&
\\
{\K^q(\BG)}\ar[d]^{(5)}              &{\K_0(\mathcal{R}\otimes\C\Gamma)} \ar[dd] ^{(6)}\ar[rr] _>>>>>>>>>>>>>>{(7)}   &    &  {\HP_0(\C\Gamma)_{<1>}\ar[dd]^{(8)}}
\\
{\K^q(\EG\times_\Gamma G/B)}\ar[dddr]_{(9)}  &     & &  
\\
 & \K_0\bigl(\mathcal{R}\otimes\Gamma\allcross\CINF(G/B)\bigr)\ar[rr]^>>>>>>>{(10)}\ar[dd]^>>>>>>>>>>>{(11)}\ar[rrd]_{(15)} & &{\HP_0\bigl(\Gamma\allcross \CINF(G/B)\bigr)_{<1>}}\\
&&&\HA_0\bigl(\Gamma\allcross \CINF(G/B)\bigr)\ar[u]_{(12)}\ar[d]^{(13)}\\
&\K_0(\Gamma\rlcross C(G/B))\ar[rr]_{(14)}&&\HA_0(\Gamma\rlcross C(G/B))
}
$$
where the arrows are defined as follows:
\begin{enumerate}[(1)]
\item is the ordinary assembly map. It factors over $\K_0(\mathcal{R}\Gamma)$ by \cite{ConnesMosc} and \cite{Thom}, where $\mathcal{R}\Gamma$ is the ring of infinite matrices with rapid decay.
\item is the homological Chern character, as defined by \cite{Baumetc} (see also \cite{Valette}, lemma A.4.2).
\item is $\K$-theoretic Poincare duality. Observe that
$$
\dim \BG=\dim G/K= \dim B/T=q
$$
by the Langlands decomposition.
\item is the special case $A=\C$ of the construction carried out in the previous sections
\item is the map $\pi^*$.
\item is the map induced by the shrinking map, studied abundantly in the preceding sections.
\item is the Chern character as defined by Connes, followed by the Morita equivalence $\mathcal{R}\sim\C$ and projection onto homogeneous part. 
\item is equally induced by $G/B\to \pt$.
\item is the assembly map defined in the K-theory part and reviewed in the preceding lemma.
\item is again Connes' Chern character, followed by similar maps as in (7).
\item is the map induced by the inclusion, tensored with the inclusion $\mathcal{R}\to\Ko$ and the completion to the $C^*$-tensor product.
\item is the natural transformation from analytic to periodic cyclic homology. We have proved in the preceding section that it is an isomorphism.
\item is the map induced by the inclusion, followed by the same maps as in (11)
\item is the Chern character as defined by Meyer (\cite{Meyerthesis}), followed by Morita equivalence, given by the isomorphism with local cyclic homology.
\item is again Meyer's Chern character.
\end{enumerate}
Let us first prove that the diagram commutes. The commutativity of the upper horizontal diagram is a corollary of the Connes-Moscovici index theorem.\\
The commutativity of both the upper and lower front side rectangle diagrams is the naturality of the Chern character.\\
The commutativity of the front side triangular diagram follows from a look at Meyer's Chern character.
Only the commutativity of the left hand side requires a closer look. First, there is, following \cite{Kas}, the description of the assembly map $\mu:\K_o(\BG)\to\K_0(C^*_\mathrm{r}\Gamma)$ as Poincare duality $\K_0(\BG)\to\K^q(\BG)$ followed by right multiplication with $j(\alpha^\Gamma)$. The map $\K_0(C^*_\mathrm{r}\Gamma)\to K_0(\Gamma\ltimes_\mathrm{r}C(G/B))$ is right multiplication with $j(G/B\to\pt)$, and there is the equality $\tau_{C_0(\EG)}(G/B\to\pt)\otimes\tau_{C(G/B)}(\alpha^\Gamma)=\alpha\otimes(G/B\to\pt)$. Consider the following diagram.
$$
{\normalsize
\xymatrix@C=7em{
\K_0(\BG)\ar[rd]_\mu\ar[rrd]^{(1)=\mu_\mathrm{CM}}\ar[d]_{(3)}^\PD&&\\
\K^q(\BG)\ar[d]_{(16)}\ar[r]_{\otimes j(\alpha^\Gamma)} &\K_0(C^*_\mathrm{r}\Gamma)\ar[d]&\K_0(\mathcal{R}\Gamma)\ar[l]\ar[d]^{i_*}\\
\K^q(\EG\times_\Gamma G/B)\ar[r]_{\otimes j(\tau_{C(G/B)}(\alpha^\Gamma))}& \K_0(\Gamma\rlcross C(G/B))&\K_0(\Gamma\allcross\CINF(G/B))\ar[l]
}
}
$$
The uppermost triangle commutes because these are just two ways of writing the assembly map. The left hand triangle commutes by Kasparov's description of the assembly map. The commutativity of the left rectangle is a property of the intersection product. The commutativity of the right hand rectangle is obvious because it commutes already on the algebra level.\\
So after all we have established the commutativity of both diagrams.\\
So we have defined a "commuting cube". Let us now have a closer look at the ``backside'' diagram.
$$
\xymatrix@C=6em{
{\K_0(\BG) \ar[r]^{(2)}\ar[d]^{(5)\circ(3)}}  &  {\Homol_\ev(\Gamma)} \ar[d]^{(13)\circ(12)^{-1}\circ(8)\circ(4)}  \\
{\K^q(\EG\times_\Gamma G/B)}\ar[r]^-{(14)\circ(9)}& {\HA_0(\Gamma\rlcross C(G/B))}
}
$$
The calculations of the preceding chapter imply that the lower arrow $(9)$ is rationally injective. Consider the subgroup $((14)\circ(9))^{-1}(\mathrm{Im} \Homol(\Gamma))$ of $\K^q(\EG\times_\Gamma G/B)$. By the commutativity this subgroup contains the image of $\K_0(\BG)$ and by the split injectivity of $(5)\circ(3)$ it splits as $\mathrm{Im} \K_0(\BG)\oplus X$ for some subgroup $X$. Assume $X$ had an element of infinite order. Then its image under $(14)\circ(9)$ would be non-zero by the injectivity, and it would come from an element in $\Homol(\Gamma)$. By the rational injectivity of the topological Chern character (2), this element would come from an element in $\K_0(\BG)\otimes\C$. But in view of $(\K_0(\BG)\oplus X)\otimes \C=\K_0(\BG)\otimes \C\oplus X\otimes \C$ this is a contradiction to the fact that $X\otimes\C$ would not vanish. So $X$ is only made of torsion, and we see that the element $i_*e\in\K_0(\Gamma\rlcross C(G/B))\cong\K^q(\EG\times_\Gamma G/B)$, which maps to $\ch i_* e$ in the bottom diagram, is up to torsion in the image of $\K_0(\BG)$.\\
This shows that we have found an idempotent that maps to the same element in $\HA_0$ of the $C^*$-algebra. Since the map induced by inclusion
$$
\HA_0(\Gamma\allcross\CINF(G/B))\to \HP_0(\Gamma\rlcross C(G/B))
$$
and the natural transformation
$$
\HP_0(\Gamma\allcross\CINF(G/B))\to \HA_0(\Gamma\allcross\CINF(G/B))
$$
are isomorphisms and the arrow (8) is injective, the two classes must map even to the same element in $\HP_0(\C\Gamma)_{<1>}$. 
\end{proof}
\begin{cor}
Let $\Gamma$ be a discrete cocompact torsion-free lattice in $\SL(n,\C)$, $n\ge 2$. Then $\C\Gamma$ contains no nontrivial idempotents.
\end{cor}
\begin{proof}
We have $<e,\tau>=<\mu(D),\tau>=\Ind D\in\Z$  with the help of Atiyah`s $L^2$ index theorem. This implies the idempotent theorem as explained in the introduction.
\end{proof}

\clearpage{\pagestyle{empty}\cleardoublepage}
\begin{appendix}
\chapter[Appendices]{}
\section{Globally symmetric Riemannian spaces}\label{symmetric space}
Let us fix an arbitrary semisimple connected Lie group $G$, a maximal compact subgroup $K$ and a discrete torsion-free subgroup $\Gamma$ of $G$. We consider the right action of $K$ on $G$ and the left action of $\Gamma$ on $G$ and the homogeneous space $G/K$. The space $G/K$ is the standard example of a Riemannian globally symmetric space, see \cite{Helgason}*{IV.3.4}. This means that $G/K$ is a Riemannian manifold (such that the underlying manifold structure is exactly the one of the homogeneous space, of course), and every point of $G/K$ is an isolated fixed point (in fact, \emph{the} fixed point) of an involutive isometry of $G/K$. Let us summarize how to show this. Consider the Cartan decomposition
$$
\mathfrak{g} = \mathfrak{k}\oplus \mathfrak{m}
$$
of the Lie algebra $\mathfrak{g}$ of $G$ into the $+1$ resp.$-1$-eigenspaces of the Cartan involution $\sigma$ of $\mathfrak{g}$ (for the existence of the involution and decomposition, see \cite{HilgertNeeb}*{III.6}), i. e. $\sigma:\mathfrak{g}\to\mathfrak{g}$, $\sigma^2=1$, $\sigma = 1_\mathfrak{k}\oplus -1_\mathfrak{m}$ and the Killing form $B:\mathfrak{g}\times\mathfrak{g}\to\R,\ B(X,Y)=\Tr(\ad X\ad Y)$ (which is non-degenerate by definition of semisimplicity) is negative definite on $\mathfrak{k}\times\mathfrak{k}$, positive definite on $\mathfrak{m}\times\mathfrak{m}$. The involution $\sigma$ can be chosen in such a way that the $+1$-eigenspace is exactly the tangent space to the chosen maximal compact subgroup $K$ of $G$ at the unit element, i. e. its Lie algebra $\mathfrak{k}$.\\
As an immediate consequence, $\mathfrak{m}$ and $\mathfrak{p}$ are orthogonal with respect to the Killing form since for $X\in\mathfrak{k}$ and $Y\in\mathfrak{m}$ we have
$$
B(X,Y)=-B(\sigma X,\sigma Y)=-\Tr(\ad(\sigma X)\ad(\sigma Y))=-\Tr(\sigma\ad X\ad Y\sigma^{-1})=-B(X,Y)
$$
which therefore is zero.\\
Thus, $[\mathfrak{k},\mathfrak{k}]\subset\mathfrak{k}$ (this is also clear because $K$ is a subgroup), $[\mathfrak{k},\mathfrak{m}]\subset\mathfrak{m}$ and $[\mathfrak{m},\mathfrak{m}]\subset\mathfrak{k}]$. Let $Ad:G\to\mathcal{L}(\mathfrak{g})$ and $\ad:\mathfrak{g}\to\mathcal{L}(\mathfrak{g})$ denote the adjoint representations. Now the fact that for any subspace $\mathfrak{a}$ of $\mathfrak{g}$ the Lie algebra of $\{g\in G:\Ad (g)\mathfrak{a}\subset\mathfrak{a}\}$ is the same as $\{X\in\mathfrak{g}:\ad (X)\mathfrak{a}\subset\mathfrak{a}\}$ implies that the adjoint representation $G\to\mathcal{L}(\mathfrak{g})$ restricts to a representation of $K$ on $\mathfrak{m}$.\\
Now the map
$$
K\times\mathfrak{m}\to G,\ (k,M)\mapsto k\exp(M)
$$
is a diffeomorphism (\cite{HilgertNeeb}*{III.6.7}). In particular, $G/K$ is $K$-equivariantly diffeomorphic to a Euclidean space $\mathfrak{m}=\R^i$, and the left action of $K$ on $G/K$ is diffeomorphic to the linear action of $K$ on $\mathfrak{m}$. This fact will become essential when we set out to construct equivariant Bott elements in $\KK^K(\C,C_0(G/K))$.\\
The summand $\mathfrak{m}$ of $\mathfrak{g}$ is a Euclidean vector space with the Killing form, and the adjoint representation of $K$ acts through $\SO(\mathfrak{m})$. In fact, let $Z\in\mathfrak{k}$, $X,Y\in\mathfrak{m}$. Then it follows from
$$
B(\ad Z(X),Y)=\Tr([[Z,X],[Y,-]])=-\Tr([X,[Z,Y]])=-B(X,\ad Z(Y))
$$
that $\ad \mathfrak{k}\subset \mathfrak{so}(m)$ and thus $K\to\SO(\mathfrak{m})$. The Euclidean structure makes $\mathfrak{m}$ a Riemannian manifold, and thus also $G/K$ by the exponential map. This Riemannian structure is then $G$-left-invariant by construction.\\
 Now the global symmetry of $G/K$ around the origin, is simply given by the exponentiation of the involution $X\mapsto -X$ of $\mathfrak{m}$.\\
The space $G/K$ always has non-positive sectional curvature (see \cite{Helgason}). There is the explicit formula for the sectional curvature
$$
K(Y,Z)=-||[Y,Z]||^2,\ Y, Z\in \mathfrak{m}=T_p(G/K),\ ||Y||=||Z||=1.
$$
So between any two points of $G/K$ there exists one and only one geodesic joining them.
\section{The Euler-Poincar\'e measure}\label{EP}
In this section we give the definition and an account of the fundamental properties of the Euler-Poincar\'e measure. Using it, we prove the equivalence of the statements $(4)$, $(5)$ and $(6)$ of theorem \ref{mainKtheorem}. The content is a reformulation of parts of \cite{Serre} and \cite{HopfSamelson}.\\
The idea is that the Euler class is represented by a differential form that, if lifted to $\EG$, is \emph{$G$-left-translation invariant}. The reason is that it only depends on the curvature of the symmetrical space, which is translation invariant.\\
\begin{defn}
A left-translation-invariant measure $\mu$ on a locally compact unimodular group $G$ is called an \emph{Euler-Poincar\'e} measure if every torsion-free discrete cocompact subgroup $\Gamma$ has a finite classifying cell complex $\BG$ and the Euler characteristic of $\Gamma$, which is therefore defined, is given by the volume of the quotient:
$$
\chi(\Gamma)=\int_{\Gamma\backslash G}\mu.
$$
\end{defn}
If $G$ has an Euler-Poincar\'e measure, then it is unique. So we may speak of \emph{the} Euler-Poincar\'e measure of $G$.\\
There are some more or less obvious, but not very interesting examples. If $G$ is already discrete and $\classs G$ is finite, then all torsion-free cocompact subgroups $\Gamma$ have finite $\BG$, and $\mu(g)=\chi(G)$ is an Euler-Poincar\'e measure. In fact, if $\Gamma$ is a subgroup of index $k$ then $\chi(\BG)=k\chi(\classs G)=k\chi(G)$ because the Euler characteristic of a $k$-sheeted cover is $k$ times that of the base.\\
Second, if $G$ is compact then the Haar measure normalized to $\mu(G)=1$ is an Euler-Poincar\'e measure because there are no non-trivial discrete torsion-free subgroups.\\
An interesting and to us relevant example is provided by
\begin{lemma}
A real semisimple connected Lie group has an Euler-Poincar\'e measure.
\end{lemma}
\begin{proof}
For $\Gamma$ discrete and cocompact, $\BG=\Gamma\backslash G/K$ for a maximal compact $K$ of $G$.\\
If $q=\dim\BG$ is odd, then obviously $\mu=0$ is an Euler-Poincar\'e measure. So let us assume $q$ is even.\\
The crucial part relies on the higher Gauss-Bonnet theorem.\\
Recall that the Chern-Weil approach to characteristic classes states that the Euler class of a compact Riemannian manifold can be realized explicitly as a differential form by the following procedure.\\
Let $U$ be a coordinate patch of a not necessarily compact Riemannian manifold with orthonormal frame $e_1,\dots,e_q$. Let $R$ denote the curvature tensor, a tensor of degree $(3,1)$ which satisfies
$$
R(u,v)e_i=\sum_j\Omega_{ij}(u,v)e_j
$$
for a family $\Omega_{ij},1\le i,j\le i$ of differential forms of degree 2. The forms $\Omega_{ij}$ depend alternatingly on $(i,j)$ and commute with each other since they are of degree 2.\\
Given $n\in\N$, there exists one and up to sign only one polynomial, that we shall call $\Pfaff$, with integer coefficients which assigns to the entries of each $2n$-by-$2n$ skew-symmetric matrix $A$ over a commutative ring a ring element $\Pfaff(A)\in A$ whose square is the determinant of $A$. $\Pfaff$ has degree $n$ and satisfies $\Pfaff(BAB^t)=\Pfaff(A)\det B$. Now it makes perfect sense to plug in the $i$-by-$i$ matrix $(\Omega_{ij})$ of two-forms on $U$, and the result is a differential form $\Pfaff ((\Omega_{ij}))$ of degree $2n=q$ on $U$. By theorem, this form is closed. If we put
$$
\Omega=\frac{1}{(2\pi)^{i/2}}\Pfaff((\Omega_{ij})).
$$
then it is the content of the general Gauss-Bonnet theorem that for any compact piece $A$ with boundary $\partial A$ of the Riemannian manifold we have
$$
\chi(A)=\int_A\Omega + \int_{\partial A} \Pi,
$$
where $\Pi$ is a certain differential form of degree $q-1$ on $\partial A$.\\
If $(e_1,\dots,e_q)$ is replaced by another orthonormal frame then $\Omega$ does not change, and is replaced by $-\Omega$ if orientation is reversed. So one may glue the $\Omega$'s on the different coordinate patches together and obtain a top degree differential form on $\BG$, in other words a measure, called \emph{Gauss-Bonnet measure}. If $G$ is a Lie group then by the translation-invariance, $\mu$ is a real multiple of the volume form, $\mu=\lambda d\text{Vol},\lambda\ge 0$. The theorem of Gauss-Bonnet then takes the form
$$
\chi(\Gamma)=\int_{\BG}\Omega,
$$
so one may also call $\Omega$ an ``Euler density''. Now apply this to the Riemannian manifold $\EG=G/K$. By the homogeneity of $\EG$ the curvature is $G$-left invariant, and so is the form $\Omega$, therefore. Since $K$ is compact, there is a unique invariant measure $\mu$ on $G$ such that the image of $\mu$ under $G\to G/K$ is $\Omega$. If $\nu$ is the normalized Haar measure of $K$ such that $\nu(K)=1$ then $\mu/\nu=\Omega$.\\
Now we claim that $\mu$ is an Euler-Poincar\'e measure on $G$. To show this, let $\Gamma$ be a discrete torsion-free cocompact subgroup of $G$. Now $\chi(\Gamma)=\chi(\BG)=\int_{\Gamma\backslash G/K}\Omega$ where it is justified to use the letter $\Omega$ for the Euler density on $\BG$ as well as on $\EG$. $\int_{\Gamma\backslash G/K}\Omega$ is equal to the integral of $\Omega$ over a fundamental domain of the $\Gamma$-action on $\EG$. Now this in turn is equal to the integral of $\mu$ over the inverse image of the fundamental domain under $G\to G/K$, which is a fundamental domain for the action of $\Gamma$ on $G$. We obtain $\int_{\Gamma\backslash G}\mu$, as required.
\end{proof}
So from the defining property of the Euler-Poincar\'e measure we deduce the equivalence of the statements (4) and (5), both of which may also be formulated as $\lambda>0$ for $\lambda$ as above.\\
$(5)\Leftrightarrow (6)$
The preceding proof shows that the Euler-Poincar\'e measure on $G$ is obtained by finding the unique invariant preimage under $G\to G/K$ of the Pfaffian of the curvature of $G/K$.\\
Let $\mathfrak{m}$ be the orthogonal to $\mathfrak{k}$ in the Lie algebra $\mathfrak{g}$ of $G$ with respect to the Killing form (here we need $G$ to be semisimple). $\mathfrak{m}$ can be identified canonically with tangent space to $G/K$ in the neutral element. Then the curvature tensor $R$ in that point is a certain tensor of degree $(3,1)$, more precisely it is an element of $\Lambda^2((\mathfrak{g}/\mathfrak{k})^*)\otimes \Hom(\mathfrak{g}/\mathfrak{k},\mathfrak{g}/\mathfrak{k})$. It can be expressed in terms of the Lie bracket by
$$
R(u,v)w=[[u,v],w],\quad u,v,w\in\mathfrak{m},
$$
see for example \cite{Helgason}*{theorem IV.4.2}.\\
Now let $\mathfrak{g}_\C=\C\otimes \mathfrak{g}=\mathfrak{g}\oplus i\mathfrak{g}$ be $\mathfrak{g}$'s complexification. If $\mathfrak{g}$ already has a complex structure, we ignore it. The decomposition of $\mathfrak{g}$ in $\mathfrak{k}\oplus\mathfrak{m}$ gives a decomposition of $\mathfrak{g}_\C$:
$$
\mathfrak{g}_\C=\mathfrak{k}\oplus \mathfrak{m}\oplus i\mathfrak{k}\oplus i\mathfrak{m}.
$$
Its Lie subalgebra $\mathfrak{g}'=\mathfrak{k}\oplus i\mathfrak{m}$ is another real form of $\mathfrak{g}_\C$ (recall that $[\mathfrak{m},\mathfrak{m}]\subset\mathfrak{k}$). $\mathfrak{g}'$ is a \emph{compact} Lie algebra. In order to verify this, we have to assure that the Killing form is strictly negative definite. We use the fact that if $\mathfrak{h}$ is a complex Lie algebra with Killing form $\kappa$, and $\kappa^\R$ is the Killing form of $\mathfrak{h}$ viewed as a real Lie algebra, then $\kappa$ and $\kappa^\R$ are related by the formula
$$
\kappa^\R(X,Y)=2\thinspace\thinspace\Re (\kappa(X,Y)),\quad X,Y\in \mathfrak{h}.
$$
Now let $X\in\mathfrak{k}$, $Y\in\mathfrak{m}$,$X\ne 0,Y\ne 0$. Then $\kappa^\R(X+iY,X+iY)=2\thinspace\thinspace\Re(\kappa(X+iY,X+iY))=2\thinspace\thinspace\Re(\kappa(X,X)-\kappa(Y,Y))=2\kappa(X,X)-2\kappa(Y,Y)<0$, so $\mathfrak{g}'=\mathfrak{k}+i\mathfrak{m}$ is a compact form of $\mathfrak{g}_{\C}$. Let $G'$ be the unique simply connected Lie group with algebra $\mathfrak{g}'$. Furthermore, let $K'$ be the subgroup of $G'$ corresponding to $\mathfrak{k}\subset \mathfrak{g}'$. $G'$ and $K'$ are compact, and $G'/K'$ is sometimes called the dual of $G/K$ (see \cite{Helgason}, V.2). $i\mathfrak{m}$ is the tangent space to $G'/K'$ at the neutral element, and the curvature tensor $R'$ of $G'/K'$ is given by the same formula
$$
R(u,v)w=[[u,v],w],\quad u,v,w\in i\mathfrak{m}.
$$
It follows that the bijection $u\mapsto iu$ maps $R$ in $-R'$. Let $\Omega$ and $\Omega'$ be the value of the Gauss-Bonnet form of $G/K$ resp. $G'/K'$. Now $u\mapsto iu$ transforms $\Omega$ to $\sign\Omega'$. But the integral of $\Omega'$ over $G'/K'$ is equal to $\chi(G'/K')$. Now the statement follows from the following old result of Hopf-Samelson (\cite{HopfSamelson}).
\begin{thm}
If $G$ is a \emph{compact} connected Lie group, then the Euler characteristic of a homogeneous space $G/U$ is non-zero if and only if the ranks of $G$ and $U$ agree. In that case, it is in fact positive.
\end{thm}
We are now going to give an account of their proof.
\begin{lemma}
Let $W$ be a homogeneous space of the compact connected Lie group $G$. Then any element of $G$ whose action on $W$ has only finitely many fixed points in $W$ has in fact $\chi(W)$ many.
\end{lemma}
In particular, if there is such a group element, then the Euler characteristic is non-negative.
\begin{proof}
Since $G$ acts by isometries, any isolated fixpoint of $g\in G$ can only have negative eigenvalues (a positive one would have to be 1, thus the fixpoint wouldn't be isolated). So the index of $g$ at the fixpoint is $(-1)^{dim W}$. Here, the index of a smooth map $f:M\to M$ at $p\in M$ is defined to be the integer by which $f_*:\Homol_{\dim M}(M,M-p;\Z)\cong\Z\to \Homol_{\dim M}(M,M-f(p);\Z)\cong\Z$ is multiplication. On the other hand, the sum of the indices of a continuous self-map of a finite complex $P$ that has only finitely many fixed points is $(-1)^{\dim P}\chi(P)$, if the map is homotopic to the identity. So $\sum_{\# \text{ fixpoints}} (-1)^{\dim W}=(-1)^{\dim W}\chi(W)$. We deduce that the number of fixpoints equals the Euler characteristic. If the dimension of $W$ is odd then its Euler characteristic vanishes. So the statement follows.
\end{proof}
{\it Proof} (of the Hopf-Samelson theorem).
We use the following facts on the tori of $G$ (\cite{BorelTori}). The automorphism group of a torus is discrete, and every torus has ``generating'' elements, i. e. elements whose powers are dense in the torus. Furthermore, a maximal torus has finite index in its normalizer, the quotient is called the Weyl group. Furthermore, the conjugates $gTg^{-1},\ g\in G$ cover $G$. Any two maximal tori are conjugate. So the rank of $G$ is well-defined as the dimension of a maximal torus.\\
We are now going to show the existence of a group element with only finitely many fixpoints, and then count them. Changing the viewpoint from the homogeneous space to the Lie group, we have to show that for any closed subgroup $U$ of $G$ there is $a\in G$ such that we have $xax^{-1}\in U$ for only finitely many cosets $xU$. We claim that a generating element $a\in G$ of a maximal torus $T\subset G$ does the job. So we have to show that $x^{-1}Tx\subset U$ for only finitely many cosets $xU$.\\
Assume that the rank of $U$ is smaller than that of $G$. Then not one of the cosets $xU$ satisfies the relation, $\chi(G/U)=0$ and we are done.\\
Assume that the ranks are equal. We may assume that $T\subset U$. We claim that $xU$ is a fixpoint of $a\in G$ iff it contains an element of $N_T$, the normalizer of $T$. Let $x$ be an element as required. Then $T$ and $x^{-1}Tx$ are conjugate in $U$, so there is even $u\in U$ with $uTu^{-1}=x^{-1}Tx$, so $xu\in N_T$. If, on the other hand, $y\in xU\cap N_T$, then $x$ conjugates the torus into $U$ by $x^{-1}Tx=u^{-1}y^{-1}Tyu=u^{-1}Tu\subset U$.\\
Now the finiteness of the Weyl group $N_T/T$ implies that there are only finitely many cosets $yT$ where $y\in T$. Consequently there are only finitely many cosets $xU$ with the required property when $y\in xU$. So we are done and note that in that case $\chi(G/U)\ge 1$ because $U$ itself represents a fixpoint.\qed\\
\section{The foliation index theorem}\label{foliation}
In this appendix we show how Connes' foliation index theorem (\cite{Co2}, III.$7.\gamma$) would give the integrality if it applied. Since this theorem is a very advanced machinery, and it is not clear if its conditions hold, we proceed differently in the body of the text. This section is therefore unnecessary for the logical coherence of the thesis. We merely show how Connes' foliation index theorem would yield the integrality, if it could be applied (and if we knew that the K-theory classes in the crossed product $\Gamma\rlcross C(G/B)$ that can come from $\C\Gamma$ are in the image of $\K^*(\BG)$ in $\K^*(\EG\times_\Gamma G/B)\cong\K^(\BG)\otimes K^*(G/B)$, which is an extremely natural, but seemingly hard to prove statement).\\
Suppose that $\BG$ is compact. A suitable formulation of the latter is given by Lott (\cite{Lott}, 1.2). His formula, boiled down to our notations, reads
$$
<\Ind \Dirac_{W},\Phi(\eta)>=\int_{\EG\times_\Gamma G/B}\pi^*\Aroof(\BG)\cup\ch(W)\cup\eta.
$$
where $\Dirac$ is the Dirac operator on $\EG\times_\Gamma G/B$ (we prove in remark \ref{SPINN2} that a Spin-structure exists), $\Dirac_W$ is the Dirac operator twisted with a complex vector bundle $W$ on $\EG\times_\Gamma G/B$, $\Ind\Dirac_W$ is its index in $\K_0(\mathcal{A})$ (where $\Gamma\allcross\CINF(G/B)\subset\mathcal{A}\subset\Gamma\rlcross C(G/B)$, see below), $\eta\in\Homol^{\dim G/B}(\EG\times_\Gamma G/B)$ is a cohomology class, $\Phi(\eta)\in\HC^0(\mathcal{A})$ is its associated cyclic cohomology class and $\Aroof$ is the $\Aroof$-genus. If we knew that all K-theory classes in the smooth crossed product that can come from $\K_0(\C\Gamma)$ are those which are indices of Dirac operators twisted with a vector bundle of the form $W=\pi^* V$ then we would have, taking for $\eta$ the preimage $\eta=\PD(\pi^*[\BG])$ of the trace $\tau\in\HP^0(\C\Gamma)$:
\begin{align*}
\int_{\EG\times_\Gamma G/B}\pi^*\Aroof(\BG)\cup\ch(\pi^*V)\cup\eta&=\int_{\EG\times_\Gamma G/B}\pi^*\Aroof(\BG)\cup\ch(W)\cup\PD\pi^*[\BG]\\
&=\int_{\BG}\Aroof(\BG)\cup\ch(V)\cup\pi_!(\PD(\pi^*[\BG]))\\
&=\int_{\BG}\Aroof(\BG)\cup\ch(W)\cup 1=\Ind(\Dirac^{\BG}_V)
\end{align*}
where the latter index is the ordinary Atiyah-Singer index of the Dirac operator $\Dirac^{\BG}$ on $\BG$ twisted with $V$, which is an integer.\\
However, it is not clear if this theorem is applicable because it is difficult to decide whether there are any smooth subalgebras $\Gamma\allcross\CINF(G/B)\subset\mathcal{A}\subset\Gamma\rlcross C(G/B)$ between the algebraic and $C^*$-algebraic crossed products such that $\mathcal{A}$ is closed under holomorphic functional calculus in its $C^*$-algebra completion $\Gamma\rlcross C(G/B)$ (such that both algebras have the same K-theory), and such that the cyclic cocycle $\tau'$ lifts to $\mathcal{A}$.
\section[A deformation version of the smooth Dirac-dual Dirac method]{An alternative proof of theorem \ref{XXX} by a deformation version of the Dirac-dual Dirac construction}
The goal of this appendix is to give an alternative proof of theorem \ref{XXX}. It is based to a large extent on a deformation proof in K-theory of $C^*$-algebras, as it has been described by Higson and Roe (\cite{Summerschool}). The difference is that we use smooth subalgebras and $\HA$, which requires some extra work.\\
Let $A$ be an arbitrary bornological algebra endowed with a bounded action of the torus $T$. Our goal is the definition of Dirac- and dual Dirac-morphisms
$$
\beta\in\HA(A\alcross T,(A\protens\CCINF \R^2)\alcross T)\\
\alpha\in\HA((A\protens\CCINF \R^2)\alcross T,A\alcross T)
$$
such that $\beta\otimes\alpha=1$.
\subsection{The definition of $\beta$}
This requires some preliminary work. Recall that for a smooth non-compact manifold $M$ the algebra $\CCINF M$ of smooth functions with compact supports is as a bornological algebra defined as follows: A subset $W$ of $\CCINF M$ is called bounded if there is a compact $K\subset M$ such that all functions $f\in W$ have support in $K$ and all derivatives $f^{(n)}$, $f\in W$ are uniformly bounded on $K$.\\
So $\CCINF M$ is the direct limit
$$
\lim_{K\subset M} SK
$$
taken over all smooth compact submanifolds with boundary, of the space $SK$ of functions that vanish on $\partial K$ in arbitrary order. $\CCINF M$ is an LF-algebra (\cite{Meyerthesis}).\\
Consider now the algebra $\S=\CCINF \R$ with its $\Z_2$-grading $(\epsilon.f)(t)=f(-t)$ as well as the associated crossed product $\S\alcross \Z_2$.
\begin{lemma}
Denote by $\Cone$ the algebra of smooth functions on $[0,\infty]$ which assume in 0 arbitrary values, without a restriction on the derivatives, and which have compact support. Endow $\Cone$ with a bornology analogous to that of $\S$. Then $\Cone$ is $\HA$-contractible.\\
Furthermore, $\S$ is, as expected, in analytic cyclic homology the reduced cohomology of the circle.
\end{lemma}
\begin{proof}
$\Cone$, thus defined, is the bornological direct limit $\lim_K\S[0,K]$ of the algebras of functions on compact intervals all of whose derivatives vanish at K, but not in 0. As $\HA$ is continuous with respect to such a limit, it suffices to compute $\HA \S[0,K]$. But $\S[0,K]$ is easily seen to be contractible through a smooth path. For the second statement, observe that there is an extension of $\C$ by  $\S[0,K]$ with quotient $\S$. Thus the statement follows from excision.
\end{proof}
\begin{thm}\mbox{ }\\
\vspace{-1.5em}
\begin{enumerate}
\item
There is an invertible in $\HA(\S\alcross\Z_2,\C)$ defined by the composition
$$
\xymatrix{\S\alcross\Z_2\ar[r]^>>>>>{\ev_0}&\C\alcross \Z_2\cong \C\oplus \C\ar[r]^<<<<{\pi_1}&\C}
$$
where $\pi_1$ is the projection onto the first factor, uniquely defined by the projector $\frac{1}{2}(\delta_0+\delta_1)\in\C\alcross\Z_2=\C\delta_0\oplus\C\delta_1$ (in fact, we have to fix a choice, but it is arbitrary) onto the even elements.
\item
The same is valid for $A\protens (\pi_1\circ (\ev_0\alcross \Z_2))$, i. e. the corresponding element in $\HA(A\protens \S\alcross\Z_2,A)$ is invertible.
\end{enumerate}
\end{thm}
\begin{proof}\mbox{ }\\
\vspace{-1.5em}
\begin{enumerate}
\item
Recall that $\mathcal{C}$ is the algebra of smooth complex-valued functions on the half line $\R^+$ whose derivatives in $0$ assume all possible values. Then the algebra $\S\alcross\Z_2$ is isomorphic to
$$
\{f\in\mathcal{C}M_2:f(0) \text{ has diagonal form}\}
$$
in such a way that
$$\ev_0\alcross\Z_2:\S\alcross\Z_2\to\C\Z_2\cong\C\oplus \C$$
corresponds to evaluation at $0$ taking values in $\diag(\C,\C)\cong\C\oplus\C$.\\
Let us prove the following claim: The kernel of $\pi_1\circ (\ev_0\alcross\Z_2)$ is $\HA$-contractible. In fact, under the above isomorphism it corresponds to the algebra
$$
V=\{f\in\mathcal{C}M_2:f(0)\text{ has form } \diag(0,*)\}.
$$
 By definition, there is an inclusion $V\hookrightarrow \mathcal{C}\otimes M_2$. Consider the composition, called $i_1$,
$$
V\hookrightarrow \Cone M_2\hookrightarrow
\biggl\{
f\in\Cone M_4:f(0)\text{ has the form }
{\tiny
\left(
\begin{array}{cccc}
* & * &  0&0\\
* & * & 0 &0\\
0 & 0 & 0 &0\\
0&0&0&0
\end{array}
\right)
}
\biggr\}
$$
as well as the standard upper-left-corner embedding $i_2:V\to M_2 V$. Now it is easily seen, using a smooth path $\varphi_t$ inside the smooth compact connected manifold $U(4)$ between $\varphi_0={\tiny
\left(
\begin{array}{cccc}
1&0&0&0\\
0 & 0&1&0\\
0 & 1&0&0\\
0 & 0&0&1
\end{array}
\right)
}$ and $\varphi_1=\id_{\C^4}$ that $i_1$ is smoothly homotopic to $i_2$. In fact, upper-left-corner-embed $V$ into $M_2V\subset \Cone M_4$, conjugate at time $t\in[0,1]$ $V$ inside $\Cone M_4$ with$\varphi_t$ and observe that the result is still contained in $M_2V$ (which is not conjugated). For $t=0$, this gives $i_1$, for $t=1$, it gives $i_2$.\\
Now $i_1$ induces $0\in\HA(V,M_2V)$ because it factors through the contractible algebra $\Cone M_2$. On the other hand, by Morita equivalence, $i_2$ induces an $\HA$-equivalence, hence the claim.\\
By excision, the result is thus proved.
\item
The proof of the first part of the proposition is written in such a way that it applies in the case of $A$ replacing everywhere $\C$ as well.
\end{enumerate}
\end{proof}
Consider now the algebra $\CCINF\R^4\otimes M_4$ and view it as the algebra of smooth compactly supported sections of the trivial bundle over $\R^4=T^*\R^2$ with fiber $\End \Lambda^*_\C T^*\R^2$ where $T^*\R^2$ is viewed as $\C^ 2$ - with its canonical hermitian scalar product - by taking the horizontal directions as real, the vertical ones as imaginary.Now let us suppose that a fixed linear action of the torus $T=\mathbb{S}^1\times\dots\times\mathbb{S}^1$ on $\R^2=\C$ is given. Then $\alTor$ acts on $\R^4=\C^2$ on each factor $\C$ separately, preserving the scalar product, and on $\End \Lambda^*_\C T^*\R^2=\Cliff \R^4$ by the universal property of the Clifford algebra. We define the action on $\CCINF\R^4\otimes M_4$ to be the tensor product action.\\
\begin{lemmaanddef}[The dual-Dirac-element]\label{dD}
For $\theta\in\R^4=T^*\R^2$, let $\theta\wedge-:\Lambda^*_\C T^*\R^2\to\Lambda^{*+1}_\C T^*\R^2$ denote the exterior multiplication operator on $\Lambda^*_\C T^*\R^2$, and let $\theta\vee-:$ be its adjoint with respect to the obvious hermitian scalar product on $\Lambda^*_\C T^*\R^2=\C^4$.\\
Then the formula
$$
b(\theta)(x)=\theta\wedge x+\theta\vee x, \theta \in T^*\R^2, x\in
\Lambda^*_\C T^*\R^2
$$
then defines a real-linear map 
$$
\R^4\to \End \Lambda^*_\C T^*\R^2.
$$
It passes to an unbounded, formally self-adjoint multiplier of the $C^*$-alge-braic completion $C_0(\R^4\otimes M_4)S$ of $\CCINF(\R^4)\otimes M_4$. This means that functional calculus 
$$C_0(\R)\to C_0(\R^4)\otimes M_4,\ f\mapsto f(b)
$$
is well-defined. It is standard that $b(\theta)^2=||b||^2 \id$ which implies that functional calculus restricts to a map
$$
\CCINF(\R)\to\CCINF \R^4\otimes M_4.
$$
Endow furthermore the algebra $M_4=\End \Lambda^*_\C T^*\R^2$ with a $\Z_2$-grading by declaring those endomorphisms which preserve $\Lambda^\ev$ and $\Lambda^\odd$ to be even, those which interchange them to be odd. Then the map $\CCINF\R\to \CCINF\R^4\otimes M_4$ defined above is $\Z_2$- as well as $T$-equivariant (i. e. it is graded and the image is $T$-invariant) and thus defines a homomorphism
$$
\alcross T\protens\S\alcross\Z_2\cong A\protens\S\alcross\Z_2\alcross
T\to A\protens \CCINF\R^4\otimes M_4\alcross\Z_2\alcross T
$$
and therefore, via the invertible in $\HA(A\alcross T\protens\S\alcross\Z_2,A\alcross T)$ an element
$$
\beta'\in\HA(A\alcross T,A\protens\CCINF\R^4\otimes M_4\alcross\Z_2\alcross T).
$$
Furthermore, $A\protens\CCINF\R^4\otimes M_4\alcross\Z_2\alcross T$ is isomorphic to $A\protens\CCINF\R^4\alcross\Z_2\otimes M_4\alcross T$ which in turn is Morita-equivalent to $(A\protens\CCINF\R^4\alcross T)^{\oplus 2}$.\\
By definition, the composition of $\alpha'$ with this isomorphism, the Morita equivalence isomorphism, and the projection onto the first factor (the proof will make clear that, analogous to the case appearing in the previous proposition, there is a canonical way to choose the projection) gives the dual-Dirac element.
$$
\beta\in\HA(A\alcross T,A\protens\CCINF\R^4\alcross T).
$$
\end{lemmaanddef}
\begin{rmk}
At first sight it might appear that this dual Dirac element belongs to a weaker statement then Bott periodicity in that it seems to yield only a $\Z_4$ instead of a $\Z_2$-periodicity, since $\alpha$ as such is not defined in $\HA(A\alcross T,A\protens\CCINF\R^2\alcross T)$. It is, however, easy to derive starting from $\alpha$ a cycle in the latter group by the following consideration: There is an isomorphism $\R^4\cong\C^2$ such that $T$ acts on the first factor complex-linearly, and trivially on the second.\\
Thus, $A\protens \CCINF\R^4\alcross T$ is isomorphic to $(A\protens\CCINF\C\alcross T)\protens \CCINF\C.$ Now observe that for any complete bornological algebra $B$ there is the suspension extension
$$
0\to\S\protens B\to\Cone\protens B\to B\to 0
$$
which yields an invertible in $HA_1(\S\protens B, B)$ because $\Cone\protens B$ is contractible, and thus an invertible in
$\HA_0(\CCINF\R^2\protens B,B)$, which therefore reduces $\alpha$ to a cycle in $\HA_0(A\alcross T, A\protens\CCINF\R^2\alcross T)$.
\end{rmk}
\begin{proof}(of the statements made during the definition.)\\
Let us show that the map $\S\to \CCINF\R^4\otimes M_4$ is
$\Z_2$-equivariant.\\
We may assume that $||\theta||=1$. If $f$ is even, then, since there is only a $+1$ and a $-1$-eigenspace, $f(b(\theta))=f(1)\id$, which is an even operator. Observe that $b(\theta)$ must have the form ${\tiny  \left(\begin{array}{cc}0&U^*\\U&0\end{array}\right)}$ for a unitary isomorphism $U:\Hil^+\to\Hil^-$. If $f$ is odd, then $f$ maps the $+1$-eigenspace $\{(\xi,U\xi)/\xi\in\Hil^+\}$ to $f(1)$ times itself, and the $-1$-eigenspace $\{(\xi,-U\xi)/\xi\in\Hil^+\}$ to $-f(1)$ times itself, from which it follows that in fact $\Hil^+=(\xi,0)$ is mapped to $\Hil^-$ and vice versa.\\
Let us now show that the map $\S\to\CCINF\R^4\otimes M_4$ is $T$-equivariant, i. e. that the image is invariant. The universal property of the Clifford algebra means that to any map $V\to W$ of Euclidean spaces there is an induced algebra map $\Cliff V\to\Cliff W$ such that
$$\xymatrix{
V\ar[r]\ar[d]&\Cliff V\ar[d]\\
W\ar[r]&\Cliff W
}
$$
commutes. Applied to the representation maps $\R^4\to\R^4$ for an element of $T$, this yields $b(t\inv x)=t.b(x)$, since $b:\R^4=T^*\R^2\to\End \Lambda^*_\C T^*\R^2=\Cliff \R^4$ is exactly the inclusion $\R^4\hookrightarrow \Cliff\R^4$. So $b\in C(\R^4)$ is $T$-invariant, and so is any operator obtained from it by functional calculus.\\Now let us show that there is a $T$-equivariant Morita-equivalence between $\Cliff V$ and $\C$. Since the action of $t$ is complex-linear, it preserves a $\Spin^c$-structure on $\R^4$, which implies that $\Cliff V\cong\End\Spinor$, $V$ being the trivia vector bundle over $\R^4$ with fiber $\R^4$. Now $T$ acts by conjugation on $\End\Spinor$, and so it does it act on $\Spinor\oplus$ the trivial one-dimensional bundle with trivial action. The smooth sections with compact supports of this direct sum provide therefore an imprimitivity bimodule.\\
The dual-Dirac-element is thus defined.
\end{proof}
\subsection{The definition of $\alpha$}
In order to define the Dirac-element, consider the $t\in[0,1]$-dependent action of $x\in\R^2$ on itself, given by $x.y=tx+y, y\in\R^2$. The algebra $\CCINF\R^4\otimes\CINF[0,1]$, where now $\CINF[0,1]$ is the algebra of all functions all of whose derivatives except possibly the zeroeth vanish at both endpoints.\\
Let $\Sch \R^n$ denote the Schwartz-algebra of rapidly-decreasing functions $\R^n$.\\
The vector space $\Sch\R^2\protens\CINF[0,1]\protens\Sch\R^2$ can be endowed with a crossed product, yielding a complete bornological algebra
$$
\Defo=\Sch\R^2\protens\CINF[0,1]\rtimes\R^2
$$
where now $\rtimes$ denotes a slightly larger crossed product than the algebraic $\CCINF$ one.\\
In other words, consider the smooth convolution algebra of the groupoid $(\R^2\times\R)\rtimes \R^2$, where now $t$ runs over the real line, enlarge it to have $\Sch\R^5$ as underlying vector space, take the subalgebra of functions $f$ which satisfy $\frac{\partial^n}{\partial t^n}f=0$ on the set $\R^2\times\{0,1\}\times\R^2$ (a subalgebra which is isomorphic to $\Sch\R^5$ itself using a suitable diffeomorphism of $\R$) and take the quotient with respect to the ideal of all functions which vanish on $\R^2\times[0,1]\times\R^2\subset\R^5$.\\
\begin{lemma}
There is an extension
$$
0\to \smKo\protens\Sch(0,1]\to\Defo\to\Sch\R^4\to 0
$$
where the kernel is the algebra of smooth compact operators (those compact operators whose sequence of characteristic values is rapidly decreasing) tensored with the algebra of rapidly decreasing functions on $(0,1]$.
\end{lemma}
\begin{proof}
We define the extension by the ideal of functions inside $\Defo$ which vanish on $\R^2\times\{0\}\times\R^2$. This ideal is the algebra of sections of an algebra bundle of Schwartz-convolution algebras on the groupoid $\R^2\rtimes_t\R^2$ with the action over the parameter $t$. Over $t=1$, we obtain $\Sch\R^2\rtimes\R^2\cong\smKo$ by the standard matrix multiplication or integral kernel representation on the Hilbert space $L^2\R^2$. However, each single fiber is isomorphic to the fiber over $t=1$ by rescaling. In fact, it is the algebra of smooth compact operators on the Hilbert space $L^2\R^2$ with the measure $d\mu'=t^2d\mu$ instead of the Lebesgue measure $d\mu$. But all these Hilbert spaces are unitarily equivalent. Denote for further use the $(0,1]$-family of fiber rescaling isomorphisms by $\varkappa_t$.\\
The quotient is $\Sch\R^2\protens\Sch\R^2$, now even as an algebra because for $t=0$ the action is trivial, where the second $\Sch\R^2$ has convolution product. This is, however, isomorphic to $\Sch\R^4$ with pointwise product by Fourier transformation.
\end{proof}
We now want to endow all algebras with torus actions.\\
The notion of smooth action of Lie group on a groupoid is defined in the natural way. The group of automorphisms of the groupoid $G$ is defined to be the group of all natural isomorphisms $G\to G$ which are smooth maps on the morphism and object manifolds. A representation of a compact Lie group is a smooth group homomorphism into the automorphism group.\\
The groupoid $(\R^2\times\R)\rtimes R$ is endowed with a $T$-action by setting
$$
\gamma(x,t)=(\gamma x,t)
$$
on the object space $\R^2\times\R$ and $$\gamma(x,t,y)=(\gamma x,t,\gamma y)$$ on the morphism space $\R^2\times \R\times \R^2$. After checking that this defines in fact a action on the groupoid in the above sense, one can immediately endow all of the above algebras with a $T$-action.\\
We thus may consider the induction algebra $A\protens\Defo\alcross T$ for any complete bornological $T$-algebra $A$. $A\protens\Defo\alcross T$ is an extension because $\Defo$ is nuclear, being itself an extension of two nuclear algebras. There is a Morita equivalence $A\protens\smKo\alcross T\equiv A\alcross T$ given by the smooth groupoid $(\R^2\rtimes \R^2)\coprod \pt$ (disjoint union with a point).
\begin{defn}[The Dirac element, deformation version]
There is an element, called the Dirac element
$$
\alpha\in\HA(A\protens \CCINF\R^4\alcross T,A\alcross T)
$$
given by the composition of the following elements.
\begin{enumerate}
\item
The inclusion $A\protens\CCINF\R^4\subset A\protens \Sch \R^4$
\item
The inverse of the elements
$$
A\protens p \alcross T\in\HA(A\protens\Defo\alcross T,A\protens
\Sch\R^4\alcross T)
$$
which is invertible since the ideal is contractible.
\item
The element $$A\protens\ev_1\alcross T\in\HA(A\protens\Defo\alcross T,A\protens\smKo\alcross T)$$
\item
The (invertible) Morita element $\in \HA(A\protens\smKo\alcross T,A\alcross T)$.
\end{enumerate}
\end{defn}
\subsection{Bott periodicity, deformation version}
The contents of this section are very close to Higson and Roe's lecture notes (\cite{Summerschool}).
\begin{thm}\label{Bottdefo}
$\beta\otimes\alpha=1\in\HA(A\alcross T,A\alcross T)$
\end{thm}
\begin{proof}
For a $C^*$-algebra $B$, define $AB$ to be $C_b(0,1]\otimes B$ ($C^*$-algebraic tensor product, here being unique), the algebra of bounded functions from $(0,1]$ to $B$, and $CB$ to be the ideal $C_0(0,1]\otimes B$, i. e. the cone over $B$.\\
Consider the $C^*$-completion
$$
0\to C\Ko\to C_0(\R^2\times[0,1])\rtimes \R^2\to C_0(\R^4)\to 0
$$ 
of the extension introduced above, and choose a linear and continuous section $\varsigma$ which restricts to a linear bounded section of the smooth extension.\\
Define a linear asymptotic morphism $\alpha_t:C_0(\R^4)\to\Ko$ by the formula $\alpha_t=\varkappa_t\circ\ev_t\circ \varsigma$, where $$\varkappa_t:C_0(\R^2)\rtimes_t\R^2\to C_0(\R^2)\rtimes \R^2\cong\Ko$$  is the ($C^*$-completion of) the rescaling isomorphism introduced above, and $\ev_t$ is, of course, evaluation at $t\in(0,1]$
$$
C_0(\R^2\times[0,1])\rtimes \R^2\to C_0(\R^2)\rtimes_t\R^2
$$
to the crossed product involving the $t$-action. The definition does not depend on the choice of the section.\\
As does any asymptotic morphism, $\alpha_t$ defines a $*$-homomorphism $C_0(\R^4)\to A\Ko/C\Ko$ by composing the linear $\alpha_t$, viewed as a linear map $C_0(\R^4)\to A\Ko$, with the quotient map. Consider also the complete bornological tensor product $A\protens B$ resp. $C\protens B$ of the $C^*$-algebras $A$ resp. $C$, endowed with the precompact bornology, with some complete bornological algebra $B$. Identify $\Ko\otimes M_4$ with the compact operators on $L^2(\R^2,\Lambda^*_\C T^*\R^2)$, the graded Hilbert space of $L^2$-sections of the constant hermitian bundle with fiber $\C^4=\Lambda^*_\C T^*\R^2$ over $\R^2$.\\
Define a family of elliptic order one differential operators on this bundle by setting
$$
D_t=d_t+d_t^*,\ t\in(0,1]
$$
where $d_t\omega=\eta\wedge\omega+td\omega$, $\eta$ the differential of the function $\frac{1}{2}||x||^2$ on $\R^2$. Now Higson shows that
$$
(\alpha_t)\circ\beta:C_0(\R)\to C_0(\R^4)\otimes M_4\to A\Ko/ C\Ko
\otimes M_4
$$
lifts to the $*$-homomorphism $C_0(\R)\to A\Ko\otimes M_4=A(\Ko\otimes M_4)$ defined by functional calculus $f\mapsto f(D_t)$. In other words,
$$
\xymatrix{
&A\Ko\ar[d]\\
C_0(\R)\ar[ru]^{f\mapsto f(D_t)}\ar[r]_{\alpha\circ\beta}&A\Ko/C\Ko}
$$
commutes.\\
Fix a diffeomorphism from $(0,1]$ to $(-\infty,0]$.\\
Define furthermore a smooth version of $A$ for a complete bornological algebra $B$ as $\CINF_b(0,1]\protens B$, where $\CINF_b(0,1]$ is the algebra of smooth and bounded functions on $(0,1]$, whose lim sup of the derivatives grow at most polynomially under the diffeomorphism. Denote that algebra by $\smA$. Analogously, a smooth cone ideal $\Cone B$ is defined by functions all whose derivatives vanish at $0$. This is in fact an ideal since under the diffeomorphism this algebra obeys exactly the Schwartz decay condition. (Previously, $\Cone$ denoted a slightly different version of the cone.)\\
Now define a smooth version of the above diagram as follows.
$$
\xymatrix{&&\smA\smKo\ar[d]\\
\S\ar[rru]^{f\mapsto f(D_t)}\ar[r]_>>>>>{\beta}&\CCINF\R^4\otimes
M_4\ar[r]_>>>>{\alpha_t}&\smA\smKo/ \Cone\smKo\otimes M_4
}
$$
Let us first see how the maps are defined and then why it commutes.\\
$\beta$ is the map already defined. In order to make the asymptotic morphism $(\alpha_t)$ descend on a map as indicated, recall that there is chosen a smooth linear lift of
$$
\Sch(\R^2\times[0,1])\rtimes \R^2\to \Sch \R^4
$$
(which exists because the left hand side algebra is just $\Sch(\R^4\times [0,1])$. Then the composition $\CCINF\R^4\hookrightarrow \Sch\R^4$ with this lift, $\ev_t$ and $\varkappa_t$ gives the desired linear map $\CCINF \R^4\to \smA\smKo$.\\
For the definition of functional calculus $f\mapsto f(D_t)$, observe that $D_t$ has discrete
spectrum of finite multiplicities.\\
In other words, we have established the existence of the last diagram. Regarding its commutativity, there is a map from it to the $C^*$-diagram, which shows that it commutes up to the kernel $(\smA\smKo\cap C\Ko)/\Cone\smKo$ of this map. Here, the intersection is taken in the $C^*$-algebra $A\Ko$.\\
We now claim that it commutes in fact strictly, not just up to that subalgebra. So we want to show that
$$
\alpha_t\circ \beta(f)=\pi(f(D_t)).
$$
where $\pi$ is the obvious projection. More precisely, we demand the formula $\alpha_t\circ\beta(f)\cong f(D_t)$ to hold in $\CINF_b(0,1]\protens\smKo$ modulo the ideal $\Cone\protens\smKo$.\\
We recall that we cannot demand the equality strictly since $\alpha_t$ is only defined up to the ideal $\Cone\protens\smKo$ by $\pi\circ\varkappa\circ\zeta$.\\
In other words, we have to prove that $\kappa_t^{-1}(f(D_t))$ can be completed to a section $\zeta$ in $\ev_0:\Defo\otimes M_4\to\Sch(\R^4)\otimes M_4$ of (the Fourier transform in two variables of) $\beta(f)=f(b)\in\Sch(\R^4)\otimes M_4$, more
precisely that there exists an element $E_t\in \Defo$ with the properties
$$
E_0=f(b)
$$
and
$$
\varkappa_t(E_t)=f(D_t)
$$
for $t\in(0,1]$.\\
We shall exhibit such a family $E_t$ in $\Defo$ in the case of half dimensions, i. e. where $f(b)\in\Sch(\R^2)\otimes M_2$ etc., check the announced properties and then indicate which changes are necessary for the higher dimensional cases.\\
It is clear, however, that the central issue, namely that there is the kernel
$$
(C_0(0,1]\otimes\Ko\cap\CINF_b(0,1]\protens\smKo)/\Cone\protens \smKo
$$
of the map from the smooth diagram to the $C^*$-algebraic one, already arises in the low-dimensional case.\\
So define a family $E_t$ of operators on the real line by the smooth expression
$$
E_{\Psi^{-1}(t)}=f\left(
\begin{array}{cc}
0&m+t(x-m)+d\\
m+t(x-m)-d&0
\end{array}
\right)
$$
$\Psi$ being a suitable diffeomorphism of $[0,1]$ that assures that the derivatives vanish, $d$ being $d/dx$, $x$ and $m$ being the coordinates constituting the crossed product. $x$ belongs to the first, $m$ to the second $\R$. Now under the standard diffeomorphism from $\R\times [0,1]\times\R$ to the tangent groupoid $\mathbb{TR}$ of $\R$, restricted to $t\in[0,1]$, $E_{\Psi^{-1}(t)}$ is transformed to
$
f\left(\begin{smallmatrix}0&x-td\\x+td&0\end{smallmatrix}\right)
$
 for $t\ne 0$ and $f\left(\begin{smallmatrix}0&m+d\\m-d&0\end{smallmatrix}\right)$ for $t=0$. Now standard theorems from partial differential equations theory say that all $f(E_t)$ are represented by smooth kernel functions in $\Sch(\R)\rtimes_t \R$. In double dimensions just replace $\left(\begin{smallmatrix}0&x-td\\x+td&0\end{smallmatrix}\right)$ etc. by $D_t=d_t+d_t^*$. This shows that the diagram commutes strictly.\\
Consider furthermore the commutative diagram 
$$
\xymatrix{
\Sch(\R^4)\ar[d]^=&&\Sch(\R^2\times[0,1])\rtimes\R^2\ar[d]^\varkappa\ar[ll]\ar[r]^>>>>{\ev_1}&\smKo\ar[d]^=\\
\Sch(\R^4)\ar[r]^>>>>\alpha&\smA\smKo/\Cone\smKo&\smA\smKo\ar[l]\ar[r]&\smKo
}
$$
which shows that the deformation of the compacts defines the same map $\alpha$, viewed as a morphism in $\HA$, as the asymptotic morphism.\\
Altogether, the second last diagram commutes, too. Now look at the composition
\begin{multline*}
\S\to\CCINF\R^4\otimes M_4\hookrightarrow \Sch\R^4\otimes M_4=\\
\Sch(\R^2\times[0,1])\rtimes\R^2/\smKo(0,1]\otimes M_4\xleftarrow{\cong}\Sch(\R^2\times[0,1])\rtimes\R^2\xrightarrow{\ev_1}\smKo\otimes M_4
\end{multline*}
which defines $\alpha\circ\beta$ (after taking $(A\protens -)\rtimes \Z_2\alcross T$). We have proved that the composition equals $f\mapsto f(D_1)$. But $D_1$ is a surjective operator with 1-dimensional kernel spanned by the 0-form $\exp{(-\frac{1}{2}||x||^2)}$. So $f\mapsto f(D_1)$ is smoothly and $\Z_2$-equivariantly homotopic to $f\mapsto f(0)p_{\Ker D_1}$ (by making $f$'s support smaller and smaller, using that the spectrum is discrete, in other words using the formula $f(tD_1)$). So the diagram
$$
\xymatrix{ \S\alcross\Z_2\ar[rd]_{\ev_0\alcross\Z_2}\ar[r]&\CCINF\R^4\otimes M_4\ar[r]&\smKo\otimes M_4\alcross \Z_2\\
&\C\alcross\Z_2\ar[ru]_{\ \ \ \ \ (z\mapsto p_{\Ker D_1})\alcross\Z_2}
}
$$
commutes up to smooth homotopy, and after comparing with the projection onto one summand produced by $M_4\alcross\Z_2$, we obtain the commutativity of the following diagram.
$$
\xymatrix{
A\protens\S\alcross\Z_2\alcross
T\ar[rd]\ar[r]^{\alpha\circ\beta}&(A\alcross T)^{\oplus
  2}\ar[r]&A\alcross T\\
&A\otimes\C\alcross\Z_2\alcross T\ar@{=}[u]\ar[ur]
}
$$
which implies that $\alpha\circ \beta$ is an $\HA$-equivalence, which is the claim.
\thispagestyle{Selbergemacht5}
\end{proof}
\end{appendix}
\pagestyle{Selbergemacht4}
\clearpage{\pagestyle{empty}\cleardoublepage}

\newpage
\thispagestyle{empty}
\vspace{7cm}
\fontfamily{cmr}\fontseries{bx}
L'homme est n\'e libre, et partout il est dans les fers.\\
\vspace{2em}
\it{Jean-Jacques Rousseau}
\end{document}